\documentclass[letterpaper]{amsart}

\pdfoutput=1

\theoremstyle{definition}

\theoremstyle{remark}

\usepackage{hyperref}
\usepackage{graphicx}
\usepackage{color}
\usepackage{amsmath}
\usepackage{amsfonts}
\usepackage{euler}
\usepackage{amssymb}
\usepackage{epsfig}
\usepackage{rotating}
\usepackage{graphics}

\definecolor{blue}{cmyk}{1.,1.,0.,0.63}
\definecolor{red}{cmyk}{0.,1.,1.,0.63}
\definecolor{green}{cmyk}{1.,0.,1.,0.63}
\definecolor{black}{cmyk}{1.,1.,1.,1.}

\usepackage{t1enc}

\usepackage{ifthen}

\numberwithin{equation}{section}

\newcounter{romenumi}
\newcommand{\labelromenumi}{(\roman{romenumi})}

\renewcommand{\Im}{{\rm Im}}

\renewcommand{\Re}{{\rm Re}}

\newcommand{\isqrt}{{\scriptstyle{\sqrt{-1}}}}

\newcommand{\oaux}{{\text{\usefont{T1}{qcs}{m}{sl}o}}}

\begin{document}

\title[Infinitesimal CR Symmetries of Accidental CR Structures]{
Infinitesimal CR Symmetries 
\\
of Accidental CR Structures}

\vskip 1.truecm

\author{C. Denson Hill} 
\address{Department of Mathematics, Stony Brook University, 
Stony Brook NY 11794, USA}
\email{Dhill@math.stonybrook.edu}

\author{Jo\"el Merker}
\address{Laboratoire de Math\'ematiques d'Orsay, 
Universit\'e Paris-Saclay, Facult\'e des Sciences, 
91405 Orsay Cedex, France}
\email{joel.merker@universite-paris-saclay.fr}

\author{Zhaohu Nie}
\address{Department of Mathematics and Statistics, 
Utah State University, Logan, UT 84322-3900, USA}
\email{zhaohu.nie@usu.edu}

\author{Pawe\l~ Nurowski} 
\address{Center for Theoretical Physics, Polish Academy of Sciences, 
Al. Lotnik\'ow 32/46, 02-668 Warszawa, Poland}
\email{Pawel.Nurowski@fuw.edu.pl}

\date{\today}

\begin{abstract}
In this companion paper to our article {\em Accidental 
CR structures} (arxiv.org, January 2023),
thought of as an appendix not submitted for publication,
we provide complete explicit lists of infinitesimal
CR automorphisms for the concerned CR models 
having respective Lie algebra structures:
\[
{\bf E}_{II},
\ \ \ \ \ \ \ \ \ \ \ \ \ \ \ \ \ \ \ \
{\bf E}_{III},
\ \ \ \ \ \ \ \ \ \ \ \ \ \ \ \ \ \ \ \
\mathfrak{so}(\ell-1,\ell+1),
\ \ \ \ \ \ \ \ \ \ \ \ \ \ \ \ \ \ \ \
\mathfrak{su}(p,q).
\]

We start from our lists of {\em quadric}
CR submanifolds $M^{2n+c} \subset \mathbb{C}^{n+c}$
of codimension $c >1$ 
which are shown to be {\em accidental},
in the sense that their CR symmetry groups are {\em equal to}
(and not smaller than) 
the symmetry groups of the underlying
real distribution structures\,\,---\,\,after forgetting 
the complex structure.

Thanks to intensive symbolic computer explorations,  
we then determine embedded 
vector field generators of these CR symmetries Lie algebras,
and we express them in {\em extrinsic} holomorphic coordinates, 
because intrinsic formulas would be too extended to be shown.
\end{abstract}

\maketitle

\section{Introduction}
\label{introduction}

Any smooth manifold $M$ can be locally equipped with a vector
subdistribution (subbundle) 
$H \subset TM$ of any rank and of any codimension.
Such a pair $(M, H)$ is called {\sl almost Cauchy-Riemann}
(almost CR for short) if $H$ can be `{\sl decorated}'
{\cite{Hill-Merker-Nie-Nurowski-2023}} by a linear
operator $J \colon H \longrightarrow H$
satisfying $J^2 = -\, {\rm Id}$. 
If $J$ is {\sl integrable} in the sense that
for any two sections $X, Y$ of $H$:
\[
\big([X,Y]-[JX,JY]\big)
\,\in\,H,
\ \ \ \ \ \ \ \ \ \ \ \ \ \ \ \ \ \ \ \
J\big([X,Y]-[JX,JY]\big)=[JX,Y]+[X,JY],
\]
then $J$ is said to be {\sl integrable}, and 
$(M, H, J)$ is said to be a {\sl CR manifold}.
One also calls $H =: T^cM$ the {\sl complex tangent bundle},
necessarily of even rank.
We will assume positivity of CR dimension and of codimension:
\[
1
\,\leqslant\,
{\rm CRdim}\,M
\,:=\,
\tfrac{1}{2}\,
{\rm rank}\,T^cM,
\ \ \ \ \ \ \ \ \ \ \ \ \ \ \ \ \ \ \ \
\dim\,M-2\,n
\,=:\,
c
\,\geqslant\,
1.
\]

Both for a naked distribution structure $(M, H)$ and
for a CR structure $(M, T^cM, J)$, 
the two groups $G$ and $G_J$ of 
automorphisms and their respective Lie algebras
$\mathfrak{g}$ and $\mathfrak{g}_J$ 
of infinitesimal symmetries,
are presented and defined precisely 
in~{\cite{Hill-Merker-Nie-Nurowski-2023}}.
Of course:
\[
G_J
\subset
G,
\ \ \ \ \ \ \ \ \ \ \ \ \ \ \ \ \ \ \ \
\mathfrak{g}_J
\,\subset\,
\mathfrak{g},
\]
the typical situation being that $G_J \subsetneqq G$ and
$\mathfrak{g}_J \subsetneqq \mathfrak{g}$ 
are {\em proper inclusions}.

But (very) exceptionally,
equalities $G_J = G$ and 
$\mathfrak{g}_J = \mathfrak{g}$ can occur.
Therefore, 
in~{\cite{Hill-Merker-Nie-Nurowski-2023}}, such CR
structures $(M, T^cM, J)$ for which $G_J = G$ and 
$\mathfrak{g}_J = \mathfrak{g}$
are called
{\sl accidental}, in the sense that
the generic situation that $G_J \subsetneqq G$
is a {\em proper} subgroup is
"by accident" replaced with the "exceptional" equality $G_J = G$. 

Applying purely algebraic Lie-theoretic methods,
Medori-Nacinovich~\cite{Medori-Nacinovich-1998}
classified abstractly such accidental 
CR structures with $\mathfrak{g}$ {\em simple}.
We would like to mention that
not all embeddings of such CR structures 
into some complex space $\mathbb{C}^N$
are known,
especially none is known in codimension $c>1$. 

In this paper, which consists of {\em advanced extended formulas},
our main aim is to 
uncover {\em explicit algebraic faces of accidental 
CR structures},
by exhibiting complete lists of infinitesimal CR symmetries.

Under the form of a long `appendix' not submitted for publication
in a journal,
we therefore complement our 
article~{\cite{Hill-Merker-Nie-Nurowski-2023}},
which provided
a full list {\em embedded CR structures}
in appropriate $\mathbb{C}^N$'s satisfying:

\smallskip\noindent$\bullet$\,
they have real codimension $c>1$;

\smallskip\noindent$\bullet$\,
the complex-tangential distribution $T^cM$ 
is such that $[T^cM, T^cM] + T^cM = TM$;

\smallskip\noindent$\bullet$\,
the local group $G_J$ of CR automorphisms 
of the CR structure $(M,T^cM,J)$ is simple, 
acts transitively on $M$,
and has isotropy $P$ being a parabolic subgroup in $G$;

\smallskip\noindent$\bullet$\,
the local symmetry group $G$ of the naked real vector distribution 
$T^cM$ on $M$ coincides with the group $G_J = G$ 
of CR automorphisms of $(M,T^cM,J)$, or, what is 
the same, $\mathfrak{g}_J = \mathfrak{g}$.

\medskip

In particular, the article~{\cite{Hill-Merker-Nie-Nurowski-2023}}
provided {\em striking 
geometric realizations}
of real simple Lie groups, 
such as for example 
two real forms ${\bf E_{II}}$ 
and ${\bf E_{III}}$ of the
exceptional simple complex Lie group ${\bf E_6}$. 
These two
real simple exceptional Lie groups of dimension 78 are realized
as {\em groups of CR automorphisms of the corresponding CR
structures}, of CR dimension 8 and of codimension 8.

Thus, in the present paper `appended' 
to~{\cite{Hill-Merker-Nie-Nurowski-2023}}, 
we show in particular two collections of 78 
{\em explicit} holomorphic vector fields
in $\mathbb{C}^{8 + 8}$ which generate the respective
Lie algebras of infinitesimal CR symmetries 
of these two embedded CR realizations of ${\bf E_{II}}$ 
and of ${\bf E_{III}}$. We also show much more.

\section{Background about Embedded CR Manifolds}
\label{background-embedded-CR-manifolds}

It is known~{\cite{Merker-Porten-2006}}
that every transitive group action can be equipped
with some {\em real analytic} manifold structure.
It is also known that every real analytic (abstract) CR structure
$(M, T^cM, J)$ of CR dimension $n \geqslant 1$ and
of codimension $c \geqslant 1$
can be {\em embedded} as a CR-generic submanifold:
\[
M^{2n+c}
\,\subset\,
\mathbb{C}^{n+c},
\]
{\sl genericity} meaning that:
\[
TM
+
JTM
\,=\,
T\mathbb{C}^{n+c}
\big\vert_M.
\]
Clearly:
\[
T^cM
\,=\,
TM
\cap
JTM.
\]

It is also known that there exist complex coordinates: 
\[
\big(
z_1,\dots,z_n,\,w_1,\dots,w_c
\big),
\ \ \ \ \ \ \ \ 
z_i
\,=\,
x_i+\isqrt\,y_i,
\ \ \ \ \ \ \ \ 
w_j
\,=\,
u_j+\isqrt\,v_j,
\]
in which such a CR-generic submanifold
$M^{2n+c} \subset \mathbb{C}^{n+c}$ is 
locally graphed as:
\[
2\,u_j
\,=\,
\varphi_j\big(z,\overline{z},v\big)
\eqno
{\scriptstyle{(1\,\leqslant\,j\,\leqslant\,c)}},
\]
with the $\varphi_j$'s being real analytic in some neighborhood of 
the origin:
\[
\varphi_j
\,=\,
\sum_{\alpha\in\mathbb{N}^n}\,
\sum_{\beta\in\mathbb{N}^n}\,
\sum_{\tau\in\mathbb{N}^c}\,
\varphi_{j\alpha\beta\tau}\,
z^\alpha
\overline{z}^\beta
v^\tau,
\ \ \ \ \ \ \ \ \ \ \ \ \ \ \ \ \ \ \ \
\big\vert
\varphi_{j\alpha\beta\tau}
\big\vert
\,\leqslant\,
C^{\vert\alpha\vert+\vert\beta\vert+\vert\tau\vert}.
\]

When the CR structure is embedded,
the infinitesimal CR symmetries of $M$
take the form
of {\em holomorphic} vector fields:
\[
L
\,=\,
\sum_{i=1}^n\,
Z_i(z,w)\,
\partial_{z_i}
+
\sum_{j=1}^c\,
W_j(z,w)\,
\partial_{w_j},
\]
whose (double) real part is {\em tangent} 
to the (extrinsic) CR-generic submanifold:
\[
\aligned
\mathfrak{hol}(M)
&
\,:=\,
\big\{
L\,\colon\,\,
\big(L+\overline{L}\big)\big\vert_M\,\,
\text{is tangent to}\,\,
M
\big\},
\\
\mathfrak{aut}_{CR}(M)
&
\,:=\,
2\,\Re\,\mathfrak{hol}(M).
\endaligned
\]
This makes two uniquely defined {\em real} Lie algebras
of vector fields which are {\em isomorphic}. Usually, 
only the holomorphic 
$L$ is written, with the understanding that it comes
with its `$+\overline{L}$', 
as it is the case for {\em all} of the explicit formulas
exhibited {\em infra}.

For all accidental CR models
of~{\cite{Hill-Merker-Nie-Nurowski-2023}} 
which have {\em quadratic} (see also below)
right-hand sides $\varphi_j = Q_j(z, \overline{z})$ 
{\em not depending} on $v$:
\[
w_j+\overline{w}_j
\,=\,
Q_j(z,\overline{z}),
\]
it is natural to attribute the weights:
\[
[z^i]:=1=:[\overline{z}^i],
\ \ \ \ \
\big[\partial_{z^i}\big]
:=-1=:
\big[\partial_{\overline{z}^i}\big],
\ \ \ \ \ 
[w^j]:=2=:[\overline{w}^j],
\ \ \ \ \ 
\big[\partial_{w^j}\big]
:=-2=:
\big[\partial_{\overline{w}^j}\big],
\]
and it is easy to verify that $\mathfrak{hol}(M)$ then splits into
a direct sum:
\[
\mathfrak{hol}(M)
\,=\,
\bigoplus_{\oaux\in\mathbb{Z}}\,
\mathfrak{g}_{\oaux},
\]
with the graded pieces:
\[
\mathfrak{g}_{\oaux}
\,=\,
\big\{
L\in\mathfrak{hol}(M)
\colon\,\,
[L]
=
\oaux
\big\},
\]
so that the elements of $\mathfrak{g}_{\oaux}$ are those vector fields
having only terms of constant weight $\oaux$:
\[
L_{\oaux}
\,:=\,
\sum_{i=1}^n\,
\bigg(
\sum_{\vert\alpha\vert+\vert\tau\vert=\oaux+1}\,
Z_{j\alpha\tau}\,
z^\alpha\,
w^\tau
\bigg)
\frac{\partial}{\partial{z_i}}
+
\sum_{j=1}^c\,
\bigg(
\sum_{\vert\alpha\vert+\vert\tau\vert=\oaux+2}\,
W_{j\alpha\tau}\,
z^\alpha\,
w^\tau
\bigg)
\frac{\partial}{\partial{w_j}},
\]
which satisfy that $L_{\oaux} + \overline{L}_{\oaux}$
is tangent to $M$.

In fact, for all accidental CR models 
of~{\cite{Hill-Merker-Nie-Nurowski-2023}},
all $\mathfrak{g}_{\oaux} = 0$ are zero for $\vert \oaux \vert
\geqslant 3$, whence:
\[
\mathfrak{hol}(M)
\,=\,
\mathfrak{g}_{-2}
\oplus
\mathfrak{g}_{-1}
\oplus
\mathfrak{g}_{0}
\oplus
\mathfrak{g}_{1}
\oplus
\mathfrak{g}_{2}
\ \ \ \ \ \ \ \ \ 
\text{with}
\ \ \ \ \ \ \ \ \
\dim\,\mathfrak{g}_{-\oaux}
\,=\,
\dim\,\mathfrak{g}_{\oaux}.
\]

Therefore, for the CR-generic models 
of~{\cite{Hill-Merker-Nie-Nurowski-2023}},
having simple accidental CR geometry,
we determined the full graded pieces
$\mathfrak{g}_{\oaux}$ with $-2 \leqslant \oaux \leqslant 2$ 
of $\mathfrak{hol}(M)$. 

\smallskip

{\em About 454 exploration files on Maple were cooked up}.

\smallskip

{\em Welcome!}

\section{$E_{II}$ CR Models}
\label{E-II-CR-models}

In $\mathbb{C}^{8 + 8}$, consider:
\[
\aligned
\Re\,w_1
&
\,=\,
\Re\,
\big(
z_1\,\overline{z}_4
+
z_2\,\overline{z}_3
\big),
\\
\Re\,w_2
&
\,=\,
\Re\,
\big(
z_1\,\overline{z}_6
+
z_2\,\overline{z}_5
\big),
\\
\Im\,w_3
&
\,=\,
\Im\,
\big(
z_1\,\overline{z}_7
+
z_5\,\overline{z}_3
\big),
\\
\Im\,w_4
&
\,=\,
\Im\,
\big(
z_2\,\overline{z}_7
+
z_3\,\overline{z}_6
-
z_5\,\overline{z}_4
-
z_8\,\overline{z}_1
\big),
\\
\Re\,w_5
&
\,=\,
\Re\,
\big(
z_2\,\overline{z}_7
+
z_3\,\overline{z}_6
-
z_5\,\overline{z}_4
-
z_8\,\overline{z}_1
\big),
\\
\Im\,w_6
&
\,=\,
\Im\,
\big(
z_2\,\overline{z}_8
+
z_6\,\overline{z}_4
\big),
\\
\Re\,w_7
&
\,=\,
\Re\,
\big(
z_3\,\overline{z}_8
+
z_4\,\overline{z}_7
\big),
\\
\Re\,w_8
&
\,=\,
\Re\,
\big(
z_5\,\overline{z}_8
+
z_6\,\overline{z}_7
\big).
\endaligned
\]

The $3 + 5$ generators of $\mathfrak{g}_{-2}$ are:
\[
\aligned
L_{w_3}^{-2}
&
\,:=\,
\partial_{w_3},
\\
L_{w_4}^{-2}
&
\,:=\,
\partial_{w_4},
\\
L_{w_6}^{-2}
&
\,:=\,
\partial_{w_6},
\endaligned
\ \ \ \ \ \ \ \ \ \ \ \ \ \ \ \ \ \ \ \ \ \ \ \ \ \
\aligned
IL_{w_1}^{-2}
&
\,:=\,
\isqrt\,\partial_{w_1},
\\
IL_{w_2}^{-2}
&
\,:=\,
\isqrt\,\partial_{w_2},
\\
IL_{w_5}^{-2}
&
\,:=\,
\isqrt\,\partial_{w_5},
\\
IL_{w_7}^{-2}
&
\,:=\,
\isqrt\,\partial_{w_7},
\\
IL_{w_8}^{-2}
&
\,:=\,
\isqrt\,\partial_{w_8}.
\endaligned
\]

The $8 + 8$ generators of $\mathfrak{g}_{-1}$ are:
\[
\aligned
L_{z_1}^{-1}
&
\,:=\,
\partial_{z_1}
+
z_4\,\partial_{w_1}
+
z_6\,\partial_{w_2}
-
z_7\,\partial_{w_3}
-
z_8\,\partial_{w_4}
-
z_8\,\partial_{w_5},
\\
L_{z_2}^{-1}
&
\,:=\,
\partial_{z_2}
+
z_3\,\partial_{w_1}
+
z_5\,\partial_{w_2}
-
z_7\,\partial_{w_4}
+
z_7\,\partial_{w_5}
-
z_8\,\partial_{w_6},
\\
L_{z_3}^{-1}
&
\,:=\,
\partial_{z_3}
+
z_2\,\partial_{w_1}
+
z_5\,\partial_{w_3}
-
z_6\,\partial_{w_4}
+
z_6\,\partial_{w_5}
+
z_8\,\partial_{w_7},
\\
L_{z_4}^{-1}
&
\,:=\,
\partial_{z_4}
+
z_1\,\partial_{w_1}
-
z_5\,\partial_{w_4}
-
z_5\,\partial_{w_5}
+
z_6\,\partial_{w_6}
+
z_7\,\partial_{w_7},
\\
L_{z_5}^{-1}
&
\,:=\,
\partial_{z_5}
+
z_2\,\partial_{w_2}
-
z_3\,\partial_{w_3}
+
z_4\,\partial_{w_4}
-
z_4\,\partial_{w_5}
+
z_8\,\partial_{w_8},
\\
L_{z_6}^{-1}
&
\,:=\,
\partial_{z_6}
+
z_1\,\partial_{w_2}
+
z_3\,\partial_{w_4}
+
z_3\,\partial_{w_5}
-
z_4\,\partial_{w_6}
+
z_7\,\partial_{w_8},
\\
L_{z_7}^{-1}
&
\,:=\,
\partial_{z_7}
+
z_1\,\partial_{w_3}
+
z_2\,\partial_{w_4}
+
z_2\,\partial_{w_5}
+
z_4\,\partial_{w_7}
+
z_6\,\partial_{w_8},
\\
L_{z_8}^{-1}
&
\,:=\,
\partial_{z_8}
+
z_1\,\partial_{w_4}
-
z_1\,\partial_{w_5}
+
z_2\,\partial_{w_6}
+
z_3\,\partial_{w_7}
+
z_5\,\partial_{w_8},
\endaligned
\]
\[
\aligned
IL_{z_1}^{-1}
&
\,:=\,
\isqrt\,
\Big(
\partial_{z_1}
-
z_4\,\partial_{w_1}
-
z_6\,\partial_{w_2}
+
z_7\,\partial_{w_3}
+
z_8\,\partial_{w_4}
+
z_8\,\partial_{w_5}
\Big),
\\
IL_{z_2}^{-1}
&
\,:=\,
\isqrt\,
\Big(
\partial_{z_2}
-
z_3\,\partial_{w_1}
-
z_5\,\partial_{w_2}
+
z_7\,\partial_{w_4}
-
z_7\,\partial_{w_5}
+
z_8\,\partial_{w_6}
\Big),
\\
IL_{z_3}^{-1}
&
\,:=\,
\isqrt\,
\Big(
\partial_{z_3}
-
z_2\,\partial_{w_1}
-
z_5\,\partial_{w_3}
+
z_6\,\partial_{w_4}
-
z_6\,\partial_{w_5}
-
z_8\,\partial_{w_7}
\Big),
\\
IL_{z_4}^{-1}
&
\,:=\,
\isqrt\,
\Big(
\partial_{z_4}
-
z_1\,\partial_{w_1}
+
z_5\,\partial_{w_4}
+
z_5\,\partial_{w_5}
-
z_6\,\partial_{w_6}
-
z_7\,\partial_{w_7}
\Big),
\\
IL_{z_5}^{-1}
&
\,:=\,
\isqrt\,
\Big(
\partial_{z_5}
-
z_2\,\partial_{w_2}
+
z_3\,\partial_{w_3}
-
z_4\,\partial_{w_4}
+
z_4\,\partial_{w_5}
-
z_8\,\partial_{w_8}
\Big),
\\
IL_{z_6}^{-1}
&
\,:=\,
\isqrt\,
\Big(
\partial_{z_6}
-
z_1\,\partial_{w_2}
-
z_3\,\partial_{w_4}
-
z_3\,\partial_{w_5}
+
z_4\,\partial_{w_6}
-
z_7\,\partial_{w_8}
\Big),
\\
IL_{z_7}^{-1}
&
\,:=\,
\isqrt\,
\Big(
\partial_{z_7}
-
z_1\,\partial_{w_3}
-
z_2\,\partial_{w_4}
-
z_2\,\partial_{w_5}
-
z_4\,\partial_{w_7}
-
z_6\,\partial_{w_8}
\Big),
\\
IL_{z_8}^{-1}
&
\,:=\,
\isqrt\,
\Big(
\partial_{z_8}
-
z_1\,\partial_{w_4}
+
z_1\,\partial_{w_5}
-
z_2\,\partial_{w_6}
-
z_3\,\partial_{w_7}
-
z_5\,\partial_{w_8}
\Big).
\endaligned
\]

The $14 + 16$ generators of $\mathfrak{g}_0$ are:
\[
\aligned
L_1^0
&
\,:=\,
-z_5\,\partial_{z_1}+z_6\,\partial_{z_2}+z_7\,\partial_{z_3}-z_8\,\partial_{z_4}+w_5\,\partial_{w_1}+2\,w_8\,\partial_{w_5},
\\
L_2^0
&
\,:=\,
z_3\,\partial_{z_1}-z_4\,\partial_{z_2}+z_7\,\partial_{z_5}-z_8\,\partial_{z_6}+w_5\,\partial_{w_2}-2\,w_7\,\partial_{w_5},
\\
L_3^0
&
\,:=\,
z_2\,\partial_{z_1}-z_4\,\partial_{z_3}-z_6\,\partial_{z_5}+z_8\,\partial_{z_7}+w_4\,\partial_{w_3}+2\,w_6\,\partial_{w_4},
\\
L_4^0
&
\,:=\,
z_2\,\partial_{z_2}+z_4\,\partial_{z_4}+z_6\,\partial_{z_6}+z_8\,\partial_{z_8}+w_1\,\partial_{w_1}+w_2\,\partial_{w_2}+w_4\,\partial_{w_4}
\\
&
\ \ \ \ \ 
+w_5\,\partial_{w_5}+2\,w_6\,\partial_{w_6}+w_7\,\partial_{w_7}+w_8\,\partial_{w_8},
\\
L_5^0
&
\,:=\,
z_1\,\partial_{z_2}-z_3\,\partial_{z_4}-z_5\,\partial_{z_6}+z_7\,\partial_{z_8}+2\,w_3\,\partial_{w_4}+w_4\,\partial_{w_6},
\\
L_6^0
&
\,:=\,
-z_1\,\partial_{z_3}+z_2\,\partial_{z_4}-z_5\,\partial_{z_7}+z_6\,\partial_{z_8}-2\,w_2\,\partial_{w_5}+w_5\,\partial_{w_7},
\\
L_7^0
&
\,:=\,
-z_1\,\partial_{z_5}+z_2\,\partial_{z_6}+z_3\,\partial_{z_7}-z_4\,\partial_{z_8}+2\,w_1\,\partial_{w_5}+w_5\,\partial_{w_8},
\\
L_8^0
&
\,:=\,
z_1\,\partial_{z_1}+z_2\,\partial_{z_2}-z_7\,\partial_{z_7}-z_8\,\partial_{z_8}+w_1\,\partial_{w_1}+w_2\,\partial_{w_2}-w_7\,\partial_{w_7}-w_8\,\partial_{w_8},
\\
L_9^0
&
\,:=\,
z_7\,\partial_{z_1}+z_8\,\partial_{z_2}+w_7\,\partial_{w_1}+w_8\,\partial_{w_2},
\\
L_{10}^0
&
\,:=\,
-z_2\,\partial_{z_2}+z_3\,\partial_{z_3}-z_6\,\partial_{z_6}+z_7\,\partial_{z_7}-w_2\,\partial_{w_2}+w_3\,\partial_{w_3}-w_6\,\partial_{w_6}+w_7\,\partial_{w_7},
\\
L_{11}^0
&
\,:=\,
z_5\,\partial_{z_3}+z_6\,\partial_{z_4}+w_2\,\partial_{w_1}+w_8\,\partial_{w_7},
\\
L_{12}^0
&
\,:=\,
-z_2\,\partial_{z_2}-z_4\,\partial_{z_4}+z_5\,\partial_{z_5}+z_7\,\partial_{z_7}-w_1\,\partial_{w_1}+w_3\,\partial_{w_3}-w_6\,\partial_{w_6}+w_8\,\partial_{w_8},
\\
L_{13}^0
&
\,:=\,
z_3\,\partial_{z_5}+z_4\,\partial_{z_6}+w_1\,\partial_{w_2}+w_7\,\partial_{w_8},
\\
L_{14}^0
&
\,:=\,
z_1\,\partial_{z_7}+z_2\,\partial_{z_8}+w_1\,\partial_{w_7}+w_2\,\partial_{w_8}.
\endaligned
\]
\[
\aligned
IL_1^0
&
\,:=\,
\isqrt\,
\Big(
z_5\,\partial_{z_1}+z_6\,\partial_{z_2}+z_7\,\partial_{z_3}+z_8\,\partial_{z_4}-w_4\,\partial_{w_1}+2\,w_8\,\partial_{w_4},
\Big),
\\
IL_2^0
&
\,:=\,
\isqrt\,
\Big(
z_3\,\partial_{z_1}+z_4\,\partial_{z_2}-z_7\,\partial_{z_5}-z_8\,\partial_{z_6}+w_4\,\partial_{w_2}+2\,w_7\,\partial_{w_4},
\Big),
\\
IL_3^0
&
\,:=\,
\isqrt\,
\Big(
z_1\,\partial_{z_3}+z_2\,\partial_{z_4}-z_5\,\partial_{z_7}-z_6\,\partial_{z_8}+2\,w_2\,\partial_{w_4}+w_4\,\partial_{w_7},
\Big),
\\
IL_4^0
&
\,:=\,
\isqrt\,
\Big(
z_1\,\partial_{z_5}+z_2\,\partial_{z_6}+z_3\,\partial_{z_7}+z_4\,\partial_{z_8}-2\,w_1\,\partial_{w_4}+w_4\,\partial_{w_8},
\Big),
\\
IL_5^0
&
\,:=\,
\isqrt\,
\Big(
z_6\,\partial_{z_1}-z_8\,\partial_{z_3}+w_6\,\partial_{w_1}+w_8\,\partial_{w_3},
\Big),
\\
IL_6^0
&
\,:=\,
\isqrt\,
\Big(
z_5\,\partial_{z_2}-z_7\,\partial_{z_4}+w_3\,\partial_{w_1}+w_8\,\partial_{w_6},
\Big),
\\
IL_7^0
&
\,:=\,
\isqrt\,
\Big(
z_2\,\partial_{z_2}+z_3\,\partial_{z_3}+z_5\,\partial_{z_5}+z_8\,\partial_{z_8}+w_5\,\partial_{w_4}+w_4\,\partial_{w_5},
\Big),
\\
IL_8^0
&
\,:=\,
\isqrt\,
\Big(
z_4\,\partial_{z_1}+z_8\,\partial_{z_5}-w_6\,\partial_{w_2}+w_7\,\partial_{w_3},
\Big),
\\
IL_9^0
&
\,:=\,
\isqrt\,
\Big(
z_1\,\partial_{z_1}+z_4\,\partial_{z_4}+z_6\,\partial_{z_6}+z_7\,\partial_{z_7}-w_5\,\partial_{w_4}-w_4\,\partial_{w_5},
\Big),
\\
IL_{10}^0
&
\,:=\,
\isqrt\,
\Big(
z_3\,\partial_{z_2}+z_7\,\partial_{z_6}-w_3\,\partial_{w_2}+w_7\,\partial_{w_6},
\Big),
\\
IL_{11}^0
&
\,:=\,
\isqrt\,
\Big(
z_2\,\partial_{z_5}-z_4\,\partial_{z_7}+w_1\,\partial_{w_3}+w_6\,\partial_{w_8},
\Big),
\\
IL_{12}^0
&
\,:=\,
\isqrt\,
\Big(
z_2\,\partial_{z_3}+z_6\,\partial_{z_7}-w_2\,\partial_{w_3}+w_6\,\partial_{w_7},
\Big),
\\
IL_{13}^0
&
\,:=\,
\isqrt\,
\Big(
z_2\,\partial_{z_1}+z_4\,\partial_{z_3}+z_6\,\partial_{z_5}+z_8\,\partial_{z_7}+w_5\,\partial_{w_3}-2\,w_6\,\partial_{w_5},
\Big),
\\
IL_{14}^0
&
\,:=\,
\isqrt\,
\Big(
z_1\,\partial_{z_6}-z_3\,\partial_{z_8}+w_1\,\partial_{w_6}+w_3\,\partial_{w_8},
\Big),
\\
IL_{15}^0
&
\,:=\,
\isqrt\,
\Big(
z_1\,\partial_{z_4}+z_5\,\partial_{z_8}-w_2\,\partial_{w_6}+w_3\,\partial_{w_7},
\Big),
\\
IL_{16}^0
&
\,:=\,
\isqrt\,
\Big(
z_1\,\partial_{z_2}+z_3\,\partial_{z_4}+z_5\,\partial_{z_6}+z_7\,\partial_{z_8}+2\,w_3\,\partial_{w_5}-w_5\,\partial_{w_6}.
\Big).
\endaligned
\]

The $8 + 8$ generators of $\mathfrak{g}_1$ are:
\[
\aligned
L_{z_1z_1}^1
&
\,:=\,
2\,z_1^2\,\partial_{z_1}+2\,z_1\,z_2\,\partial_{z_2}+2\,z_1\,z_3\,\partial_{z_3}
-(-2\,z_1\,z_4+2\,w_1)\,\partial_{z_4}+2\,z_1\,z_5\,\partial_{z_5}
\\
&
\ \ \ \ \
-(-2\,z_1\,z_6+2\,w_2)\,\partial_{z_6}+(2\,z_1\,z_7+2\,w_3)\,\partial_{z_7}+(2\,z_2\,z_7-2\,z_3\,z_6+2\,z_4\,z_5+w_4+w_5)\,\partial_{z_8}
\\
&
\ \ \ \ \
+2\,w_1\,z_1\,\partial_{w_1}+2\,w_2\,z_1\,\partial_{w_2}+2\,w_3\,z_1\,\partial_{w_3}-(2\,w_1\,z_5-2\,w_2\,z_3-2\,w_3\,z_2-w_4\,z_1+w_5\,z_1)\,\partial_{w_4}
\\
&
\ \ \ \ \
-(2\,w_1\,z_5-2\,w_2\,z_3-2\,w_3\,z_2+w_4\,z_1-w_5\,z_1)\,\partial_{w_5}+(2\,w_1\,z_6-2\,w_2\,z_4+w_4\,z_2-w_5\,z_2)\,\partial_{w_6}
\\
&
\ \ \ \ \
+(2\,w_1\,z_7+2\,w_3\,z_4+w_4\,z_3-w_5\,z_3)\,\partial_{w_7}+(2\,w_2\,z_7+2\,w_3\,z_6+w_4\,z_5-w_5\,z_5)\,\partial_{w_8},
\endaligned
\]
\[
\aligned
L_{z_2z_2}^1
&
\,:=\,
-2\,z_1\,z_2\,\partial_{z_1}-2\,z_2^2\,\partial_{z_2}+(-2\,z_2\,z_3+2\,w_1)\,\partial_{z_3}-2\,z_2\,z_4\,\partial_{z_4}+(-2\,z_2\,z_5+2\,w_2)\,\partial_{z_5}
\\
&
\ \ \ \ \
-2\,z_2\,z_6\,\partial_{z_6}-(2\,z_1\,z_8+2\,z_3\,z_6-2\,z_4\,z_5+w_4-w_5)\,\partial_{z_7}-(2\,z_2\,z_8+2\,w_6)\,\partial_{z_8}
\\
&
\ \ \ \ \
-2\,w_1\,z_2\,\partial_{w_1}-2\,w_2\,z_2\,\partial_{w_2}-(2\,w_1\,z_5-2\,w_2\,z_3+w_4\,z_1+w_5\,z_1)\,\partial_{w_3}
\\
&
\ \ \ \ \
+(2\,w_1\,z_6-2\,w_2\,z_4-w_4\,z_2-w_5\,z_2-2\,w_6\,z_1)\,\partial_{w_4}
\\
&
\ \ \ \ \
-(2\,w_1\,z_6-2\,w_2\,z_4+w_4\,z_2+w_5\,z_2-2\,w_6\,z_1)\,\partial_{w_5}-2\,w_6\,z_2\,\partial_{w_6}
\\
&
\ \ \ \ \
-(2\,w_1\,z_8+w_4\,z_4+w_5\,z_4+2\,w_6\,z_3)\,\partial_{w_7}-(2\,w_2\,z_8+w_4\,z_6+w_5\,z_6+2\,w_6\,z_5)\,\partial_{w_8},
\endaligned
\]
\[
\aligned
L_{z_3z_3}^1
&
\,:=\,
-2\,z_1\,z_3\,\partial_{z_1}+(-2\,z_2\,z_3+2\,w_1)\,\partial_{z_2}-2\,z_3^2\,\partial_{z_3}-2\,z_3\,z_4\,\partial_{z_4}+(-2\,z_3\,z_5+2\,w_3)\,\partial_{z_5}
\\
&
\ \ \ \ \
-(-2\,z_1\,z_8+2\,z_2\,z_7+2\,z_4\,z_5+w_4-w_5)\,\partial_{z_6}-2\,z_3\,z_7\,\partial_{z_7}+(-2\,z_3\,z_8+2\,w_7)\,\partial_{z_8}
\\
&
\ \ \ \ \
-2\,w_1\,z_3\,\partial_{w_1}-(2\,w_1\,z_5-2\,w_3\,z_2+w_4\,z_1+w_5\,z_1)\,\partial_{w_2}-2\,w_3\,z_3\,\partial_{w_3}
\\
&
\ \ \ \ \
+(2\,w_1\,z_7+2\,w_3\,z_4-w_4\,z_3-w_5\,z_3-2\,w_7\,z_1)\,\partial_{w_4}
\\
&
\ \ \ \ \
-(2\,w_1\,z_7+2\,w_3\,z_4+w_4\,z_3+w_5\,z_3-2\,w_7\,z_1)\,\partial_{w_5}
\\
&
\ \ \ \ \
+(2\,w_1\,z_8+w_4\,z_4+w_5\,z_4-2\,w_7\,z_2)\,\partial_{w_6}-2\,w_7\,z_3\,\partial_{w_7}
\\
&
\ \ \ \ \
+(2\,w_3\,z_8-w_4\,z_7-w_5\,z_7-2\,w_7\,z_5)\,\partial_{w_8},
\endaligned
\]
\[
\aligned
L_{z_4z_4}^1
&
\,:=\,
-(-2\,z_1\,z_4+2\,w_1)\,\partial_{z_1}+2\,z_2\,z_4\,\partial_{z_2}+2\,z_3\,z_4\,\partial_{z_3}+2\,z_4^2\,\partial_{z_4}
\\
&
\ \ \ \ \
+(2\,z_1\,z_8-2\,z_2\,z_7+2\,z_3\,z_6+w_4+w_5)\,\partial_{z_5}-(-2\,z_4\,z_6+2\,w_6)\,\partial_{z_6}
\\
&
\ \ \ \ \
-(-2\,z_4\,z_7+2\,w_7)\,\partial_{z_7}+2\,z_4\,z_8\,\partial_{z_8}+2\,w_1\,z_4\,\partial_{w_1}
\\
&
\ \ \ \ \
+(2\,w_1\,z_6+w_4\,z_2-w_5\,z_2-2\,w_6\,z_1)\,\partial_{w_2}-(2\,w_1\,z_7+w_4\,z_3-w_5\,z_3-2\,w_7\,z_1)\,\partial_{w_3}
\\
&
\ \ \ \ \
-(2\,w_1\,z_8-w_4\,z_4+w_5\,z_4+2\,w_6\,z_3-2\,w_7\,z_2)\,\partial_{w_4}
\\
&
\ \ \ \ \
-(2\,w_1\,z_8+w_4\,z_4-w_5\,z_4+2\,w_6\,z_3-2\,w_7\,z_2)\,\partial_{w_5}+2\,w_6\,z_4\,\partial_{w_6}+2\,w_7\,z_4\,\partial_{w_7}
\\
&
\ \ \ \ \
+(w_4\,z_8-w_5\,z_8-2\,w_6\,z_7+2\,w_7\,z_6)\,\partial_{w_8},
\endaligned
\]
\[
\aligned
L_{z_5z_5}^1
&
\,:=\,
2\,z_1\,z_5\,\partial_{z_1}-(-2\,z_2\,z_5+2\,w_2)\,\partial_{z_2}+(2\,z_3\,z_5+2\,w_3)\,\partial_{z_3}
\\
&
\ \ \ \ \
-(-2\,z_1\,z_8+2\,z_2\,z_7-2\,z_3\,z_6+w_4-w_5)\,\partial_{z_4}+2\,z_5^2\,\partial_{z_5}+2\,z_5\,z_6\,\partial_{z_6}+2\,z_5\,z_7\,\partial_{z_7}
\\
&
\ \ \ \ \
-(-2\,z_5\,z_8+2\,w_8)\,\partial_{z_8}+(2\,w_2\,z_3+2\,w_3\,z_2-w_4\,z_1-w_5\,z_1)\,\partial_{w_1}+2\,w_2\,z_5\,\partial_{w_2}
\\
&
\ \ \ \ \
+2\,w_3\,z_5\,\partial_{w_3}-(2\,w_2\,z_7+2\,w_3\,z_6-w_4\,z_5-w_5\,z_5-2\,w_8\,z_1)\,\partial_{w_4}
\\
&
\ \ \ \ \
+(2\,w_2\,z_7+2\,w_3\,z_6+w_4\,z_5+w_5\,z_5-2\,w_8\,z_1)\,\partial_{w_5}-(2\,w_2\,z_8+w_4\,z_6+w_5\,z_6-2\,w_8\,z_2)\,\partial_{w_6}
\\
&
\ \ \ \ \
+(2\,w_3\,z_8-w_4\,z_7-w_5\,z_7+2\,w_8\,z_3)\,\partial_{w_7}+2\,w_8\,z_5\,\partial_{w_8},
\endaligned
\]
\[
\aligned
L_{z_6z_6}^1
&
\,:=\,
(-2\,z_1\,z_6+2\,w_2)\,\partial_{z_1}-2\,z_2\,z_6\,\partial_{z_2}+(2\,z_1\,z_8-2\,z_2\,z_7-2\,z_4\,z_5+w_4+w_5)\,\partial_{z_3}
\\
&
\ \ \ \ \
-(2\,z_4\,z_6+2\,w_6)\,\partial_{z_4}-2\,z_5\,z_6\,\partial_{z_5}-2\,z_6^2\,\partial_{z_6}+(-2\,z_6\,z_7+2\,w_8)\,\partial_{z_7}
\\
&
\ \ \ \ \
-2\,z_6\,z_8\,\partial_{z_8}-(2\,w_2\,z_4-w_4\,z_2+w_5\,z_2+2\,w_6\,z_1)\,\partial_{w_1}-2\,w_2\,z_6\,\partial_{w_2}
\\
&
\ \ \ \ \
+(2\,w_2\,z_7+w_4\,z_5-w_5\,z_5-2\,w_8\,z_1)\,\partial_{w_3}+(2\,w_2\,z_8-w_4\,z_6+w_5\,z_6+2\,w_6\,z_5-2\,w_8\,z_2)\,\partial_{w_4}
\\
&
\ \ \ \ \
+(2\,w_2\,z_8+w_4\,z_6-w_5\,z_6+2\,w_6\,z_5-2\,w_8\,z_2)\,\partial_{w_5}-2\,w_6\,z_6\,\partial_{w_6}
\\
&
\ \ \ \ \
+(w_4\,z_8-w_5\,z_8-2\,w_6\,z_7-2\,w_8\,z_4)\,\partial_{w_7}-2\,w_8\,z_6\,\partial_{w_8},
\endaligned
\]
\[
\aligned
L_{z_7z_7}^1
&
\,:=\,
(-2\,z_1\,z_7+2\,w_3)\,\partial_{z_1}+(-2\,z_1\,z_8-2\,z_3\,z_6+2\,z_4\,z_5+w_4+w_5)\,\partial_{z_2}-2\,z_3\,z_7\,\partial_{z_3}
\\
&
\ \ \ \ \
+(-2\,z_4\,z_7+2\,w_7)\,\partial_{z_4}-2\,z_5\,z_7\,\partial_{z_5}+(-2\,z_6\,z_7+2\,w_8)\,\partial_{z_6}-2\,z_7^2\,\partial_{z_7}-2\,z_7\,z_8\,\partial_{z_8}
\\
&
\ \ \ \ \
+(2\,w_3\,z_4+w_4\,z_3-w_5\,z_3-2\,w_7\,z_1)\,\partial_{w_1}+(2\,w_3\,z_6+w_4\,z_5-w_5\,z_5-2\,w_8\,z_1)\,\partial_{w_2}
\\
&
\ \ \ \ \
-2\,w_3\,z_7\,\partial_{w_3}-(2\,w_3\,z_8+w_4\,z_7-w_5\,z_7-2\,w_7\,z_5+2\,w_8\,z_3)\,\partial_{w_4}
\\
&
\ \ \ \ \
-(2\,w_3\,z_8-w_4\,z_7+w_5\,z_7-2\,w_7\,z_5+2\,w_8\,z_3)\,\partial_{w_5}-(w_4\,z_8-w_5\,z_8+2\,w_7\,z_6-2\,w_8\,z_4)\,\partial_{w_6}
\\
&
\ \ \ \ \
-2\,w_7\,z_7\,\partial_{w_7}-2\,w_8\,z_7\,\partial_{w_8},
\endaligned
\]
\[
\aligned
L_{z_8z_8}^1
&
\,:=\,
-(-2\,z_2\,z_7+2\,z_3\,z_6-2\,z_4\,z_5+w_4-w_5)\,\partial_{z_1}-(-2\,z_2\,z_8+2\,w_6)\,\partial_{z_2}-(-2\,z_3\,z_8+2\,w_7)\,\partial_{z_3}
\\
&
\ \ \ \ \
+2\,z_4\,z_8\,\partial_{z_4}-(-2\,z_5\,z_8+2\,w_8)\,\partial_{z_5}+2\,z_6\,z_8\,\partial_{z_6}+2\,z_7\,z_8\,\partial_{z_7}+2\,z_8^2\,\partial_{z_8}
\\
&
\ \ \ \ \
-(w_4\,z_4+w_5\,z_4+2\,w_6\,z_3-2\,w_7\,z_2)\,\partial_{w_1}-(w_4\,z_6+w_5\,z_6+2\,w_6\,z_5-2\,w_8\,z_2)\,\partial_{w_2}
\\
&
\ \ \ \ \
+(w_4\,z_7+w_5\,z_7+2\,w_7\,z_5-2\,w_8\,z_3)\,\partial_{w_3}+(w_4\,z_8+w_5\,z_8+2\,w_6\,z_7-2\,w_7\,z_6+2\,w_8\,z_4)\,\partial_{w_4}
\\
&
\ \ \ \ \
+(w_4\,z_8+w_5\,z_8-2\,w_6\,z_7+2\,w_7\,z_6-2\,w_8\,z_4)\,\partial_{w_5}+2\,w_6\,z_8\,\partial_{w_6}+2\,w_7\,z_8\,\partial_{w_7}+2\,w_8\,z_8\,\partial_{w_8},
\endaligned
\]
\[
\aligned
IL_{z_1z_1}^1
&
\,:=\,
\isqrt\,
\Big(
2\,z_1^2\,\partial_{z_1}+2\,z_1\,z_2\,\partial_{z_2}+2\,z_1\,z_3\,\partial_{z_3}+(2\,z_1\,z_4+2\,w_1)\,\partial_{z_4}+2\,z_1\,z_5\,\partial_{z_5}+(2\,z_1\,z_6+2\,w_2)\,\partial_{z_6}
\\
&
\ \ \ \ \ \ \ \ \ \ \ \ \
-(-2\,z_1\,z_7+2\,w_3)\,\partial_{z_7}-(-2\,z_2\,z_7+2\,z_3\,z_6-2\,z_4\,z_5+w_4+w_5)\,\partial_{z_8}+2\,w_1\,z_1\,\partial_{w_1}
\\
&
\ \ \ \ \ \ \ \ \ \ \ \ \
+2\,w_2\,z_1\,\partial_{w_2}+2\,w_3\,z_1\,\partial_{w_3}-(2\,w_1\,z_5-2\,w_2\,z_3-2\,w_3\,z_2-w_4\,z_1+w_5\,z_1)\,\partial_{w_4}
\\
&
\ \ \ \ \ \ \ \ \ \ \ \ \
-(2\,w_1\,z_5-2\,w_2\,z_3-2\,w_3\,z_2+w_4\,z_1-w_5\,z_1)\,\partial_{w_5}
\\
&
\ \ \ \ \ \ \ \ \ \ \ \ \
+(2\,w_1\,z_6-2\,w_2\,z_4+w_4\,z_2-w_5\,z_2)\,\partial_{w_6}+(2\,w_1\,z_7+2\,w_3\,z_4+w_4\,z_3-w_5\,z_3)\,\partial_{w_7}
\\
&
\ \ \ \ \ \ \ \ \ \ \ \ \
+(2\,w_2\,z_7+2\,w_3\,z_6+w_4\,z_5-w_5\,z_5)\,\partial_{w_8}
\Big),
\endaligned
\]
\[
\aligned
IL_{z_2z_2}^1
&
\,:=\,
\isqrt\,
\Big(
2\,z_1\,z_2\,\partial_{z_1}+2\,z_2^2\,\partial_{z_2}+(2\,z_2\,z_3+2\,w_1)\,\partial_{z_3}+2\,z_2\,z_4\,\partial_{z_4}+(2\,z_2\,z_5+2\,w_2)\,\partial_{z_5}+2\,z_2\,z_6\,\partial_{z_6}
\\
&
\ \ \ \ \ \ \ \ \ \ \ \ \
-(-2\,z_1\,z_8-2\,z_3\,z_6+2\,z_4\,z_5+w_4-w_5)\,\partial_{z_7}-(-2\,z_2\,z_8+2\,w_6)\,\partial_{z_8}+2\,w_1\,z_2\,\partial_{w_1}
\\
&
\ \ \ \ \ \ \ \ \ \ \ \ \
+2\,w_2\,z_2\,\partial_{w_2}+(2\,w_1\,z_5-2\,w_2\,z_3+w_4\,z_1+w_5\,z_1)\,\partial_{w_3}
\\
&
\ \ \ \ \ \ \ \ \ \ \ \ \
-(2\,w_1\,z_6-2\,w_2\,z_4
-w_4\,z_2-w_5\,z_2-2\,w_6\,z_1)\,\partial_{w_4}
\\
&
\ \ \ \ \ \ \ \ \ \ \ \ \
+(2\,w_1\,z_6-2\,w_2\,z_4+w_4\,z_2+w_5\,z_2-2\,w_6\,z_1)\,\partial_{w_5}+2\,w_6\,z_2\,\partial_{w_6}
\\
&
\ \ \ \ \ \ \ \ \ \ \ \ \
+(2\,w_1\,z_8+w_4\,z_4+w_5\,z_4+2\,w_6\,z_3)\,\partial_{w_7}+(2\,w_2\,z_8+w_4\,z_6+w_5\,z_6+2\,w_6\,z_5)\,\partial_{w_8}
\Big),
\endaligned
\]
\[
\aligned
IL_{z_3z_3}^1
&
\,:=\,
\isqrt\,
\Big(
2\,z_1\,z_3\,\partial_{z_1}+(2\,z_2\,z_3+2\,w_1)\,\partial_{z_2}+2\,z_3^2\,\partial_{z_3}+2\,z_3\,z_4\,\partial_{z_4}+(2\,z_3\,z_5+2\,w_3)\,\partial_{z_5}
\\
&
\ \ \ \ \ \ \ \ \ \ \ \ \
-(2\,z_1\,z_8-2\,z_2\,z_7-2\,z_4\,z_5+w_4-w_5)\,\partial_{z_6}+2\,z_3\,z_7\,\partial_{z_7}+(2\,z_3\,z_8+2\,w_7)\,\partial_{z_8}
\\
&
\ \ \ \ \ \ \ \ \ \ \ \ \
+2\,w_1\,z_3\,\partial_{w_1}+(2\,w_1\,z_5-2\,w_3\,z_2+w_4\,z_1+w_5\,z_1)\,\partial_{w_2}+2\,w_3\,z_3\,\partial_{w_3}
\\
&
\ \ \ \ \ \ \ \ \ \ \ \ \
-(2\,w_1\,z_7+2\,w_3\,z_4-w_4\,z_3-w_5\,z_3-2\,w_7\,z_1)\,\partial_{w_4}
\\
&
\ \ \ \ \ \ \ \ \ \ \ \ \
+(2\,w_1\,z_7+2\,w_3\,z_4+w_4\,z_3+w_5\,z_3-2\,w_7\,z_1)\,\partial_{w_5}
\\
&
\ \ \ \ \ \ \ \ \ \ \ \ \
-(2\,w_1\,z_8+w_4\,z_4+w_5\,z_4-2\,w_7\,z_2)\,\partial_{w_6}+2\,w_7\,z_3\,\partial_{w_7}
\\
&
\ \ \ \ \ \ \ \ \ \ \ \ \
-(2\,w_3\,z_8-w_4\,z_7-w_5\,z_7-2\,w_7\,z_5)\,\partial_{w_8}
\Big),
\endaligned
\]
\[
\aligned
IL_{z_4z_4}^1
&
\,:=\,
\isqrt\,
\Big(
(2\,z_1\,z_4+2\,w_1)\,\partial_{z_1}+2\,z_2\,z_4\,\partial_{z_2}+2\,z_3\,z_4\,\partial_{z_3}+2\,z_4^2\,\partial_{z_4}
\\
&
\ \ \ \ \ \ \ \ \ \ \ \ \
-(-2\,z_1\,z_8+2\,z_2\,z_7-2\,z_3\,z_6+w_4+w_5)\,\partial_{z_5}+(2\,z_4\,z_6+2\,w_6)\,\partial_{z_6}
\\
&
\ \ \ \ \ \ \ \ \ \ \ \ \
+(2\,z_4\,z_7+2\,w_7)\,\partial_{z_7}+2\,z_4\,z_8\,\partial_{z_8}+2\,w_1\,z_4\,\partial_{w_1}
\\
&
\ \ \ \ \ \ \ \ \ \ \ \ \
+(2\,w_1\,z_6+w_4\,z_2-w_5\,z_2-2\,w_6\,z_1)\,\partial_{w_2}-(2\,w_1\,z_7+w_4\,z_3-w_5\,z_3-2\,w_7\,z_1)\,\partial_{w_3}
\\
&
\ \ \ \ \ \ \ \ \ \ \ \ \
-(2\,w_1\,z_8-w_4\,z_4+w_5\,z_4+2\,w_6\,z_3-2\,w_7\,z_2)\,\partial_{w_4}
\\
&
\ \ \ \ \ \ \ \ \ \ \ \ \
-(2\,w_1\,z_8+w_4\,z_4-w_5\,z_4+2\,w_6\,z_3-2\,w_7\,z_2)\,\partial_{w_5}+2\,w_6\,z_4\,\partial_{w_6}+2\,w_7\,z_4\,\partial_{w_7}
\\
&
\ \ \ \ \ \ \ \ \ \ \ \ \
+(w_4\,z_8-w_5\,z_8-2\,w_6\,z_7+2\,w_7\,z_6)\,\partial_{w_8}
\Big),
\endaligned
\]
\[
\aligned
IL_{z_5z_5}^1
&
\,:=\,
\isqrt\,
\Big(
2\,z_1\,z_5\,\partial_{z_1}+(2\,z_2\,z_5+2\,w_2)\,\partial_{z_2}-(-2\,z_3\,z_5+2\,w_3)\,\partial_{z_3}
\\
&
\ \ \ \ \ \ \ \ \ \ \ \ \
+(2\,z_1\,z_8-2\,z_2\,z_7+2\,z_3\,z_6+w_4-w_5)\,\partial_{z_4}+2\,z_5^2\,\partial_{z_5}+2\,z_5\,z_6\,\partial_{z_6}+2\,z_5\,z_7\,\partial_{z_7}
\\
&
\ \ \ \ \ \ \ \ \ \ \ \ \
+(2\,z_5\,z_8+2\,w_8)\,\partial_{z_8}+(2\,w_2\,z_3+2\,w_3\,z_2-w_4\,z_1-w_5\,z_1)\,\partial_{w_1}+2\,w_2\,z_5\,\partial_{w_2}
\\
&
\ \ \ \ \ \ \ \ \ \ \ \ \
+2\,w_3\,z_5\,\partial_{w_3}-(2\,w_2\,z_7+2\,w_3\,z_6-w_4\,z_5-w_5\,z_5-2\,w_8\,z_1)\,\partial_{w_4}
\\
&
\ \ \ \ \ \ \ \ \ \ \ \ \
+(2\,w_2\,z_7+2\,w_3\,z_6+w_4\,z_5+w_5\,z_5-2\,w_8\,z_1)\,\partial_{w_5}
\\
&
\ \ \ \ \ \ \ \ \ \ \ \ \
-(2\,w_2\,z_8+w_4\,z_6+w_5\,z_6-2\,w_8\,z_2)\,\partial_{w_6}
\\
&
\ \ \ \ \ \ \ \ \ \ \ \ \
+(2\,w_3\,z_8-w_4\,z_7-w_5\,z_7+2\,w_8\,z_3)\,\partial_{w_7}+2\,w_8\,z_5\,\partial_{w_8}
\Big),
\endaligned
\]
\[
\aligned
IL_{z_6z_6}^1
&
\,:=\,
\isqrt\,
\Big(
(2\,z_1\,z_6+2\,w_2)\,\partial_{z_1}+2\,z_2\,z_6\,\partial_{z_2}+(-2\,z_1\,z_8+2\,z_2\,z_7+2\,z_4\,z_5+w_4+w_5)\,\partial_{z_3}
\\
&
\ \ \ \ \ \ \ \ \ \ \ \ \
-(-2\,z_4\,z_6+2\,w_6)\,\partial_{z_4}+2\,z_5\,z_6\,\partial_{z_5}+2\,z_6^2\,\partial_{z_6}+(2\,z_6\,z_7+2\,w_8)\,\partial_{z_7}+2\,z_6\,z_8\,\partial_{z_8}
\\
&
\ \ \ \ \ \ \ \ \ \ \ \ \
+(2\,w_2\,z_4-w_4\,z_2+w_5\,z_2+2\,w_6\,z_1)\,\partial_{w_1}+2\,w_2\,z_6\,\partial_{w_2}
\\
&
\ \ \ \ \ \ \ \ \ \ \ \ \
-(2\,w_2\,z_7+w_4\,z_5-w_5\,z_5-2\,w_8\,z_1)\,\partial_{w_3}
\\
&
\ \ \ \ \ \ \ \ \ \ \ \ \
-(2\,w_2\,z_8-w_4\,z_6+w_5\,z_6+2\,w_6\,z_5-2\,w_8\,z_2)\,\partial_{w_4}
\\
&
\ \ \ \ \ \ \ \ \ \ \ \ \
-(2\,w_2\,z_8+w_4\,z_6-w_5\,z_6+2\,w_6\,z_5-2\,w_8\,z_2)\,\partial_{w_5}+2\,w_6\,z_6\,\partial_{w_6}
\\
&
\ \ \ \ \ \ \ \ \ \ \ \ \
-(w_4\,z_8-w_5\,z_8-2\,w_6\,z_7-2\,w_8\,z_4)\,\partial_{w_7}+2\,w_8\,z_6\,\partial_{w_8}
\Big),
\endaligned
\]
\[
\aligned
IL_{z_7z_7}^1
&
\,:=\,
\isqrt\,
\Big(
(2\,z_1\,z_7+2\,w_3)\,\partial_{z_1}+(2\,z_1\,z_8+2\,z_3\,z_6-2\,z_4\,z_5+w_4+w_5)\,\partial_{z_2}+2\,z_3\,z_7\,\partial_{z_3}
\\
&
\ \ \ \ \ \ \ \ \ \ \ \ \
+(2\,z_4\,z_7+2\,w_7)\,\partial_{z_4}+2\,z_5\,z_7\,\partial_{z_5}+(2\,z_6\,z_7+2\,w_8)\,\partial_{z_6}+2\,z_7^2\,\partial_{z_7}+2\,z_7\,z_8\,\partial_{z_8}
\\
&
\ \ \ \ \ \ \ \ \ \ \ \ \
-(2\,w_3\,z_4+w_4\,z_3-w_5\,z_3-2\,w_7\,z_1)\,\partial_{w_1}-(2\,w_3\,z_6+w_4\,z_5-w_5\,z_5-2\,w_8\,z_1)\,\partial_{w_2}
\\
&
\ \ \ \ \ \ \ \ \ \ \ \ \
+2\,w_3\,z_7\,\partial_{w_3}+(2\,w_3\,z_8+w_4\,z_7-w_5\,z_7-2\,w_7\,z_5+2\,w_8\,z_3)\,\partial_{w_4}
\\
&
\ \ \ \ \ \ \ \ \ \ \ \ \
+(2\,w_3\,z_8-w_4\,z_7+w_5\,z_7-2\,w_7\,z_5+2\,w_8\,z_3)\,\partial_{w_5}
\\
&
\ \ \ \ \ \ \ \ \ \ \ \ \
+(w_4\,z_8-w_5\,z_8+2\,w_7\,z_6-2\,w_8\,z_4)\,\partial_{w_6}+2\,w_7\,z_7\,\partial_{w_7}+2\,w_8\,z_7\,\partial_{w_8}
\Big),
\endaligned
\]
\[
\aligned
IL_{z_8z_8}^1
&
\,:=\,
\isqrt\,
\Big(
(2\,z_2\,z_7-2\,z_3\,z_6+2\,z_4\,z_5+w_4-w_5)\,\partial_{z_1}+(2\,z_2\,z_8+2\,w_6)\,\partial_{z_2}+(2\,z_3\,z_8+2\,w_7)\,\partial_{z_3}
\\
&
\ \ \ \ \ \ \ \ \ \ \ \ \
+2\,z_4\,z_8\,\partial_{z_4}+(2\,z_5\,z_8+2\,w_8)\,\partial_{z_5}+2\,z_6\,z_8\,\partial_{z_6}+2\,z_7\,z_8\,\partial_{z_7}+2\,z_8^2\,\partial_{z_8}
\\
&
\ \ \ \ \ \ \ \ \ \ \ \ \
-(w_4\,z_4+w_5\,z_4+2\,w_6\,z_3-2\,w_7\,z_2)\,\partial_{w_1}-(w_4\,z_6+w_5\,z_6+2\,w_6\,z_5-2\,w_8\,z_2)\,\partial_{w_2}
\\
&
\ \ \ \ \ \ \ \ \ \ \ \ \
+(w_4\,z_7+w_5\,z_7+2\,w_7\,z_5-2\,w_8\,z_3)\,\partial_{w_3}
\\
&
\ \ \ \ \ \ \ \ \ \ \ \ \
+(w_4\,z_8+w_5\,z_8+2\,w_6\,z_7-2\,w_7\,z_6+2\,w_8\,z_4)\,\partial_{w_4}
\\
&
\ \ \ \ \ \ \ \ \ \ \ \ \
+(w_4\,z_8+w_5\,z_8-2\,w_6\,z_7+2\,w_7\,z_6-2\,w_8\,z_4)\,\partial_{w_5}+2\,w_6\,z_8\,\partial_{w_6}
\\
&
\ \ \ \ \ \ \ \ \ \ \ \ \
+2\,w_7\,z_8\,\partial_{w_7}+2\,w_8\,z_8\,\partial_{w_8}
\Big).
\endaligned
\]

The $3 + 5$ generators of $\mathfrak{g}_2$ are:
\[
\aligned
L_{w_3w_3}^2
&
\,:=\,
4\,z_1\,w_3\,\partial_{z_1}+(4\,w_1\,z_5-4\,w_2\,z_3+2\,w_4\,z_1+2\,w_5\,z_1)\,\partial_{z_2}+4\,z_3\,w_3\,\partial_{z_3}
\\
&
\ \ \ \ \
-(4\,w_1\,z_7+2\,w_4\,z_3-2\,w_5\,z_3-4\,w_7\,z_1)\,\partial_{z_4}+4\,z_5\,w_3\,\partial_{z_5}-(4\,w_2\,z_7+2\,w_4\,z_5-2\,w_5\,z_5
\\
&
\ \ \ \ \
-4\,w_8\,z_1)\,\partial_{z_6}+4\,z_7\,w_3\,\partial_{z_7}+(2\,w_4\,z_7+2\,w_5\,z_7+4\,w_7\,z_5-4\,w_8\,z_3)\,\partial_{z_8}
\\
&
\ \ \ \ \
+4\,w_1\,w_3\,\partial_{w_1}+4\,w_2\,w_3\,\partial_{w_2}+4\,w_3^2\,\partial_{w_3}+4\,w_3\,w_4\,\partial_{w_4}+4\,w_3\,w_5\,\partial_{w_5}
\\
&
\ \ \ \ \
+(4\,w_1\,w_8-4\,w_2\,w_7+w_4^2-w_5^2)\,\partial_{w_6}+4\,w_3\,w_7\,\partial_{w_7}+4\,w_3\,w_8\,\partial_{w_8},
\endaligned
\]
\[
\aligned
L_{w_4w_4}^2
&
\,:=\,
(2\,w_1\,z_5-2\,w_2\,z_3-2\,w_3\,z_2-w_4\,z_1+w_5\,z_1)\,\partial_{z_1}
\\
&
\ \ \ \ \
+(2\,w_1\,z_6-2\,w_2\,z_4-w_4\,z_2-w_5\,z_2-2\,w_6\,z_1)\,\partial_{z_2}
\\
&
\ \ \ \ \
+(2\,w_1\,z_7+2\,w_3\,z_4-w_4\,z_3-w_5\,z_3-2\,w_7\,z_1)\,\partial_{z_3}
\\
&
\ \ \ \ \
+(2\,w_1\,z_8-w_4\,z_4+w_5\,z_4+2\,w_6\,z_3-2\,w_7\,z_2)\,\partial_{z_4}
\\
&
\ \ \ \ \
+(2\,w_2\,z_7+2\,w_3\,z_6-w_4\,z_5-w_5\,z_5-2\,w_8\,z_1)\,\partial_{z_5}
\\
&
\ \ \ \ \
+(2\,w_2\,z_8-w_4\,z_6+w_5\,z_6+2\,w_6\,z_5-2\,w_8\,z_2)\,\partial_{z_6}
\\
&
\ \ \ \ \
-(2\,w_3\,z_8+w_4\,z_7-w_5\,z_7-2\,w_7\,z_5+2\,w_8\,z_3)\,\partial_{z_7}
\\
&
\ \ \ \ \
-(w_4\,z_8+w_5\,z_8+2\,w_6\,z_7-2\,w_7\,z_6+2\,w_8\,z_4)\,\partial_{z_8}
\\
&
\ \ \ \ \
-2\,w_1\,w_4\,\partial_{w_1}-2\,w_2\,w_4\,\partial_{w_2}-2\,w_3\,w_4\,\partial_{w_3}
\\
&
\ \ \ \ \
+(4\,w_1\,w_8-4\,w_2\,w_7-4\,w_3\,w_6-w_4^2-w_5^2)\,\partial_{w_4}
\\
&
\ \ \ \ \
-2\,w_4\,w_5\,\partial_{w_5}-2\,w_4\,w_6\,\partial_{w_6}
-2\,w_4\,w_7\,\partial_{w_7}-2\,w_4\,w_8\,\partial_{w_8},
\endaligned
\]
\[
\aligned
L_{w_6w_6}^2
&
\,:=\,
(4\,w_1\,z_6-4\,w_2\,z_4+2\,w_4\,z_2-2\,w_5\,z_2)\,\partial_{z_1}+4\,z_2\,w_6\,\partial_{z_2}
\\
&
\ \ \ \ \
-(4\,w_1\,z_8+2\,w_4\,z_4+2\,w_5\,z_4-4\,w_7\,z_2)\,\partial_{z_3}+4\,z_4\,w_6\,\partial_{z_4}
\\
&
\ \ \ \ \
-(4\,w_2\,z_8+2\,w_4\,z_6+2\,w_5\,z_6-4\,w_8\,z_2)\,\partial_{z_5}+4\,z_6\,w_6\,\partial_{z_6}
\\
&
\ \ \ \ \
+(2\,w_4\,z_8-2\,w_5\,z_8+4\,w_7\,z_6-4\,w_8\,z_4)\,\partial_{z_7}+4\,z_8\,w_6\,\partial_{z_8}+4\,w_1\,w_6\,\partial_{w_1}
\\
&
\ \ \ \ \
+4\,w_2\,w_6\,\partial_{w_2}+(4\,w_1\,w_8-4\,w_2\,w_7+w_4^2-w_5^2)\,\partial_{w_3}+4\,w_4\,w_6\,\partial_{w_4}
\\
&
\ \ \ \ \
+4\,w_5\,w_6\,\partial_{w_5}
+4\,w_6^2\,\partial_{w_6}+4\,w_6\,w_7\,\partial_{w_7}+4\,w_6\,w_8\,\partial_{w_8},
\endaligned
\]

\[
\aligned
IL_{w_1w_1}^2
&
\,:=\,
\isqrt\,
\Big(
4\,z_1\,w_1\,\partial_{z_1}+4\,z_2\,w_1\,\partial_{z_2}+4\,z_3\,w_1\,\partial_{z_3}+4\,z_4\,w_1\,\partial_{z_4}
\\
&
\ \ \ \ \ \ \ \ \ \ \ \ \
+(4\,w_2\,z_3+4\,w_3\,z_2-2\,w_4\,z_1-2\,w_5\,z_1)\,\partial_{z_5}+(4\,w_2\,z_4
\\
&
\ \ \ \ \ \ \ \ \ \ \ \ \
-2\,w_4\,z_2+2\,w_5\,z_2+4\,w_6\,z_1)\,\partial_{z_6}-(4\,w_3\,z_4+2\,w_4\,z_3-2\,w_5\,z_3-4\,w_7\,z_1)\,\partial_{z_7}
\\
&
\ \ \ \ \ \ \ \ \ \ \ \ \
-(2\,w_4\,z_4+2\,w_5\,z_4+4\,w_6\,z_3-4\,w_7\,z_2)\,\partial_{z_8}+4\,w_1^2\,\partial_{w_1}+4\,w_1\,w_2\,\partial_{w_2}
\\
&
\ \ \ \ \ \ \ \ \ \ \ \ \
+4\,w_1\,w_3\,\partial_{w_3}+4\,w_1\,w_4\,\partial_{w_4}+4\,w_1\,w_5\,\partial_{w_5}+4\,w_1\,w_6\,\partial_{w_6}
\\
&
\ \ \ \ \ \ \ \ \ \ \ \ \
+4\,w_1\,w_7\,\partial_{w_7}+(4\,w_2\,w_7+4\,w_3\,w_6-w_4^2+w_5^2)\,\partial_{w_8}
\Big),
\endaligned
\]
\[
\aligned
IL_{w_2w_2}^2
&
\,:=\,
\isqrt\,
\Big(
4\,z_1\,w_2\,\partial_{z_1}+4\,z_2\,w_2\,\partial_{z_2}+(4\,w_1\,z_5-4\,w_3\,z_2+2\,w_4\,z_1+2\,w_5\,z_1)\,\partial_{z_3}
\\
&
\ \ \ \ \ \ \ \ \ \ \ \ \
+(4\,w_1\,z_6+2\,w_4\,z_2-2\,w_5\,z_2-4\,w_6\,z_1)\,\partial_{z_4}+4\,z_5\,w_2\,\partial_{z_5}+4\,z_6\,w_2\,\partial_{z_6}
\\
&
\ \ \ \ \ \ \ \ \ \ \ \ \
-(4\,w_3\,z_6+2\,w_4\,z_5-2\,w_5\,z_5-4\,w_8\,z_1)\,\partial_{z_7}-(2\,w_4\,z_6+2\,w_5\,z_6+4\,w_6\,z_5-4\,w_8\,z_2)\,\partial_{z_8}
\\
&
\ \ \ \ \ \ \ \ \ \ \ \ \
+4\,w_1\,w_2\,\partial_{w_1}+4\,w_2^2\,\partial_{w_2}+4\,w_2\,w_3\,\partial_{w_3}+4\,w_2\,w_4\,\partial_{w_4}
\\
&
\ \ \ \ \ \ \ \ \ \ \ \ \
+4\,w_2\,w_5\,\partial_{w_5}+4\,w_2\,w_6\,\partial_{w_6}+(4\,w_1\,w_8-4\,w_3\,w_6+w_4^2-w_5^2)\,\partial_{w_7}+4\,w_2\,w_8\,\partial_{w_8}
\Big),
\endaligned
\]
\[
\aligned
IL_{w_5w_5}^2
&
\,:=\,
\isqrt\,
\Big(
(2\,w_1\,z_5-2\,w_2\,z_3-2\,w_3\,z_2+w_4\,z_1-w_5\,z_1)\,\partial_{z_1}
\\
&
\ \ \ \ \ \ \ \ \ \ \ \ \
-(2\,w_1\,z_6-2\,w_2\,z_4+w_4\,z_2+w_5\,z_2-2\,w_6\,z_1)\,\partial_{z_2}
\\
&
\ \ \ \ \ \ \ \ \ \ \ \ \
-(2\,w_1\,z_7+2\,w_3\,z_4+w_4\,z_3+w_5\,z_3-2\,w_7\,z_1)\,\partial_{z_3}
\\
&
\ \ \ \ \ \ \ \ \ \ \ \ \
+(2\,w_1\,z_8+w_4\,z_4-w_5\,z_4+2\,w_6\,z_3-2\,w_7\,z_2)\,\partial_{z_4}
\\
&
\ \ \ \ \ \ \ \ \ \ \ \ \
-(2\,w_2\,z_7+2\,w_3\,z_6+w_4\,z_5+w_5\,z_5-2\,w_8\,z_1)\,\partial_{z_5}
\\
&
\ \ \ \ \ \ \ \ \ \ \ \ \
+(2\,w_2\,z_8+w_4\,z_6-w_5\,z_6+2\,w_6\,z_5-2\,w_8\,z_2)\,\partial_{z_6}
\\
&
\ \ \ \ \ \ \ \ \ \ \ \ \
-(2\,w_3\,z_8-w_4\,z_7+w_5\,z_7-2\,w_7\,z_5+2\,w_8\,z_3)\,\partial_{z_7}
\\
&
\ \ \ \ \ \ \ \ \ \ \ \ \
-(w_4\,z_8+w_5\,z_8-2\,w_6\,z_7+2\,w_7\,z_6-2\,w_8\,z_4)\,\partial_{z_8}
\\
&
\ \ \ \ \ \ \ \ \ \ \ \ \
-2\,w_1\,w_5\,\partial_{w_1}-2\,w_2\,w_5\,\partial_{w_2}-2\,w_3\,w_5\,\partial_{w_3}-2\,w_4\,w_5\,\partial_{w_4}
\\
&
\ \ \ \ \ \ \ \ \ \ \ \ \
-(4\,w_1\,w_8-4\,w_2\,w_7-4\,w_3\,w_6+w_4^2+w_5^2)\,\partial_{w_5}
\\
&
\ \ \ \ \ \ \ \ \ \ \ \ \
-2\,w_5\,w_6\,\partial_{w_6}-2\,w_5\,w_7\,\partial_{w_7}-2\,w_5\,w_8\,\partial_{w_8}
\Big),
\endaligned
\]
\[
\aligned
IL_{w_7w_7}^2
&
\,:=\,
\isqrt\,
\Big(
(4\,w_1\,z_7+4\,w_3\,z_4+2\,w_4\,z_3-2\,w_5\,z_3)\,\partial_{z_1}+(4\,w_1\,z_8+2\,w_4\,z_4+2\,w_5\,z_4+4\,w_6\,z_3)\,\partial_{z_2}
\\
&
\ \ \ \ \ \ \ \ \ \ \ \ \
+4\,z_3\,w_7\,\partial_{z_3}+4\,z_4\,w_7\,\partial_{z_4}+(4\,w_3\,z_8-2\,w_4\,z_7-2\,w_5\,z_7+4\,w_8\,z_3)\,\partial_{z_5}
\\
&
\ \ \ \ \ \ \ \ \ \ \ \ \
-(2\,w_4\,z_8-2\,w_5\,z_8-4\,w_6\,z_7-4\,w_8\,z_4)\,\partial_{z_6}+4\,z_7\,w_7\,\partial_{z_7}+4\,z_8\,w_7\,\partial_{z_8}
\\
&
\ \ \ \ \ \ \ \ \ \ \ \ \
+4\,w_1\,w_7\,\partial_{w_1}+(4\,w_1\,w_8-4\,w_3\,w_6+w_4^2-w_5^2)\,\partial_{w_2}+4\,w_3\,w_7\,\partial_{w_3}
\\
&
\ \ \ \ \ \ \ \ \ \ \ \ \
+4\,w_4\,w_7\,\partial_{w_4}+4\,w_5\,w_7\,\partial_{w_5}+4\,w_6\,w_7\,\partial_{w_6}+4\,w_7^2\,\partial_{w_7}+4\,w_7\,w_8\,\partial_{w_8}
\Big),
\endaligned
\]
\[
\aligned
IL_{w_8w_8}^2
&
\,:=\,
\isqrt\,
\Big(
(4\,w_2\,z_7+4\,w_3\,z_6+2\,w_4\,z_5-2\,w_5\,z_5)\,\partial_{z_1}+(4\,w_2\,z_8+2\,w_4\,z_6+2\,w_5\,z_6+4\,w_6\,z_5)\,\partial_{z_2}
\\
&
\ \ \ \ \ \ \ \ \ \ \ \ \
-(4\,w_3\,z_8-2\,w_4\,z_7-2\,w_5\,z_7-4\,w_7\,z_5)\,\partial_{z_3}+(2\,w_4\,z_8-2\,w_5\,z_8-4\,w_6\,z_7+4\,w_7\,z_6)\,\partial_{z_4}
\\
&
\ \ \ \ \ \ \ \ \ \ \ \ \
+4\,z_5\,w_8\,\partial_{z_5}+4\,z_6\,w_8\,\partial_{z_6}+4\,z_7\,w_8\,\partial_{z_7}+4\,z_8\,w_8\,\partial_{z_8}
\\
&
\ \ \ \ \ \ \ \ \ \ \ \ \
+(4\,w_2\,w_7+4\,w_3\,w_6-w_4^2+w_5^2)\,\partial_{w_1}+4\,w_2\,w_8\,\partial_{w_2}+4\,w_3\,w_8\,\partial_{w_3}
\\
&
\ \ \ \ \ \ \ \ \ \ \ \ \
+4\,w_4\,w_8\,\partial_{w_4}
+4\,w_5\,w_8\,\partial_{w_5}+4\,w_6\,w_8\,\partial_{w_6}+4\,w_7\,w_8\,\partial_{w_7}+4\,w_8^2\,\partial_{w_8}
\Big).
\endaligned
\]

\section{$E_{III}$ CR Models}
\label{E-III-CR-models}

In $\mathbb{C}^{8 + 8}$, consider:
\[
\aligned
\Re\,w_1
&
\,=\,
\Re\,
\big(
\vert z_1\vert^2
+
\vert z_2\vert^2
+
\vert z_3\vert^2
+
\vert z_4\vert^2
\big),
\\
\Re\,w_2
&
\,=\,
\Re\,
\big(
\vert z_5\vert^2
+
\vert z_6\vert^2
+
\vert z_7\vert^2
+
\vert z_8\vert^2
\big),
\\
\Im\,w_3
&
\,=\,
\Im\,
\big(
z_1\,\overline{z}_7
+
z_2\,\overline{z}_8
+
z_5\,\overline{z}_3
+
z_6\,\overline{z}_4
\big),
\\
\Re\,w_4
&
\,=\,
\Re\,
\big(
z_1\,\overline{z}_7
+
z_2\,\overline{z}_8
+
z_5\,\overline{z}_3
+
z_6\,\overline{z}_4
\big),
\\
\Im\,w_5
&
\,=\,
\Im\,
\big(
z_1\,\overline{z}_6
-
z_3\,\overline{z}_8
+
z_5\,\overline{z}_2
-
z_7\,\overline{z}_4
\big),
\\
\Re\,w_6
&
\,=\,
\Re\,
\big(
z_1\,\overline{z}_6
-
z_3\,\overline{z}_8
+
z_5\,\overline{z}_2
-
z_7\,\overline{z}_4
\big),
\\
\Im\,w_7
&
\,=\,
\Im\,
\big(
z_2\,\overline{z}_6
+
z_3\,\overline{z}_7
-
z_5\,\overline{z}_1
-
z_8\,\overline{z}_4
\big),
\\
\Re\,w_8
&
\,=\,
\Re\,
\big(
z_2\,\overline{z}_6
+
z_3\,\overline{z}_7
-
z_5\,\overline{z}_1
-
z_8\,\overline{z}_4
\big).
\endaligned
\]

The $3 + 5$ generators of $\mathfrak{g}_{-2}$ are:
\[
\aligned
L_{w_3}^{-2}
&
\,:=\,
\partial_{w_3},
\\
L_{w_5}^{-2}
&
\,:=\,
\partial_{w_5},
\\
L_{w_7}^{-2}
&
\,:=\,
\partial_{w_7},
\endaligned
\ \ \ \ \ \ \ \ \ \ \ \ \ \ \ \ \ \ \ \ \ \ \ \ \ \
\aligned
IL_{w_1}^{-2}
&
\,:=\,
\isqrt\,\partial_{w_1},
\\
IL_{w_2}^{-2}
&
\,:=\,
\isqrt\,\partial_{w_2},
\\
IL_{w_4}^{-2}
&
\,:=\,
\isqrt\,\partial_{w_4},
\\
IL_{w_6}^{-2}
&
\,:=\,
\isqrt\,\partial_{w_6},
\\
IL_{w_8}^{-2}
&
\,:=\,
\isqrt\,\partial_{w_8}.
\endaligned
\]

The $8 + 8$ generators of $\mathfrak{g}_{-1}$ are:
\[
\aligned
L_{z_1}^{-1}
&
\,:=\,
\partial_{z_1}+z_1\,\partial_{w_1}-z_7\,\partial_{w_3}+z_7\,\partial_{w_4}-z_6\,\partial_{w_5}+z_6\,\partial_{w_6}-z_5\,\partial_{w_7}-z_5\,\partial_{w_8},
\\
L_{z_2}^{-1}
&
\,:=\,
\partial_{z_2}+z_2\,\partial_{w_1}-z_8\,\partial_{w_3}+z_8\,\partial_{w_4}+z_5\,\partial_{w_5}+z_5\,\partial_{w_6}-z_6\,\partial_{w_7}+z_6\,\partial_{w_8},
\\
L_{z_3}^{-1}
&
\,:=\,
\partial_{z_3}+z_3\,\partial_{w_1}+z_5\,\partial_{w_3}+z_5\,\partial_{w_4}+z_8\,\partial_{w_5}-z_8\,\partial_{w_6}-z_7\,\partial_{w_7}+z_7\,\partial_{w_8},
\\
L_{z_4}^{-1}
&
\,:=\,
\partial_{z_4}+z_4\,\partial_{w_1}+z_6\,\partial_{w_3}+z_6\,\partial_{w_4}-z_7\,\partial_{w_5}-z_7\,\partial_{w_6}-z_8\,\partial_{w_7}-z_8\,\partial_{w_8},
\\
L_{z_5}^{-1}
&
\,:=\,
\partial_{z_5}+z_5\,\partial_{w_2}-z_3\,\partial_{w_3}+z_3\,\partial_{w_4}-z_2\,\partial_{w_5}+z_2\,\partial_{w_6}+z_1\,\partial_{w_7}-z_1\,\partial_{w_8},
\\
L_{z_6}^{-1}
&
\,:=\,
\partial_{z_6}+z_6\,\partial_{w_2}-z_4\,\partial_{w_3}+z_4\,\partial_{w_4}+z_1\,\partial_{w_5}+z_1\,\partial_{w_6}+z_2\,\partial_{w_7}+z_2\,\partial_{w_8},
\\
L_{z_7}^{-1}
&
\,:=\,
\partial_{z_7}+z_7\,\partial_{w_2}+z_1\,\partial_{w_3}+z_1\,\partial_{w_4}+z_4\,\partial_{w_5}-z_4\,\partial_{w_6}+z_3\,\partial_{w_7}+z_3\,\partial_{w_8},
\\
L_{z_8}^{-1}
&
\,:=\,
\partial_{z_8}+z_8\,\partial_{w_2}+z_2\,\partial_{w_3}+z_2\,\partial_{w_4}-z_3\,\partial_{w_5}-z_3\,\partial_{w_6}+z_4\,\partial_{w_7}-z_4\,\partial_{w_8},
\endaligned
\]
\[
\aligned
IL_{z_1}^{-1}
&
\,:=\,
\isqrt\,
\Big(
\partial_{z_1}-z_1\,\partial_{w_1}+z_7\,\partial_{w_3}-z_7\,\partial_{w_4}+z_6\,\partial_{w_5}-z_6\,\partial_{w_6}+z_5\,\partial_{w_7}+z_5\,\partial_{w_8}
\Big),
\\
IL_{z_2}^{-1}
&
\,:=\,
\isqrt\,
\Big(
\partial_{z_2}-z_2\,\partial_{w_1}+z_8\,\partial_{w_3}-z_8\,\partial_{w_4}-z_5\,\partial_{w_5}-z_5\,\partial_{w_6}+z_6\,\partial_{w_7}-z_6\,\partial_{w_8}
\Big),
\\
IL_{z_3}^{-1}
&
\,:=\,
\isqrt\,
\Big(
\partial_{z_3}-z_3\,\partial_{w_1}-z_5\,\partial_{w_3}-z_5\,\partial_{w_4}-z_8\,\partial_{w_5}+z_8\,\partial_{w_6}+z_7\,\partial_{w_7}-z_7\,\partial_{w_8}
\Big),
\\
IL_{z_4}^{-1}
&
\,:=\,
\isqrt\,
\Big(
\partial_{z_4}-z_4\,\partial_{w_1}-z_6\,\partial_{w_3}-z_6\,\partial_{w_4}+z_7\,\partial_{w_5}+z_7\,\partial_{w_6}+z_8\,\partial_{w_7}+z_8\,\partial_{w_8}
\Big),
\\
IL_{z_5}^{-1}
&
\,:=\,
\isqrt\,
\Big(
\partial_{z_5}-z_5\,\partial_{w_2}+z_3\,\partial_{w_3}-z_3\,\partial_{w_4}+z_2\,\partial_{w_5}-z_2\,\partial_{w_6}-z_1\,\partial_{w_7}+z_1\,\partial_{w_8}
\Big),
\\
IL_{z_6}^{-1}
&
\,:=\,
\isqrt\,
\Big(
\partial_{z_6}-z_6\,\partial_{w_2}+z_4\,\partial_{w_3}-z_4\,\partial_{w_4}-z_1\,\partial_{w_5}-z_1\,\partial_{w_6}-z_2\,\partial_{w_7}-z_2\,\partial_{w_8}
\Big),
\\
IL_{z_7}^{-1}
&
\,:=\,
\isqrt\,
\Big(
\partial_{z_7}-z_7\,\partial_{w_2}-z_1\,\partial_{w_3}-z_1\,\partial_{w_4}-z_4\,\partial_{w_5}+z_4\,\partial_{w_6}-z_3\,\partial_{w_7}-z_3\,\partial_{w_8}
\Big),
\\
IL_{z_8}^{-1}
&
\,:=\,
\isqrt\,
\Big(
\partial_{z_8}-z_8\,\partial_{w_2}-z_2\,\partial_{w_3}-z_2\,\partial_{w_4}+z_3\,\partial_{w_5}+z_3\,\partial_{w_6}-z_4\,\partial_{w_7}+z_4\,\partial_{w_8}
\Big).
\endaligned
\]

The $14 + 16$ generators of $\mathfrak{g}_0$ are:
\[
\aligned
L_1^0
&
\,:=\,
z_7\,\partial_{z_1}+z_8\,\partial_{z_2}+z_5\,\partial_{z_3}+z_6\,\partial_{z_4}+w_4\,\partial_{w_1}+2\,w_2\,\partial_{w_4},
\\
L_2^0
&
\,:=\,
z_6\,\partial_{z_1}+z_5\,\partial_{z_2}-z_8\,\partial_{z_3}-z_7\,\partial_{z_4}+w_6\,\partial_{w_1}+2\,w_2\,\partial_{w_6},
\\
L_3^0
&
\,:=\,
-z_5\,\partial_{z_1}+z_6\,\partial_{z_2}+z_7\,\partial_{z_3}-z_8\,\partial_{z_4}+w_8\,\partial_{w_1}+2\,w_2\,\partial_{w_8},
\\
L_4^0
&
\,:=\,
z_3\,\partial_{z_5}+z_4\,\partial_{z_6}+z_1\,\partial_{z_7}+z_2\,\partial_{z_8}+w_4\,\partial_{w_2}+2\,w_1\,\partial_{w_4},
\\
L_5^0
&
\,:=\,
z_2\,\partial_{z_5}+z_1\,\partial_{z_6}-z_4\,\partial_{z_7}-z_3\,\partial_{z_8}+w_6\,\partial_{w_2}+2\,w_1\,\partial_{w_6},
\\
L_6^0
&
\,:=\,
-z_1\,\partial_{z_5}+z_2\,\partial_{z_6}+z_3\,\partial_{z_7}-z_4\,\partial_{z_8}+w_8\,\partial_{w_2}+2\,w_1\,\partial_{w_8},
\\
L_7^0
&
\,:=\,
z_1\,\partial_{z_1}+z_2\,\partial_{z_2}+z_3\,\partial_{z_3}+z_4\,\partial_{z_4}+2\,w_1\,\partial_{w_1}+w_3\,\partial_{w_3}+w_4\,\partial_{w_4}+w_5\,\partial_{w_5}
\\
&
\ \ \ \ \
+w_6\,\partial_{w_6}+w_7\,\partial_{w_7}+w_8\,\partial_{w_8},
\\
L_8^0
&
\,:=\,
z_4\,\partial_{z_1}-z_1\,\partial_{z_4}+z_8\,\partial_{z_5}-z_5\,\partial_{z_8}+w_5\,\partial_{w_3}-w_6\,\partial_{w_4}-w_3\,\partial_{w_5}+w_4\,\partial_{w_6},
\\
L_9^0
&
\,:=\,
-z_4\,\partial_{z_3}+z_3\,\partial_{z_4}-z_6\,\partial_{z_5}+z_5\,\partial_{z_6}+w_7\,\partial_{w_5}-w_8\,\partial_{w_6}-w_5\,\partial_{w_7}+w_6\,\partial_{w_8},
\\
L_{10}^0
&
\,:=\,
z_3\,\partial_{z_2}-z_2\,\partial_{z_3}+z_7\,\partial_{z_6}-z_6\,\partial_{z_7}-w_5\,\partial_{w_3}-w_6\,\partial_{w_4}+w_3\,\partial_{w_5}+w_4\,\partial_{w_6},
\\
L_{11}^0
&
\,:=\,
z_4\,\partial_{z_2}-z_2\,\partial_{z_4}-z_7\,\partial_{z_5}+z_5\,\partial_{z_7}+w_7\,\partial_{w_3}-w_8\,\partial_{w_4}-w_3\,\partial_{w_7}+w_4\,\partial_{w_8},
\\
L_{12}^0
&
\,:=\,
z_3\,\partial_{z_1}-z_1\,\partial_{z_3}-z_8\,\partial_{z_6}+z_6\,\partial_{z_8}+w_7\,\partial_{w_3}+w_8\,\partial_{w_4}-w_3\,\partial_{w_7}-w_4\,\partial_{w_8},
\\
L_{13}^0
&
\,:=\,
-z_2\,\partial_{z_1}+z_1\,\partial_{z_2}-z_8\,\partial_{z_7}+z_7\,\partial_{z_8}-w_7\,\partial_{w_5}-w_8\,\partial_{w_6}+w_5\,\partial_{w_7}+w_6\,\partial_{w_8},
\\
L_{14}^0
&
\,:=\,
z_5\,\partial_{z_5}+z_6\,\partial_{z_6}+z_7\,\partial_{z_7}+z_8\,\partial_{z_8}+2\,w_2\,\partial_{w_2}+w_3\,\partial_{w_3}+w_4\,\partial_{w_4}+w_5\,\partial_{w_5}
\\
&
\ \ \ \ \
+w_6\,\partial_{w_6}+w_7\,\partial_{w_7}+w_8\,\partial_{w_8},
\endaligned
\]
\[
\aligned
IL_1^0
&
\,:=\,
\isqrt\,
\Big(
z_7\,\partial_{z_1}+z_8\,\partial_{z_2}-z_5\,\partial_{z_3}-z_6\,\partial_{z_4}-w_3\,\partial_{w_1}+2\,w_2\,\partial_{w_3}
\Big),
\\
IL_2^0
&
\,:=\,
\isqrt\,
\Big(
z_6\,\partial_{z_1}-z_5\,\partial_{z_2}-z_8\,\partial_{z_3}+z_7\,\partial_{z_4}-w_5\,\partial_{w_1}+2\,w_2\,\partial_{w_5}
\Big),
\\
IL_3^0
&
\,:=\,
\isqrt\,
\Big(
z_5\,\partial_{z_1}+z_6\,\partial_{z_2}+z_7\,\partial_{z_3}+z_8\,\partial_{z_4}-w_7\,\partial_{w_1}+2\,w_2\,\partial_{w_7}
\Big),
\\
IL_4^0
&
\,:=\,
\isqrt\,
\Big(
z_3\,\partial_{z_5}+z_4\,\partial_{z_6}-z_1\,\partial_{z_7}-z_2\,\partial_{z_8}-w_3\,\partial_{w_2}+2\,w_1\,\partial_{w_3}
\Big),
\\
IL_5^0
&
\,:=\,
\isqrt\,
\Big(
z_2\,\partial_{z_5}-z_1\,\partial_{z_6}-z_4\,\partial_{z_7}+z_3\,\partial_{z_8}-w_5\,\partial_{w_2}+2\,w_1\,\partial_{w_5}
\Big),
\\
IL_6^0
&
\,:=\,
\isqrt\,
\Big(
z_1\,\partial_{z_5}+z_2\,\partial_{z_6}+z_3\,\partial_{z_7}+z_4\,\partial_{z_8}+w_7\,\partial_{w_2}-2\,w_1\,\partial_{w_7}
\Big),
\\
IL_7^0
&
\,:=\,
\isqrt\,
\Big(
z_1\,\partial_{z_1}-z_4\,\partial_{z_4}+z_5\,\partial_{z_5}-z_8\,\partial_{z_8}+w_4\,\partial_{w_3}+w_3\,\partial_{w_4}+w_6\,\partial_{w_5}+w_5\,\partial_{w_6}
\Big),
\\
IL_8^0
&
\,:=\,
\isqrt\,
\Big(
z_2\,\partial_{z_2}+z_4\,\partial_{z_4}+z_6\,\partial_{z_6}+z_8\,\partial_{z_8}-w_6\,\partial_{w_5}-w_5\,\partial_{w_6}
\Big),
\\
IL_9^0
&
\,:=\,
\isqrt\,
\Big(
z_3\,\partial_{z_3}-z_4\,\partial_{z_4}+z_5\,\partial_{z_5}-z_6\,\partial_{z_6}+w_6\,\partial_{w_5}+w_5\,\partial_{w_6}+w_8\,\partial_{w_7}+w_7\,\partial_{w_8}
\Big),
\\
IL_{10}^0
&
\,:=\,
\isqrt\,
\Big(
z_4\,\partial_{z_3}+z_3\,\partial_{z_4}+z_6\,\partial_{z_5}+z_5\,\partial_{z_6}+w_8\,\partial_{w_5}-w_7\,\partial_{w_6}-w_6\,\partial_{w_7}+w_5\,\partial_{w_8}
\Big),
\\
IL_{11}^0
&
\,:=\,
\isqrt\,
\Big(
z_4\,\partial_{z_2}+z_2\,\partial_{z_4}-z_7\,\partial_{z_5}-z_5\,\partial_{z_7}-w_8\,\partial_{w_3}+w_7\,\partial_{w_4}+w_4\,\partial_{w_7}-w_3\,\partial_{w_8}
\Big),
\\
IL_{12}^0
&
\,:=\,
\isqrt\,
\Big(
z_4\,\partial_{z_1}+z_1\,\partial_{z_4}+z_8\,\partial_{z_5}+z_5\,\partial_{z_8}-w_6\,\partial_{w_3}+w_5\,\partial_{w_4}+w_4\,\partial_{w_5}-w_3\,\partial_{w_6}
\Big),
\\
IL_{13}^0
&
\,:=\,
\isqrt\,
\Big(
z_3\,\partial_{z_1}+z_1\,\partial_{z_3}-z_8\,\partial_{z_6}-z_6\,\partial_{z_8}+w_8\,\partial_{w_3}+w_7\,\partial_{w_4}+w_4\,\partial_{w_7}+w_3\,\partial_{w_8}
\Big),
\\
IL_{14}^0
&
\,:=\,
\isqrt\,
\Big(
z_3\,\partial_{z_2}+z_2\,\partial_{z_3}+z_7\,\partial_{z_6}+z_6\,\partial_{z_7}-w_6\,\partial_{w_3}-w_5\,\partial_{w_4}-w_4\,\partial_{w_5}-w_3\,\partial_{w_6}
\Big),
\\
IL_{15}^0
&
\,:=\,
\isqrt\,
\Big(
2\,z_4\,\partial_{z_4}-z_5\,\partial_{z_5}+z_6\,\partial_{z_6}+z_7\,\partial_{z_7}+z_8\,\partial_{z_8}-w_4\,\partial_{w_3}-w_3\,\partial_{w_4}-w_6\,\partial_{w_5}
\\
&
\ \ \ \ \ \ \ \ \ \ \ \ \
-w_5\,\partial_{w_6}-w_8\,\partial_{w_7}-w_7\,\partial_{w_8}
\Big),
\\
IL_{16}^0
&
\,:=\,
\isqrt\,
\Big(
z_2\,\partial_{z_1}+z_1\,\partial_{z_2}+z_8\,\partial_{z_7}+z_7\,\partial_{z_8}+w_8\,\partial_{w_5}+w_7\,\partial_{w_6}+w_6\,\partial_{w_7}+w_5\,\partial_{w_8}
\Big).
\endaligned
\]

The $8 + 8$ generators of $\mathfrak{g}_1$ are:
\[
\aligned
L_{z_1z_1}^1
&
\,:=\,
-(-2\,z_1^2+2\,w_1)\,\partial_{z_1}+2\,z_1\,z_2\,\partial_{z_2}+2\,z_1\,z_3\,\partial_{z_3}+2\,z_1\,z_4\,\partial_{z_4}
\\
&
\ \ \ \ \
+(2\,z_1\,z_5+w_7+w_8)\,\partial_{z_5}+(2\,z_1\,z_6+w_5-w_6)\,\partial_{z_6}+(2\,z_1\,z_7+w_3-w_4)\,\partial_{z_7}
\\
&
\ \ \ \ \
+(2\,z_2\,z_7-2\,z_3\,z_6+2\,z_4\,z_5)\,\partial_{z_8}+2\,z_1\,w_1\,\partial_{w_1}
\\
&
\ \ \ \ \
+(w_3\,z_7+w_4\,z_7+w_5\,z_6+w_6\,z_6+w_7\,z_5-w_8\,z_5)\,\partial_{w_2}
\\
&
\ \ \ \ \
-(2\,w_1\,z_7-w_3\,z_1-w_4\,z_1+w_5\,z_4+w_6\,z_4+w_7\,z_3-w_8\,z_3)\,\partial_{w_3}
\\
&
\ \ \ \ \
+(2\,w_1\,z_7+w_3\,z_1+w_4\,z_1+w_5\,z_4+w_6\,z_4+w_7\,z_3-w_8\,z_3)\,\partial_{w_4}
\\
&
\ \ \ \ \
-(2\,w_1\,z_6-w_3\,z_4-w_4\,z_4-w_5\,z_1-w_6\,z_1+w_7\,z_2-w_8\,z_2)\,\partial_{w_5}
\\
&
\ \ \ \ \
+(2\,w_1\,z_6-w_3\,z_4-w_4\,z_4+w_5\,z_1+w_6\,z_1+w_7\,z_2-w_8\,z_2)\,\partial_{w_6}
\\
&
\ \ \ \ \
-(2\,w_1\,z_5-w_3\,z_3-w_4\,z_3-w_5\,z_2-w_6\,z_2-w_7\,z_1+w_8\,z_1)\,\partial_{w_7}
\\
&
\ \ \ \ \
-(2\,w_1\,z_5-w_3\,z_3-w_4\,z_3-w_5\,z_2-w_6\,z_2+w_7\,z_1-w_8\,z_1)\,\partial_{w_8},
\endaligned
\]
\[
\aligned
L_{z_2z_2}^1
&
\,:=\,
2\,z_1\,z_2\,\partial_{z_1}-(-2\,z_2^2+2\,w_1)\,\partial_{z_2}+2\,z_2\,z_3\,\partial_{z_3}+2\,z_2\,z_4\,\partial_{z_4}-(-2\,z_2\,z_5+w_5+w_6)\,\partial_{z_5}
\\
&
\ \ \ \ \
+(2\,z_2\,z_6+w_7-w_8)\,\partial_{z_6}+(2\,z_1\,z_8+2\,z_3\,z_6-2\,z_4\,z_5)\,\partial_{z_7}+(2\,z_2\,z_8+w_3-w_4)\,\partial_{z_8}
\\
&
\ \ \ \ \
+2\,z_2\,w_1\,\partial_{w_1}+(w_3\,z_8+w_4\,z_8-w_5\,z_5+w_6\,z_5+w_7\,z_6+w_8\,z_6)\,\partial_{w_2}
\\
&
\ \ \ \ \
-(2\,w_1\,z_8-w_3\,z_2-w_4\,z_2-w_5\,z_3+w_6\,z_3+w_7\,z_4+w_8\,z_4)\,\partial_{w_3}
\\
&
\ \ \ \ \
+(2\,w_1\,z_8+w_3\,z_2+w_4\,z_2-w_5\,z_3+w_6\,z_3+w_7\,z_4+w_8\,z_4)\,\partial_{w_4}
\\
&
\ \ \ \ \
+(2\,w_1\,z_5-w_3\,z_3-w_4\,z_3+w_5\,z_2-w_6\,z_2+w_7\,z_1+w_8\,z_1)\,\partial_{w_5}
\\
&
\ \ \ \ \
+(2\,w_1\,z_5-w_3\,z_3-w_4\,z_3-w_5\,z_2+w_6\,z_2+w_7\,z_1+w_8\,z_1)\,\partial_{w_6}
\\
&
\ \ \ \ \
-(2\,w_1\,z_6-w_3\,z_4-w_4\,z_4+w_5\,z_1-w_6\,z_1-w_7\,z_2-w_8\,z_2)\,\partial_{w_7}
\\
&
\ \ \ \ \
+(2\,w_1\,z_6-w_3\,z_4-w_4\,z_4+w_5\,z_1-w_6\,z_1+w_7\,z_2+w_8\,z_2)\,\partial_{w_8},
\endaligned
\]
\[
\aligned
L_{z_3z_3}^1
&
\,:=\,
2\,z_1\,z_3\,\partial_{z_1}+2\,z_2\,z_3\,\partial_{z_2}-(-2\,z_3^2+2\,w_1)\,\partial_{z_3}+2\,z_3\,z_4\,\partial_{z_4}-(-2\,z_3\,z_5+w_3+w_4)\,\partial_{z_5}
\\
&
\ \ \ \ \
-(2\,z_1\,z_8-2\,z_2\,z_7-2\,z_4\,z_5)\,\partial_{z_6}+(2\,z_3\,z_7+w_7-w_8)\,\partial_{z_7}-(-2\,z_3\,z_8+w_5-w_6)\,\partial_{z_8}
\\
&
\ \ \ \ \
+2\,z_3\,w_1\,\partial_{w_1}-(w_3\,z_5-w_4\,z_5+w_5\,z_8+w_6\,z_8-w_7\,z_7-w_8\,z_7)\,\partial_{w_2}
\\
&
\ \ \ \ \
+(2\,w_1\,z_5+w_3\,z_3-w_4\,z_3-w_5\,z_2-w_6\,z_2+w_7\,z_1+w_8\,z_1)\,\partial_{w_3}
\\
&
\ \ \ \ \
+(2\,w_1\,z_5-w_3\,z_3+w_4\,z_3-w_5\,z_2-w_6\,z_2+w_7\,z_1+w_8\,z_1)\,\partial_{w_4}
\\
&
\ \ \ \ \
+(2\,w_1\,z_8+w_3\,z_2-w_4\,z_2+w_5\,z_3+w_6\,z_3+w_7\,z_4+w_8\,z_4)\,\partial_{w_5}
\\
&
\ \ \ \ \
-(2\,w_1\,z_8+w_3\,z_2-w_4\,z_2-w_5\,z_3-w_6\,z_3+w_7\,z_4+w_8\,z_4)\,\partial_{w_6}
\\
&
\ \ \ \ \
-(2\,w_1\,z_7+w_3\,z_1-w_4\,z_1+w_5\,z_4+w_6\,z_4-w_7\,z_3-w_8\,z_3)\,\partial_{w_7}
\\
&
\ \ \ \ \
+(2\,w_1\,z_7+w_3\,z_1-w_4\,z_1+w_5\,z_4+w_6\,z_4+w_7\,z_3+w_8\,z_3)\,\partial_{w_8},
\endaligned
\]
\[
\aligned
L_{z_4z_4}^1
&
\,:=\,
2\,z_1\,z_4\,\partial_{z_1}+2\,z_2\,z_4\,\partial_{z_2}+2\,z_3\,z_4\,\partial_{z_3}-(-2\,z_4^2+2\,w_1)\,\partial_{z_4}+(2\,z_1\,z_8-2\,z_2\,z_7+2\,z_3\,z_6)\,\partial_{z_5}
\\
&
\ \ \ \ \
-(-2\,z_4\,z_6+w_3+w_4)\,\partial_{z_6}+(2\,z_4\,z_7+w_5+w_6)\,\partial_{z_7}+(2\,z_4\,z_8+w_7+w_8)\,\partial_{z_8}+2\,z_4\,w_1\,\partial_{w_1}
\\
&
\ \ \ \ \
-(w_3\,z_6-w_4\,z_6-w_5\,z_7+w_6\,z_7-w_7\,z_8+w_8\,z_8)\,\partial_{w_2}
\\
&
\ \ \ \ \
+(2\,w_1\,z_6+w_3\,z_4-w_4\,z_4+w_5\,z_1-w_6\,z_1+w_7\,z_2-w_8\,z_2)\,\partial_{w_3}
\\
&
\ \ \ \ \
+(2\,w_1\,z_6-w_3\,z_4+w_4\,z_4+w_5\,z_1-w_6\,z_1+w_7\,z_2-w_8\,z_2)\,\partial_{w_4}
\\
&
\ \ \ \ \
-(2\,w_1\,z_7+w_3\,z_1-w_4\,z_1-w_5\,z_4+w_6\,z_4+w_7\,z_3-w_8\,z_3)\,\partial_{w_5}
\\
&
\ \ \ \ \
-(2\,w_1\,z_7+w_3\,z_1-w_4\,z_1+w_5\,z_4-w_6\,z_4+w_7\,z_3-w_8\,z_3)\,\partial_{w_6}
\\
&
\ \ \ \ \
-(2\,w_1\,z_8+w_3\,z_2-w_4\,z_2-w_5\,z_3+w_6\,z_3-w_7\,z_4+w_8\,z_4)\,\partial_{w_7}
\\
&
\ \ \ \ \
-(2\,w_1\,z_8+w_3\,z_2-w_4\,z_2-w_5\,z_3+w_6\,z_3+w_7\,z_4-w_8\,z_4)\,\partial_{w_8},
\endaligned
\]
\[
\aligned
L_{z_5z_5}^1
&
\,:=\,
-(-2\,z_1\,z_5+w_7-w_8)\,\partial_{z_1}+(2\,z_2\,z_5+w_5-w_6)\,\partial_{z_2}+(2\,z_3\,z_5+w_3-w_4)\,\partial_{z_3}
\\
&
\ \ \ \ \
+(2\,z_1\,z_8-2\,z_2\,z_7+2\,z_3\,z_6)\,\partial_{z_4}-(-2\,z_5^2+2\,w_2)\,\partial_{z_5}+2\,z_5\,z_6\,\partial_{z_6}+2\,z_5\,z_7\,\partial_{z_7}+2\,z_5\,z_8\,\partial_{z_8}
\\
&
\ \ \ \ \
+(w_3\,z_3+w_4\,z_3+w_5\,z_2+w_6\,z_2-w_7\,z_1-w_8\,z_1)\,\partial_{w_1}+2\,w_2\,z_5\,\partial_{w_2}
\\
&
\ \ \ \ \
-(2\,w_2\,z_3-w_3\,z_5-w_4\,z_5+w_5\,z_8+w_6\,z_8-w_7\,z_7-w_8\,z_7)\,\partial_{w_3}
\\
&
\ \ \ \ \
+(2\,w_2\,z_3+w_3\,z_5+w_4\,z_5+w_5\,z_8+w_6\,z_8-w_7\,z_7-w_8\,z_7)\,\partial_{w_4}
\\
&
\ \ \ \ \
-(2\,w_2\,z_2-w_3\,z_8-w_4\,z_8-w_5\,z_5-w_6\,z_5-w_7\,z_6-w_8\,z_6)\,\partial_{w_5}
\\
&
\ \ \ \ \
+(2\,w_2\,z_2-w_3\,z_8-w_4\,z_8+w_5\,z_5+w_6\,z_5-w_7\,z_6-w_8\,z_6)\,\partial_{w_6}
\\
&
\ \ \ \ \
+(2\,w_2\,z_1-w_3\,z_7-w_4\,z_7-w_5\,z_6-w_6\,z_6+w_7\,z_5+w_8\,z_5)\,\partial_{w_7}
\\
&
\ \ \ \ \
-(2\,w_2\,z_1-w_3\,z_7-w_4\,z_7-w_5\,z_6-w_6\,z_6-w_7\,z_5-w_8\,z_5)\,\partial_{w_8},
\endaligned
\]
\[
\aligned
L_{z_6z_6}^1
&
\,:=\,
-(-2\,z_1\,z_6+w_5+w_6)\,\partial_{z_1}-(-2\,z_2\,z_6+w_7+w_8)\,\partial_{z_2}-(2\,z_1\,z_8-2\,z_2\,z_7-2\,z_4\,z_5)\,\partial_{z_3}
\\
&
\ \ \ \ \
+(2\,z_4\,z_6+w_3-w_4)\,\partial_{z_4}+2\,z_5\,z_6\,\partial_{z_5}-(-2\,z_6^2+2\,w_2)\,\partial_{z_6}+2\,z_6\,z_7\,\partial_{z_7}+2\,z_6\,z_8\,\partial_{z_8}
\\
&
\ \ \ \ \
+(w_3\,z_4+w_4\,z_4-w_5\,z_1+w_6\,z_1-w_7\,z_2+w_8\,z_2)\,\partial_{w_1}+2\,w_2\,z_6\,\partial_{w_2}
\\
&
\ \ \ \ \
-(2\,w_2\,z_4-w_3\,z_6-w_4\,z_6-w_5\,z_7+w_6\,z_7-w_7\,z_8+w_8\,z_8)\,\partial_{w_3}
\\
&
\ \ \ \ \
+(2\,w_2\,z_4+w_3\,z_6+w_4\,z_6-w_5\,z_7+w_6\,z_7-w_7\,z_8+w_8\,z_8)\,\partial_{w_4}
\\
&
\ \ \ \ \
+(2\,w_2\,z_1-w_3\,z_7-w_4\,z_7+w_5\,z_6-w_6\,z_6-w_7\,z_5+w_8\,z_5)\,\partial_{w_5}
\\
&
\ \ \ \ \
+(2\,w_2\,z_1-w_3\,z_7-w_4\,z_7-w_5\,z_6+w_6\,z_6-w_7\,z_5+w_8\,z_5)\,\partial_{w_6}
\\
&
\ \ \ \ \
+(2\,w_2\,z_2-w_3\,z_8-w_4\,z_8+w_5\,z_5-w_6\,z_5+w_7\,z_6-w_8\,z_6)\,\partial_{w_7}
\\
&
\ \ \ \ \
+(2\,w_2\,z_2-w_3\,z_8-w_4\,z_8+w_5\,z_5-w_6\,z_5-w_7\,z_6+w_8\,z_6)\,\partial_{w_8},
\endaligned
\]
\[
\aligned
L_{z_7z_7}^1
&
\,:=\,
-(-2\,z_1\,z_7+w_3+w_4)\,\partial_{z_1}+(2\,z_1\,z_8+2\,z_3\,z_6-2\,z_4\,z_5)\,\partial_{z_2}-(-2\,z_3\,z_7+w_7+w_8)\,\partial_{z_3}
\\
&
\ \ \ \ \
-(-2\,z_4\,z_7+w_5-w_6)\,\partial_{z_4}+2\,z_5\,z_7\,\partial_{z_5}+2\,z_6\,z_7\,\partial_{z_6}-(-2\,z_7^2+2\,w_2)\,\partial_{z_7}+2\,z_7\,z_8\,\partial_{z_8}
\\
&
\ \ \ \ \
-(w_3\,z_1-w_4\,z_1+w_5\,z_4+w_6\,z_4+w_7\,z_3-w_8\,z_3)\,\partial_{w_1}+2\,w_2\,z_7\,\partial_{w_2}
\\
&
\ \ \ \ \
+(2\,w_2\,z_1+w_3\,z_7-w_4\,z_7-w_5\,z_6-w_6\,z_6-w_7\,z_5+w_8\,z_5)\,\partial_{w_3}
\\
&
\ \ \ \ \
+(2\,w_2\,z_1-w_3\,z_7+w_4\,z_7-w_5\,z_6-w_6\,z_6-w_7\,z_5+w_8\,z_5)\,\partial_{w_4}
\\
&
\ \ \ \ \
+(2\,w_2\,z_4+w_3\,z_6-w_4\,z_6+w_5\,z_7+w_6\,z_7-w_7\,z_8+w_8\,z_8)\,\partial_{w_5}
\\
&
\ \ \ \ \
-(2\,w_2\,z_4+w_3\,z_6-w_4\,z_6-w_5\,z_7-w_6\,z_7-w_7\,z_8+w_8\,z_8)\,\partial_{w_6}
\\
&
\ \ \ \ \
+(2\,w_2\,z_3+w_3\,z_5-w_4\,z_5+w_5\,z_8+w_6\,z_8+w_7\,z_7-w_8\,z_7)\,\partial_{w_7}
\\
&
\ \ \ \ \
+(2\,w_2\,z_3+w_3\,z_5-w_4\,z_5+w_5\,z_8+w_6\,z_8-w_7\,z_7+w_8\,z_7)\,\partial_{w_8},
\endaligned
\]
\[
\aligned
L_{z_8z_8}^1
&
\,:=\,
(2\,z_2\,z_7-2\,z_3\,z_6+2\,z_4\,z_5)\,\partial_{z_1}-(-2\,z_2\,z_8+w_3+w_4)\,\partial_{z_2}+(2\,z_3\,z_8+w_5+w_6)\,\partial_{z_3}
\\
&
\ \ \ \ \
-(-2\,z_4\,z_8+w_7-w_8)\,\partial_{z_4}+2\,z_5\,z_8\,\partial_{z_5}+2\,z_6\,z_8\,\partial_{z_6}+2\,z_7\,z_8\,\partial_{z_7}-(-2\,z_8^2+2\,w_2)\,\partial_{z_8}
\\
&
\ \ \ \ \
-(w_3\,z_2-w_4\,z_2-w_5\,z_3+w_6\,z_3+w_7\,z_4+w_8\,z_4)\,\partial_{w_1}+2\,w_2\,z_8\,\partial_{w_2}
\\
&
\ \ \ \ \
+(2\,w_2\,z_2+w_3\,z_8-w_4\,z_8+w_5\,z_5-w_6\,z_5-w_7\,z_6-w_8\,z_6)\,\partial_{w_3}
\\
&
\ \ \ \ \
+(2\,w_2\,z_2-w_3\,z_8+w_4\,z_8+w_5\,z_5-w_6\,z_5-w_7\,z_6-w_8\,z_6)\,\partial_{w_4}
\\
&
\ \ \ \ \
-(2\,w_2\,z_3+w_3\,z_5-w_4\,z_5-w_5\,z_8+w_6\,z_8-w_7\,z_7-w_8\,z_7)\,\partial_{w_5}
\\
&
\ \ \ \ \
-(2\,w_2\,z_3+w_3\,z_5-w_4\,z_5+w_5\,z_8-w_6\,z_8-w_7\,z_7-w_8\,z_7)\,\partial_{w_6}
\\
&
\ \ \ \ \
+(2\,w_2\,z_4+w_3\,z_6-w_4\,z_6-w_5\,z_7+w_6\,z_7+w_7\,z_8+w_8\,z_8)\,\partial_{w_7}
\\
&
\ \ \ \ \
-(2\,w_2\,z_4+w_3\,z_6-w_4\,z_6-w_5\,z_7+w_6\,z_7-w_7\,z_8-w_8\,z_8)\,\partial_{w_8},
\endaligned
\]
\[
\aligned
IL_{z_1z_1}^1
&
\,:=\,
\isqrt\,
\Big(
(2\,z_1^2+2\,w_1)\,\partial_{z_1}+2\,z_1\,z_2\,\partial_{z_2}+2\,z_1\,z_3\,\partial_{z_3}+2\,z_1\,z_4\,\partial_{z_4}\\
&
\ \ \ \ \ \ \ \ \ \ \ \ \
-(-2\,z_1\,z_5+w_7+w_8)\,\partial_{z_5}-(-2\,z_1\,z_6+w_5-w_6)\,\partial_{z_6}-(-2\,z_1\,z_7+w_3-w_4)\,\partial_{z_7}
\\
&
\ \ \ \ \ \ \ \ \ \ \ \ \
+(2\,z_2\,z_7-2\,z_3\,z_6+2\,z_4\,z_5)\,\partial_{z_8}+2\,z_1\,w_1\,\partial_{w_1}
\\
&
\ \ \ \ \ \ \ \ \ \ \ \ \
+(w_3\,z_7+w_4\,z_7+w_5\,z_6+w_6\,z_6+w_7\,z_5-w_8\,z_5)\,\partial_{w_2}
\\
&
\ \ \ \ \ \ \ \ \ \ \ \ \
-(2\,w_1\,z_7-w_3\,z_1-w_4\,z_1+w_5\,z_4+w_6\,z_4+w_7\,z_3-w_8\,z_3)\,\partial_{w_3}
\\
&
\ \ \ \ \ \ \ \ \ \ \ \ \
+(2\,w_1\,z_7+w_3\,z_1+w_4\,z_1+w_5\,z_4+w_6\,z_4+w_7\,z_3-w_8\,z_3)\,\partial_{w_4}
\\
&
\ \ \ \ \ \ \ \ \ \ \ \ \
-(2\,w_1\,z_6-w_3\,z_4-w_4\,z_4-w_5\,z_1-w_6\,z_1+w_7\,z_2-w_8\,z_2)\,\partial_{w_5}
\\
&
\ \ \ \ \ \ \ \ \ \ \ \ \
+(2\,w_1\,z_6-w_3\,z_4-w_4\,z_4+w_5\,z_1+w_6\,z_1+w_7\,z_2-w_8\,z_2)\,\partial_{w_6}
\\
&
\ \ \ \ \ \ \ \ \ \ \ \ \
-(2\,w_1\,z_5-w_3\,z_3-w_4\,z_3-w_5\,z_2-w_6\,z_2-w_7\,z_1+w_8\,z_1)\,\partial_{w_7}
\\
&
\ \ \ \ \ \ \ \ \ \ \ \ \
-(2\,w_1\,z_5-w_3\,z_3-w_4\,z_3-w_5\,z_2-w_6\,z_2+w_7\,z_1-w_8\,z_1)\,\partial_{w_8}
\Big),
\endaligned
\]
\[
\aligned
IL_{z_2z_2}^1
&
\,:=\,
\isqrt\,
\Big(
2\,z_1\,z_2\,\partial_{z_1}+(2\,z_2^2+2\,w_1)\,\partial_{z_2}+2\,z_2\,z_3\,\partial_{z_3}+2\,z_2\,z_4\,\partial_{z_4}+(2\,z_2\,z_5+w_5+w_6)\,\partial_{z_5}
\\
&
\ \ \ \ \ \ \ \ \ \ \ \ \
-(-2\,z_2\,z_6+w_7-w_8)\,\partial_{z_6}+(2\,z_1\,z_8+2\,z_3\,z_6-2\,z_4\,z_5)\,\partial_{z_7}
\\
&
\ \ \ \ \ \ \ \ \ \ \ \ \
-(-2\,z_2\,z_8+w_3-w_4)\,\partial_{z_8}+2\,z_2\,w_1\,\partial_{w_1}
\\
&
\ \ \ \ \ \ \ \ \ \ \ \ \
+(w_3\,z_8+w_4\,z_8-w_5\,z_5+w_6\,z_5+w_7\,z_6+w_8\,z_6)\,\partial_{w_2}
\\
&
\ \ \ \ \ \ \ \ \ \ \ \ \
-(2\,w_1\,z_8-w_3\,z_2-w_4\,z_2-w_5\,z_3+w_6\,z_3+w_7\,z_4+w_8\,z_4)\,\partial_{w_3}
\\
&
\ \ \ \ \ \ \ \ \ \ \ \ \
+(2\,w_1\,z_8+w_3\,z_2+w_4\,z_2-w_5\,z_3+w_6\,z_3+w_7\,z_4+w_8\,z_4)\,\partial_{w_4}
\\
&
\ \ \ \ \ \ \ \ \ \ \ \ \
+(2\,w_1\,z_5-w_3\,z_3-w_4\,z_3+w_5\,z_2-w_6\,z_2+w_7\,z_1+w_8\,z_1)\,\partial_{w_5}
\\
&
\ \ \ \ \ \ \ \ \ \ \ \ \
+(2\,w_1\,z_5-w_3\,z_3-w_4\,z_3-w_5\,z_2+w_6\,z_2+w_7\,z_1+w_8\,z_1)\,\partial_{w_6}
\\
&
\ \ \ \ \ \ \ \ \ \ \ \ \
-(2\,w_1\,z_6-w_3\,z_4-w_4\,z_4+w_5\,z_1-w_6\,z_1-w_7\,z_2-w_8\,z_2)\,\partial_{w_7}
\\
&
\ \ \ \ \ \ \ \ \ \ \ \ \
+(2\,w_1\,z_6-w_3\,z_4-w_4\,z_4+w_5\,z_1-w_6\,z_1+w_7\,z_2+w_8\,z_2)\,\partial_{w_8}
\Big),
\endaligned
\]
\[
\aligned
IL_{z_3z_3}^1
&
\,:=\,
\isqrt\,
\Big(
2\,z_1\,z_3\,\partial_{z_1}+2\,z_2\,z_3\,\partial_{z_2}+(2\,z_3^2+2\,w_1)\,\partial_{z_3}+2\,z_3\,z_4\,\partial_{z_4}+(2\,z_3\,z_5+w_3+w_4)\,\partial_{z_5}
\\
&
\ \ \ \ \ \ \ \ \ \ \ \ \
-(2\,z_1\,z_8-2\,z_2\,z_7-2\,z_4\,z_5)\,\partial_{z_6}-(-2\,z_3\,z_7+w_7-w_8)\,\partial_{z_7}+(2\,z_3\,z_8+w_5-w_6)\,\partial_{z_8}
\\
&
\ \ \ \ \ \ \ \ \ \ \ \ \
+2\,z_3\,w_1\,\partial_{w_1}-(w_3\,z_5-w_4\,z_5+w_5\,z_8+w_6\,z_8-w_7\,z_7-w_8\,z_7)\,\partial_{w_2}
\\
&
\ \ \ \ \ \ \ \ \ \ \ \ \
+(2\,w_1\,z_5+w_3\,z_3-w_4\,z_3-w_5\,z_2-w_6\,z_2+w_7\,z_1+w_8\,z_1)\,\partial_{w_3}
\\
&
\ \ \ \ \ \ \ \ \ \ \ \ \
+(2\,w_1\,z_5-w_3\,z_3+w_4\,z_3-w_5\,z_2-w_6\,z_2+w_7\,z_1+w_8\,z_1)\,\partial_{w_4}
\\
&
\ \ \ \ \ \ \ \ \ \ \ \ \
+(2\,w_1\,z_8+w_3\,z_2-w_4\,z_2+w_5\,z_3+w_6\,z_3+w_7\,z_4+w_8\,z_4)\,\partial_{w_5}
\\
&
\ \ \ \ \ \ \ \ \ \ \ \ \
-(2\,w_1\,z_8+w_3\,z_2-w_4\,z_2-w_5\,z_3-w_6\,z_3+w_7\,z_4+w_8\,z_4)\,\partial_{w_6}
\\
&
\ \ \ \ \ \ \ \ \ \ \ \ \
-(2\,w_1\,z_7+w_3\,z_1-w_4\,z_1+w_5\,z_4+w_6\,z_4-w_7\,z_3-w_8\,z_3)\,\partial_{w_7}
\\
&
\ \ \ \ \ \ \ \ \ \ \ \ \
+(2\,w_1\,z_7+w_3\,z_1-w_4\,z_1+w_5\,z_4+w_6\,z_4+w_7\,z_3+w_8\,z_3)\,\partial_{w_8}
\Big),
\endaligned
\]
\[
\aligned
IL_{z_4z_4}^1
&
\,:=\,
\isqrt\,
\Big(
2\,z_1\,z_4\,\partial_{z_1}+2\,z_2\,z_4\,\partial_{z_2}+2\,z_3\,z_4\,\partial_{z_3}+(2\,z_4^2+2\,w_1)\,\partial_{z_4}+(2\,z_1\,z_8-2\,z_2\,z_7+2\,z_3\,z_6)\,\partial_{z_5}
\\
&
\ \ \ \ \ \ \ \ \ \ \ \ \
+(2\,z_4\,z_6+w_3+w_4)\,\partial_{z_6}-(-2\,z_4\,z_7+w_5+w_6)\,\partial_{z_7}
\\
&
\ \ \ \ \ \ \ \ \ \ \ \ \
-(-2\,z_4\,z_8+w_7+w_8)\,\partial_{z_8}+2\,z_4\,w_1\,\partial_{w_1}
\\
&
\ \ \ \ \ \ \ \ \ \ \ \ \
-(w_3\,z_6-w_4\,z_6-w_5\,z_7+w_6\,z_7-w_7\,z_8+w_8\,z_8)\,\partial_{w_2}
\\
&
\ \ \ \ \ \ \ \ \ \ \ \ \
+(2\,w_1\,z_6+w_3\,z_4-w_4\,z_4+w_5\,z_1-w_6\,z_1+w_7\,z_2-w_8\,z_2)\,\partial_{w_3}
\\
&
\ \ \ \ \ \ \ \ \ \ \ \ \
+(2\,w_1\,z_6-w_3\,z_4+w_4\,z_4+w_5\,z_1-w_6\,z_1+w_7\,z_2-w_8\,z_2)\,\partial_{w_4}
\\
&
\ \ \ \ \ \ \ \ \ \ \ \ \
-(2\,w_1\,z_7+w_3\,z_1-w_4\,z_1-w_5\,z_4+w_6\,z_4+w_7\,z_3-w_8\,z_3)\,\partial_{w_5}
\\
&
\ \ \ \ \ \ \ \ \ \ \ \ \
-(2\,w_1\,z_7+w_3\,z_1-w_4\,z_1+w_5\,z_4-w_6\,z_4+w_7\,z_3-w_8\,z_3)\,\partial_{w_6}
\\
&
\ \ \ \ \ \ \ \ \ \ \ \ \
-(2\,w_1\,z_8+w_3\,z_2-w_4\,z_2-w_5\,z_3+w_6\,z_3-w_7\,z_4+w_8\,z_4)\,\partial_{w_7}
\\
&
\ \ \ \ \ \ \ \ \ \ \ \ \
-(2\,w_1\,z_8+w_3\,z_2-w_4\,z_2-w_5\,z_3+w_6\,z_3+w_7\,z_4-w_8\,z_4)\,\partial_{w_8}
\Big),
\endaligned
\]
\[
\aligned
IL_{z_5z_5}^1
&
\,:=\,
\isqrt\,
\Big(
(2\,z_1\,z_5+w_7-w_8)\,\partial_{z_1}-(-2\,z_2\,z_5+w_5-w_6)\,\partial_{z_2}-(-2\,z_3\,z_5+w_3-w_4)\,\partial_{z_3}
\\
&
\ \ \ \ \ \ \ \ \ \ \ \ \
+(2\,z_1\,z_8-2\,z_2\,z_7+2\,z_3\,z_6)\,\partial_{z_4}+(2\,z_5^2+2\,w_2)\,\partial_{z_5}+2\,z_5\,z_6\,\partial_{z_6}+2\,z_5\,z_7\,\partial_{z_7}+2\,z_5\,z_8\,\partial_{z_8}
\\
&
\ \ \ \ \ \ \ \ \ \ \ \ \
+(w_3\,z_3+w_4\,z_3+w_5\,z_2+w_6\,z_2-w_7\,z_1-w_8\,z_1)\,\partial_{w_1}+2\,z_5\,w_2\,\partial_{w_2}
\\
&
\ \ \ \ \ \ \ \ \ \ \ \ \
-(2\,w_2\,z_3-w_3\,z_5-w_4\,z_5+w_5\,z_8+w_6\,z_8-w_7\,z_7-w_8\,z_7)\,\partial_{w_3}
\\
&
\ \ \ \ \ \ \ \ \ \ \ \ \
+(2\,w_2\,z_3+w_3\,z_5+w_4\,z_5+w_5\,z_8+w_6\,z_8-w_7\,z_7-w_8\,z_7)\,\partial_{w_4}
\\
&
\ \ \ \ \ \ \ \ \ \ \ \ \
-(2\,w_2\,z_2-w_3\,z_8-w_4\,z_8-w_5\,z_5-w_6\,z_5-w_7\,z_6-w_8\,z_6)\,\partial_{w_5}
\\
&
\ \ \ \ \ \ \ \ \ \ \ \ \
+(2\,w_2\,z_2-w_3\,z_8-w_4\,z_8+w_5\,z_5+w_6\,z_5-w_7\,z_6-w_8\,z_6)\,\partial_{w_6}
\\
&
\ \ \ \ \ \ \ \ \ \ \ \ \
+(2\,w_2\,z_1-w_3\,z_7-w_4\,z_7-w_5\,z_6-w_6\,z_6+w_7\,z_5+w_8\,z_5)\,\partial_{w_7}
\\
&
\ \ \ \ \ \ \ \ \ \ \ \ \
-(2\,w_2\,z_1-w_3\,z_7-w_4\,z_7-w_5\,z_6-w_6\,z_6-w_7\,z_5-w_8\,z_5)\,\partial_{w_8}
\Big),
\endaligned
\]
\[
\aligned
IL_{z_6z_6}^1
&
\,:=\,
\isqrt\,
\Big(
(2\,z_1\,z_6+w_5+w_6)\,\partial_{z_1}+(2\,z_2\,z_6+w_7+w_8)\,\partial_{z_2}-(2\,z_1\,z_8-2\,z_2\,z_7-2\,z_4\,z_5)\,\partial_{z_3}
\\
&
\ \ \ \ \ \ \ \ \ \ \ \ \
-(-2\,z_4\,z_6+w_3-w_4)\,\partial_{z_4}+2\,z_5\,z_6\,\partial_{z_5}+(2\,z_6^2+2\,w_2)\,\partial_{z_6}+2\,z_6\,z_7\,\partial_{z_7}+2\,z_6\,z_8\,\partial_{z_8}
\\
&
\ \ \ \ \ \ \ \ \ \ \ \ \
+(w_3\,z_4+w_4\,z_4-w_5\,z_1+w_6\,z_1-w_7\,z_2+w_8\,z_2)\,\partial_{w_1}+2\,z_6\,w_2\,\partial_{w_2}
\\
&
\ \ \ \ \ \ \ \ \ \ \ \ \
-(2\,w_2\,z_4-w_3\,z_6-w_4\,z_6-w_5\,z_7+w_6\,z_7-w_7\,z_8+w_8\,z_8)\,\partial_{w_3}
\\
&
\ \ \ \ \ \ \ \ \ \ \ \ \
+(2\,w_2\,z_4+w_3\,z_6+w_4\,z_6-w_5\,z_7+w_6\,z_7-w_7\,z_8+w_8\,z_8)\,\partial_{w_4}
\\
&
\ \ \ \ \ \ \ \ \ \ \ \ \
+(2\,w_2\,z_1-w_3\,z_7-w_4\,z_7+w_5\,z_6-w_6\,z_6-w_7\,z_5+w_8\,z_5)\,\partial_{w_5}
\\
&
\ \ \ \ \ \ \ \ \ \ \ \ \
+(2\,w_2\,z_1-w_3\,z_7-w_4\,z_7-w_5\,z_6+w_6\,z_6-w_7\,z_5+w_8\,z_5)\,\partial_{w_6}
\\
&
\ \ \ \ \ \ \ \ \ \ \ \ \
+(2\,w_2\,z_2-w_3\,z_8-w_4\,z_8+w_5\,z_5-w_6\,z_5+w_7\,z_6-w_8\,z_6)\,\partial_{w_7}
\\
&
\ \ \ \ \ \ \ \ \ \ \ \ \
+(2\,w_2\,z_2-w_3\,z_8-w_4\,z_8+w_5\,z_5-w_6\,z_5-w_7\,z_6+w_8\,z_6)\,\partial_{w_8}
\Big),
\endaligned
\]
\[
\aligned
IL_{z_7z_7}^1
&
\,:=\,
\isqrt\,
\Big(
(2\,z_1\,z_7+w_3+w_4)\,\partial_{z_1}+(2\,z_1\,z_8+2\,z_3\,z_6-2\,z_4\,z_5)\,\partial_{z_2}+(2\,z_3\,z_7+w_7+w_8)\,\partial_{z_3}
\\
&
\ \ \ \ \ \ \ \ \ \ \ \ \
+(2\,z_4\,z_7+w_5-w_6)\,\partial_{z_4}+2\,z_5\,z_7\,\partial_{z_5}+2\,z_6\,z_7\,\partial_{z_6}+(2\,z_7^2+2\,w_2)\,\partial_{z_7}+2\,z_7\,z_8\,\partial_{z_8}
\\
&
\ \ \ \ \ \ \ \ \ \ \ \ \
-(w_3\,z_1-w_4\,z_1+w_5\,z_4+w_6\,z_4+w_7\,z_3-w_8\,z_3)\,\partial_{w_1}+2\,z_7\,w_2\,\partial_{w_2}
\\
&
\ \ \ \ \ \ \ \ \ \ \ \ \
+(2\,w_2\,z_1+w_3\,z_7-w_4\,z_7-w_5\,z_6-w_6\,z_6-w_7\,z_5+w_8\,z_5)\,\partial_{w_3}
\\
&
\ \ \ \ \ \ \ \ \ \ \ \ \
+(2\,w_2\,z_1-w_3\,z_7+w_4\,z_7-w_5\,z_6-w_6\,z_6-w_7\,z_5+w_8\,z_5)\,\partial_{w_4}
\\
&
\ \ \ \ \ \ \ \ \ \ \ \ \
+(2\,w_2\,z_4+w_3\,z_6-w_4\,z_6+w_5\,z_7+w_6\,z_7-w_7\,z_8+w_8\,z_8)\,\partial_{w_5}
\\
&
\ \ \ \ \ \ \ \ \ \ \ \ \
-(2\,w_2\,z_4+w_3\,z_6-w_4\,z_6-w_5\,z_7-w_6\,z_7-w_7\,z_8+w_8\,z_8)\,\partial_{w_6}
\\
&
\ \ \ \ \ \ \ \ \ \ \ \ \
+(2\,w_2\,z_3+w_3\,z_5-w_4\,z_5+w_5\,z_8+w_6\,z_8+w_7\,z_7-w_8\,z_7)\,\partial_{w_7}
\\
&
\ \ \ \ \ \ \ \ \ \ \ \ \
+(2\,w_2\,z_3+w_3\,z_5-w_4\,z_5+w_5\,z_8+w_6\,z_8-w_7\,z_7+w_8\,z_7)\,\partial_{w_8}
\Big),
\endaligned
\]
\[
\aligned
IL_{z_8z_8}^1
&
\,:=\,
\isqrt\,
\Big(
(2\,z_2\,z_7-2\,z_3\,z_6+2\,z_4\,z_5)\,\partial_{z_1}+(2\,z_2\,z_8+w_3+w_4)\,\partial_{z_2}-(-2\,z_3\,z_8+w_5+w_6)\,\partial_{z_3}
\\
&
\ \ \ \ \ \ \ \ \ \ \ \ \
+(2\,z_4\,z_8+w_7-w_8)\,\partial_{z_4}+2\,z_5\,z_8\,\partial_{z_5}+2\,z_6\,z_8\,\partial_{z_6}+2\,z_7\,z_8\,\partial_{z_7}+(2\,z_8^2+2\,w_2)\,\partial_{z_8}
\\
&
\ \ \ \ \ \ \ \ \ \ \ \ \
-(w_3\,z_2-w_4\,z_2-w_5\,z_3+w_6\,z_3+w_7\,z_4+w_8\,z_4)\,\partial_{w_1}+2\,z_8\,w_2\,\partial_{w_2}
\\
&
\ \ \ \ \ \ \ \ \ \ \ \ \
+(2\,w_2\,z_2+w_3\,z_8-w_4\,z_8+w_5\,z_5-w_6\,z_5-w_7\,z_6-w_8\,z_6)\,\partial_{w_3}
\\
&
\ \ \ \ \ \ \ \ \ \ \ \ \
+(2\,w_2\,z_2-w_3\,z_8+w_4\,z_8+w_5\,z_5-w_6\,z_5-w_7\,z_6-w_8\,z_6)\,\partial_{w_4}
\\
&
\ \ \ \ \ \ \ \ \ \ \ \ \
-(2\,w_2\,z_3+w_3\,z_5-w_4\,z_5-w_5\,z_8+w_6\,z_8-w_7\,z_7-w_8\,z_7)\,\partial_{w_5}
\\
&
\ \ \ \ \ \ \ \ \ \ \ \ \
-(2\,w_2\,z_3+w_3\,z_5-w_4\,z_5+w_5\,z_8-w_6\,z_8-w_7\,z_7-w_8\,z_7)\,\partial_{w_6}
\\
&
\ \ \ \ \ \ \ \ \ \ \ \ \
+(2\,w_2\,z_4+w_3\,z_6-w_4\,z_6-w_5\,z_7+w_6\,z_7+w_7\,z_8+w_8\,z_8)\,\partial_{w_7}
\\
&
\ \ \ \ \ \ \ \ \ \ \ \ \
-(2\,w_2\,z_4+w_3\,z_6-w_4\,z_6-w_5\,z_7+w_6\,z_7-w_7\,z_8-w_8\,z_8)\,\partial_{w_8}
\Big).
\endaligned
\]

The $3 + 5$ generators of $\mathfrak{g}_2$ are:
\[
\aligned
L_{w_3w_3}^2
&
\,:=\,
(2\,w_1\,z_7-w_3\,z_1-w_4\,z_1+w_5\,z_4+w_6\,z_4+w_7\,z_3-w_8\,z_3)\,\partial_{z_1}
\\
&
\ \ \ \ \
+(2\,w_1\,z_8-w_3\,z_2-w_4\,z_2-w_5\,z_3+w_6\,z_3+w_7\,z_4+w_8\,z_4)\,\partial_{z_2}
\\
&
\ \ \ \ \
-(2\,w_1\,z_5+w_3\,z_3-w_4\,z_3-w_5\,z_2-w_6\,z_2+w_7\,z_1+w_8\,z_1)\,\partial_{z_3}
\\
&
\ \ \ \ \
-(2\,w_1\,z_6+w_3\,z_4-w_4\,z_4+w_5\,z_1-w_6\,z_1+w_7\,z_2-w_8\,z_2)\,\partial_{z_4}
\\
&
\ \ \ \ \
+(2\,w_2\,z_3-w_3\,z_5-w_4\,z_5+w_5\,z_8+w_6\,z_8-w_7\,z_7-w_8\,z_7)\,\partial_{z_5}
\\
&
\ \ \ \ \
+(2\,w_2\,z_4-w_3\,z_6-w_4\,z_6-w_5\,z_7+w_6\,z_7-w_7\,z_8+w_8\,z_8)\,\partial_{z_6}
\\
&
\ \ \ \ \
-(2\,w_2\,z_1+w_3\,z_7-w_4\,z_7-w_5\,z_6-w_6\,z_6-w_7\,z_5+w_8\,z_5)\,\partial_{z_7}
\\
&
\ \ \ \ \
-(2\,w_2\,z_2+w_3\,z_8-w_4\,z_8+w_5\,z_5-w_6\,z_5-w_7\,z_6-w_8\,z_6)\,\partial_{z_8}
\\
&
\ \ \ \ \
-2\,w_1\,w_3\,\partial_{w_1}-2\,w_2\,w_3\,\partial_{w_2}+(4\,w_1\,w_2-w_3^2-w_4^2+w_5^2-w_6^2+w_7^2-w_8^2)\,\partial_{w_3}
\\
&
\ \ \ \ \
-2\,w_3\,w_4\,\partial_{w_4}-2\,w_3\,w_5\,\partial_{w_5}-2\,w_3\,w_6\,\partial_{w_6}-2\,w_3\,w_7\,\partial_{w_7}-2\,w_3\,w_8\,\partial_{w_8},
\endaligned
\]
\[
\aligned
L_{w_5w_5}^2
&
\,:=\,
(2\,w_1\,z_6-w_3\,z_4-w_4\,z_4-w_5\,z_1-w_6\,z_1+w_7\,z_2-w_8\,z_2)\,\partial_{z_1}
\\
&
\ \ \ \ \
-(2\,w_1\,z_5-w_3\,z_3-w_4\,z_3+w_5\,z_2-w_6\,z_2+w_7\,z_1+w_8\,z_1)\,\partial_{z_2}
\\
&
\ \ \ \ \
-(2\,w_1\,z_8+w_3\,z_2-w_4\,z_2+w_5\,z_3+w_6\,z_3+w_7\,z_4+w_8\,z_4)\,\partial_{z_3}
\\
&
\ \ \ \ \
+(2\,w_1\,z_7+w_3\,z_1-w_4\,z_1-w_5\,z_4+w_6\,z_4+w_7\,z_3-w_8\,z_3)\,\partial_{z_4}
\\
&
\ \ \ \ \
+(2\,w_2\,z_2-w_3\,z_8-w_4\,z_8-w_5\,z_5-w_6\,z_5-w_7\,z_6-w_8\,z_6)\,\partial_{z_5}
\\
&
\ \ \ \ \
-(2\,w_2\,z_1-w_3\,z_7-w_4\,z_7+w_5\,z_6-w_6\,z_6-w_7\,z_5+w_8\,z_5)\,\partial_{z_6}
\\
&
\ \ \ \ \
-(2\,w_2\,z_4+w_3\,z_6-w_4\,z_6+w_5\,z_7+w_6\,z_7-w_7\,z_8+w_8\,z_8)\,\partial_{z_7}
\\
&
\ \ \ \ \
+(2\,w_2\,z_3+w_3\,z_5-w_4\,z_5-w_5\,z_8+w_6\,z_8-w_7\,z_7-w_8\,z_7)\,\partial_{z_8}
\\
&
\ \ \ \ \
-2\,w_1\,w_5\,\partial_{w_1}-2\,w_2\,w_5\,\partial_{w_2}-2\,w_3\,w_5\,\partial_{w_3}
\\
&
\ \ \ \ \
-2\,w_4\,w_5\,\partial_{w_4}+(4\,w_1\,w_2+w_3^2-w_4^2-w_5^2-w_6^2+w_7^2-w_8^2)\,\partial_{w_5}
\\
&
\ \ \ \ \
-2\,w_5\,w_6\,\partial_{w_6}-2\,w_5\,w_7\,\partial_{w_7}-2\,w_5\,w_8\,\partial_{w_8},
\endaligned
\]
\[
\aligned
L_{w_7w_7}^2
&
\,:=\,
(2\,w_1\,z_5-w_3\,z_3-w_4\,z_3-w_5\,z_2-w_6\,z_2-w_7\,z_1+w_8\,z_1)\,\partial_{z_1}
\\
&
\ \ \ \ \
+(2\,w_1\,z_6-w_3\,z_4-w_4\,z_4+w_5\,z_1-w_6\,z_1-w_7\,z_2-w_8\,z_2)\,\partial_{z_2}
\\
&
\ \ \ \ \
+(2\,w_1\,z_7+w_3\,z_1-w_4\,z_1+w_5\,z_4+w_6\,z_4-w_7\,z_3-w_8\,z_3)\,\partial_{z_3}
\\
&
\ \ \ \ \
+(2\,w_1\,z_8+w_3\,z_2-w_4\,z_2-w_5\,z_3+w_6\,z_3-w_7\,z_4+w_8\,z_4)\,\partial_{z_4}
\\
&
\ \ \ \ \
-(2\,w_2\,z_1-w_3\,z_7-w_4\,z_7-w_5\,z_6-w_6\,z_6+w_7\,z_5+w_8\,z_5)\,\partial_{z_5}
\\
&
\ \ \ \ \
-(2\,w_2\,z_2-w_3\,z_8-w_4\,z_8+w_5\,z_5-w_6\,z_5+w_7\,z_6-w_8\,z_6)\,\partial_{z_6}
\\
&
\ \ \ \ \
-(2\,w_2\,z_3+w_3\,z_5-w_4\,z_5+w_5\,z_8+w_6\,z_8+w_7\,z_7-w_8\,z_7)\,\partial_{z_7}
\\
&
\ \ \ \ \
-(2\,w_2\,z_4+w_3\,z_6-w_4\,z_6-w_5\,z_7+w_6\,z_7+w_7\,z_8+w_8\,z_8)\,\partial_{z_8}
\\
&
\ \ \ \ \
-2\,w_1\,w_7\,\partial_{w_1}-2\,w_2\,w_7\,\partial_{w_2}-2\,w_3\,w_7\,\partial_{w_3}-2\,w_4\,w_7\,\partial_{w_4}-2\,w_5\,w_7\,\partial_{w_5}
\\
&
\ \ \ \ \
-2\,w_6\,w_7\,\partial_{w_6}+(4\,w_1\,w_2+w_3^2-w_4^2+w_5^2-w_6^2-w_7^2-w_8^2)\,\partial_{w_7}-2\,w_7\,w_8\,\partial_{w_8},
\endaligned
\]

\[
\aligned
IL_{w_1w_1}^2
&
\,:=\,
\isqrt\,
\Big(
4\,w_1\,z_1\,\partial_{z_1}+4\,w_1\,z_2\,\partial_{z_2}+4\,w_1\,z_3\,\partial_{z_3}+4\,w_1\,z_4\,\partial_{z_4}
\\
&
\ \ \ \ \ \ \ \ \ \ \ \ \
+(2\,w_3\,z_3+2\,w_4\,z_3+2\,w_5\,z_2+2\,w_6\,z_2-2\,w_7\,z_1-2\,w_8\,z_1)\,\partial_{z_5}
\\
&
\ \ \ \ \ \ \ \ \ \ \ \ \
+(2\,w_3\,z_4+2\,w_4\,z_4-2\,w_5\,z_1+2\,w_6\,z_1-2\,w_7\,z_2+2\,w_8\,z_2)\,\partial_{z_6}
\\
&
\ \ \ \ \ \ \ \ \ \ \ \ \
-(2\,w_3\,z_1-2\,w_4\,z_1+2\,w_5\,z_4+2\,w_6\,z_4+2\,w_7\,z_3-2\,w_8\,z_3)\,\partial_{z_7}
\\
&
\ \ \ \ \ \ \ \ \ \ \ \ \
-(2\,w_3\,z_2-2\,w_4\,z_2-2\,w_5\,z_3+2\,w_6\,z_3+2\,w_7\,z_4+2\,w_8\,z_4)\,\partial_{z_8}
\\
&
\ \ \ \ \ \ \ \ \ \ \ \ \
+4\,w_1^2\,\partial_{w_1}-(w_3^2-w_4^2+w_5^2-w_6^2+w_7^2-w_8^2)\,\partial_{w_2}+4\,w_1\,w_3\,\partial_{w_3}
\\
&
\ \ \ \ \ \ \ \ \ \ \ \ \
+4\,w_1\,w_4\,\partial_{w_4}+4\,w_1\,w_5\,\partial_{w_5}+4\,w_1\,w_6\,\partial_{w_6}+4\,w_1\,w_7\,\partial_{w_7}+4\,w_1\,w_8\,\partial_{w_8}
\Big),
\endaligned
\]
\[
\aligned
IL_{w_2w_2}^2
&
\,:=\,
\isqrt\,
\Big(
(2\,w_3\,z_7+2\,w_4\,z_7+2\,w_5\,z_6+2\,w_6\,z_6+2\,w_7\,z_5-2\,w_8\,z_5)\,\partial_{z_1}
\\
&
\ \ \ \ \ \ \ \ \ \ \ \ \
+(2\,w_3\,z_8+2\,w_4\,z_8-2\,w_5\,z_5+2\,w_6\,z_5+2\,w_7\,z_6+2\,w_8\,z_6)\,\partial_{z_2}
\\
&
\ \ \ \ \ \ \ \ \ \ \ \ \
-(2\,w_3\,z_5-2\,w_4\,z_5+2\,w_5\,z_8+2\,w_6\,z_8-2\,w_7\,z_7-2\,w_8\,z_7)\,\partial_{z_3}
\\
&
\ \ \ \ \ \ \ \ \ \ \ \ \
-(2\,w_3\,z_6-2\,w_4\,z_6-2\,w_5\,z_7+2\,w_6\,z_7-2\,w_7\,z_8+2\,w_8\,z_8)\,\partial_{z_4}+4\,w_2\,z_5\,\partial_{z_5}
\\
&
\ \ \ \ \ \ \ \ \ \ \ \ \
+4\,w_2\,z_6\,\partial_{z_6}+4\,w_2\,z_7\,\partial_{z_7}+4\,w_2\,z_8\,\partial_{z_8}-(w_3^2-w_4^2+w_5^2-w_6^2+w_7^2-w_8^2)\,\partial_{w_1}
\\
&
\ \ \ \ \ \ \ \ \ \ \ \ \
+4\,w_2^2\,\partial_{w_2}+4\,w_2\,w_3\,\partial_{w_3}+4\,w_2\,w_4\,\partial_{w_4}+4\,w_2\,w_5\,\partial_{w_5}+4\,w_2\,w_6\,\partial_{w_6}
\\
&
\ \ \ \ \ \ \ \ \ \ \ \ \
+4\,w_2\,w_7\,\partial_{w_7}+4\,w_2\,w_8\,\partial_{w_8}
\Big),
\endaligned
\]
\[
\aligned
IL_{w_4w_4}^2
&
\,:=\,
\isqrt\,
\Big(
(2\,w_1\,z_7+w_3\,z_1+w_4\,z_1+w_5\,z_4+w_6\,z_4+w_7\,z_3-w_8\,z_3)\,\partial_{z_1}
\\
&
\ \ \ \ \ \ \ \ \ \ \ \ \
+(2\,w_1\,z_8+w_3\,z_2+w_4\,z_2-w_5\,z_3+w_6\,z_3+w_7\,z_4+w_8\,z_4)\,\partial_{z_2}
\\
&
\ \ \ \ \ \ \ \ \ \ \ \ \
+(2\,w_1\,z_5-w_3\,z_3+w_4\,z_3-w_5\,z_2-w_6\,z_2+w_7\,z_1+w_8\,z_1)\,\partial_{z_3}
\\
&
\ \ \ \ \ \ \ \ \ \ \ \ \
+(2\,w_1\,z_6-w_3\,z_4+w_4\,z_4+w_5\,z_1-w_6\,z_1+w_7\,z_2-w_8\,z_2)\,\partial_{z_4}
\\
&
\ \ \ \ \ \ \ \ \ \ \ \ \
+(2\,w_2\,z_3+w_3\,z_5+w_4\,z_5+w_5\,z_8+w_6\,z_8-w_7\,z_7-w_8\,z_7)\,\partial_{z_5}
\\
&
\ \ \ \ \ \ \ \ \ \ \ \ \
+(2\,w_2\,z_4+w_3\,z_6+w_4\,z_6-w_5\,z_7+w_6\,z_7-w_7\,z_8+w_8\,z_8)\,\partial_{z_6}
\\
&
\ \ \ \ \ \ \ \ \ \ \ \ \
+(2\,w_2\,z_1-w_3\,z_7+w_4\,z_7-w_5\,z_6-w_6\,z_6-w_7\,z_5+w_8\,z_5)\,\partial_{z_7}
\\
&
\ \ \ \ \ \ \ \ \ \ \ \ \
+(2\,w_2\,z_2-w_3\,z_8+w_4\,z_8+w_5\,z_5-w_6\,z_5-w_7\,z_6-w_8\,z_6)\,\partial_{z_8}
\\
&
\ \ \ \ \ \ \ \ \ \ \ \ \
+2\,w_1\,w_4\,\partial_{w_1}+2\,w_2\,w_4\,\partial_{w_2}+2\,w_3\,w_4\,\partial_{w_3}
\\
&
\ \ \ \ \ \ \ \ \ \ \ \ \
+(4\,w_1\,w_2+w_3^2+w_4^2+w_5^2-w_6^2+w_7^2-w_8^2)\,\partial_{w_4}+2\,w_4\,w_5\,\partial_{w_5}+2\,w_4\,w_6\,\partial_{w_6}
\\
&
\ \ \ \ \ \ \ \ \ \ \ \ \
+2\,w_4\,w_7\,\partial_{w_7}+2\,w_4\,w_8\,\partial_{w_8}
\Big),
\endaligned
\]
\[
\aligned
IL_{w_6w_6}^2
&
\,:=\,
\isqrt\,
\Big(
(2\,w_1\,z_6-w_3\,z_4-w_4\,z_4+w_5\,z_1+w_6\,z_1+w_7\,z_2-w_8\,z_2)\,\partial_{z_1}
\\
&
\ \ \ \ \ \ \ \ \ \ \ \ \
+(2\,w_1\,z_5-w_3\,z_3-w_4\,z_3-w_5\,z_2+w_6\,z_2+w_7\,z_1+w_8\,z_1)\,\partial_{z_2}
\\
&
\ \ \ \ \ \ \ \ \ \ \ \ \
-(2\,w_1\,z_8+w_3\,z_2-w_4\,z_2-w_5\,z_3-w_6\,z_3+w_7\,z_4+w_8\,z_4)\,\partial_{z_3}
\\
&
\ \ \ \ \ \ \ \ \ \ \ \ \
-(2\,w_1\,z_7+w_3\,z_1-w_4\,z_1+w_5\,z_4-w_6\,z_4+w_7\,z_3-w_8\,z_3)\,\partial_{z_4}
\\
&
\ \ \ \ \ \ \ \ \ \ \ \ \
+(2\,w_2\,z_2-w_3\,z_8-w_4\,z_8+w_5\,z_5+w_6\,z_5-w_7\,z_6-w_8\,z_6)\,\partial_{z_5}
\\
&
\ \ \ \ \ \ \ \ \ \ \ \ \
+(2\,w_2\,z_1-w_3\,z_7-w_4\,z_7-w_5\,z_6+w_6\,z_6-w_7\,z_5+w_8\,z_5)\,\partial_{z_6}
\\
&
\ \ \ \ \ \ \ \ \ \ \ \ \
-(2\,w_2\,z_4+w_3\,z_6-w_4\,z_6-w_5\,z_7-w_6\,z_7-w_7\,z_8+w_8\,z_8)\,\partial_{z_7}
\\
&
\ \ \ \ \ \ \ \ \ \ \ \ \
-(2\,w_2\,z_3+w_3\,z_5-w_4\,z_5+w_5\,z_8-w_6\,z_8-w_7\,z_7-w_8\,z_7)\,\partial_{z_8}
\\
&
\ \ \ \ \ \ \ \ \ \ \ \ \
+2\,w_1\,w_6\,\partial_{w_1}+2\,w_2\,w_6\,\partial_{w_2}+2\,w_3\,w_6\,\partial_{w_3}+2\,w_4\,w_6\,\partial_{w_4}+2\,w_5\,w_6\,\partial_{w_5}
\\
&
\ \ \ \ \ \ \ \ \ \ \ \ \
+(4\,w_1\,w_2+w_3^2-w_4^2+w_5^2+w_6^2+w_7^2-w_8^2)\,\partial_{w_6}+2\,w_6\,w_7\,\partial_{w_7}+2\,w_6\,w_8\,\partial_{w_8}
\Big),
\endaligned
\]
\[
\aligned
IL_{w_8w_8}^2
&
\,:=\,
\isqrt\,
\Big(
(2\,w_1\,z_5-w_3\,z_3-w_4\,z_3-w_5\,z_2-w_6\,z_2+w_7\,z_1-w_8\,z_1)\,\partial_{z_1}
\\
&
\ \ \ \ \ \ \ \ \ \ \ \ \
-(2\,w_1\,z_6-w_3\,z_4-w_4\,z_4+w_5\,z_1-w_6\,z_1+w_7\,z_2+w_8\,z_2)\,\partial_{z_2}
\\
&
\ \ \ \ \ \ \ \ \ \ \ \ \
-(2\,w_1\,z_7+w_3\,z_1-w_4\,z_1+w_5\,z_4+w_6\,z_4+w_7\,z_3+w_8\,z_3)\,\partial_{z_3}
\\
&
\ \ \ \ \ \ \ \ \ \ \ \ \
+(2\,w_1\,z_8+w_3\,z_2-w_4\,z_2-w_5\,z_3+w_6\,z_3+w_7\,z_4-w_8\,z_4)\,\partial_{z_4}
\\
&
\ \ \ \ \ \ \ \ \ \ \ \ \
+(2\,w_2\,z_1-w_3\,z_7-w_4\,z_7-w_5\,z_6-w_6\,z_6-w_7\,z_5-w_8\,z_5)\,\partial_{z_5}
\\
&
\ \ \ \ \ \ \ \ \ \ \ \ \
-(2\,w_2\,z_2-w_3\,z_8-w_4\,z_8+w_5\,z_5-w_6\,z_5-w_7\,z_6+w_8\,z_6)\,\partial_{z_6}
\\
&
\ \ \ \ \ \ \ \ \ \ \ \ \
-(2\,w_2\,z_3+w_3\,z_5-w_4\,z_5+w_5\,z_8+w_6\,z_8-w_7\,z_7+w_8\,z_7)\,\partial_{z_7}
\\
&
\ \ \ \ \ \ \ \ \ \ \ \ \
+(2\,w_2\,z_4+w_3\,z_6-w_4\,z_6-w_5\,z_7+w_6\,z_7-w_7\,z_8-w_8\,z_8)\,\partial_{z_8}
\\
&
\ \ \ \ \ \ \ \ \ \ \ \ \
-2\,w_1\,w_8\,\partial_{w_1}-2\,w_2\,w_8\,\partial_{w_2}-2\,w_3\,w_8\,\partial_{w_3}-2\,w_4\,w_8\,\partial_{w_4}-2\,w_5\,w_8\,\partial_{w_5}
\\
&
\ \ \ \ \ \ \ \ \ \ \ \ \
-2\,w_6\,w_8\,\partial_{w_6}-2\,w_7\,w_8\,\partial_{w_7}-(4\,w_1\,w_2+w_3^2-w_4^2+w_5^2-w_6^2+w_7^2+w_8^2)\,\partial_{w_8}
\Big).
\endaligned
\]

\section{$\mathfrak{so}(n+2,n)$ CR Models}
\label{so-n-2-n-CR-models}

In coordinates:
\[
(z_1,\dots,z_n)
\ \ \ \ \ \ \ \ \ \ \ \ \ \ \ \ \ \ \ \
\text{and}
\ \ \ \ \ \ \ \ \ \ \ \ \ \ \ \ \ \ \ \
\big(
w_{ij}
\big)_{1\leqslant i<j\leqslant n},
\]
consider:
\[
w_{ij}
-
\overline{w}_{ij}
\,=\,
z_i\,\overline{z}_j
-
\overline{z}_i\,
z_j
\eqno
{\scriptstyle{(1\,\leqslant\,i\,<\,j\,\leqslant\,n)}}.
\]
To simplify, we use $n$ instead of $\ell-1$, hence:
\[
{\rm CRdim}\,M
\,=\,
n
\ \ \ \ \ \ \ \ \ \ \ \ \ \ \ \ \ \ \ \
\text{and}
\ \ \ \ \ \ \ \ \ \ \ \ \ \ \ \ \ \ \ \
{\rm codim}\,M
\,=\,
\tfrac{n(n-1)}{2}.
\]

Expected dimensions are:
\[
\def\arraystretch{0.85}
\begin{array}{ccccc}
\,\,\,\,\,\mathfrak{g}_{-2}\,\,\,\,\, & 
\,\,\,\,\,\mathfrak{g}_{-1}\,\,\,\,\, & 
\,\,\,\,\,\mathfrak{g}_0\,\,\,\,\, & 
\,\,\,\,\,\mathfrak{g}_1\,\,\,\,\, & 
\,\,\,\,\,\mathfrak{g}_2\,\,\,\,\,
\\
\frac{n(n-1)}{2} &
2\,n &
n^2+1 &
2\,n &
\frac{n(n-1)}{2} 
\end{array}
\]

\medskip

The $\frac{n(n-1)}{2}$ generators of $\mathfrak{g}_{-2}$ are:
\[
L_{w_{ij}}^{-2}
\,:=\,
\partial_{w_{ij}}
\eqno
{\scriptstyle{(1\,\leqslant\,i\,<\,j\,\leqslant\,n)}}.
\]

The $n + n$ generators of $\mathfrak{g}_{-1}$ are:
\[
\aligned
L_{z_i}^{-1}
&
\,:=\,
\partial_{z_i}
+
\sum_{1\leqslant k<i}\,
z_k\,\partial_{w_{ki}}
-
\sum_{i<k\leqslant n}\,
z_k\,
\partial_{w_{ik}}
&
\ \ \ \ \ \ \ \ \ \ \ \ \ \ \ \ \ \ \ \
&
{\scriptstyle{(1\,\leqslant\,i\,\leqslant\,n)}},
\\
IL_{z_i}^{-1}
&
\,:=\,
\isqrt\,\partial_{z_i}
-
\sum_{1\leqslant k<i}\,
\isqrt\,z_k\,\partial_{w_{ki}}
+
\sum_{i<k\leqslant n}\,
\isqrt\,z_k\,
\partial_{w_{ik}}
&
\ \ \ \ \ \ \ \ \ \ \ \ \ \ \ \ \ \ \ \
&
{\scriptstyle{(1\,\leqslant\,i\,\leqslant\,n)}},
\endaligned
\]

The $\frac{n(n-1)}{2} + n + \frac{n(n-1)}{2} + 1$ 
generators of $\mathfrak{g}_0$ are:
\[
\aligned
L_{ij}^0
&
\,:=\,
z_i\,\partial_{z_j}
+
\sum_{1\leqslant k<i}\,
w_{ki}\,
\partial_{w_{kj}}
-
\sum_{i<k<j}\,
w_{ik}\,
\partial_{w_{kj}}
+
\sum_{j<k\leqslant n}\,
w_{ik}\,
\partial_{w_{jk}}
&
\ \ \ \ \ \ \ \ \ \ \ \ \ \ \ \ \ \ \ \
&
{\scriptstyle{(1\,\leqslant\,i\,<\,j\,\leqslant\,n)}},
\\
L_{ii}^0
&
\,:=\,
z_i\,\partial_{z_i}
+ 
\sum_{1\leqslant k<i}\,
w_{ki}\,
\partial_{w_{ki}}
+
\sum_{i<k\leqslant n}\,
w_{ik}\,
\partial_{w_{ik}}
&
\ \ \ \ \ \ \ \ \ \ \ \ \ \ \ \ \ \ \ \
&
{\scriptstyle{(1\,\leqslant\,i\,\leqslant\,n)}},
\\
L_{ij}^0
&
\,:=\,
z_i\,\partial_{z_j}
+
\sum_{1\leqslant k<j}\,
w_{ki}\,
\partial_{w_{kj}}
-
\sum_{j<k<i}\,
w_{ki}\,
\partial_{w_{jk}}
+
\sum_{i<k\leqslant n}\,
w_{ik}\,
\partial_{w_{jk}},
&
\ \ \ \ \ \ \ \ \ \ \ \ \ \ \ \ \ \ \ \
&
{\scriptstyle{(1\,\leqslant\,j\,<\,i\,\leqslant\,n)}},
\endaligned
\]
together with the rotation:
\[
R
\,:=\,
\isqrt\,z_1\,\partial_{z_1}
+\cdots+
\isqrt\,z_n\,\partial_{z_n}.
\]
The sum of the $L_{ii}^0$
equals the dilation:
\[
D
\,:=\,
\sum_{1\leqslant k\leqslant n}\,
z_k\,\partial_{z_k}
+
\sum_{1\leqslant k<\ell\leqslant n}\,
w_{k\ell}\,
\partial_{w_{k\ell}}.
\]

The $n + n$ generators of $\mathfrak{g}_1$ are:
\[
\aligned
L_{z_iz_i}^1
&
\,:=\,
\sum_{1\leqslant k<i}\,
\big(z_i\,z_k-w_{ki}\big)\,
\partial_{z_k}
+
z_i\,z_i\,\partial_{z_i}
+
\sum_{i<k\leqslant n}\,
\big(
z_i\,z_k
+
w_{ik}
\big)\,
\partial_{z_k}
\\
&
\ \ \ \ \ 
+
\sum_{1\leqslant k<i}\,
z_i\,w_{ki}\,\partial_{w_{ki}}
+
\sum_{i<\ell\leqslant n}\,
z_i\,w_{i\ell}\,\partial_{w_{i\ell}}
\\
&
\ \ \ \ \
+
\sum_{1\leqslant k<\ell<i}\,
\left\vert\!
\begin{smallmatrix}
z_k & -w_{ki}
\\
z_\ell & -w_{\ell i}
\end{smallmatrix}
\!\right\vert\,
\partial_{w_{k\ell}}
+
\sum_{1\leqslant k<i
\atop
i<\ell\leqslant n}\,
\left\vert\!
\begin{smallmatrix}
z_k & -w_{ki}
\\
z_\ell & w_{i\ell}
\end{smallmatrix}
\!\right\vert\,
\partial_{w_{k\ell}}
+
\sum_{i<k<\ell\leqslant n}\,
\left\vert\!
\begin{smallmatrix}
z_k & w_{ik}
\\
z_\ell & w_{i\ell}
\end{smallmatrix}
\!\right\vert\,
\partial_{w_{k\ell}},
\\
IL_{z_iz_i}^1
&
\,:=\,
\sum_{1\leqslant k<i}\,
\big(\isqrt\,z_i\,z_k+\isqrt\,w_{ki}\big)\,
\partial_{z_k}
+
\isqrt\,z_i\,z_i\,\partial_{z_i}
+
\sum_{i<k\leqslant n}\,
\big(
\isqrt\,z_i\,z_k
-
\isqrt\,w_{ik}
\big)\,
\partial_{z_k}
\\
&
\ \ \ \ \ 
+
\sum_{1\leqslant k<i}\,
\isqrt\,z_i\,w_{ki}\,\partial_{w_{ki}}
+
\sum_{i<\ell\leqslant n}\,
\isqrt\,z_i\,w_{i\ell}\,\partial_{w_{i\ell}}
\\
&
\ \ \ \ \
+
\sum_{1\leqslant k<\ell<i}\,
\isqrt\,\left\vert\!
\begin{smallmatrix}
z_k & -w_{ki}
\\
z_\ell & -w_{\ell i}
\end{smallmatrix}
\!\right\vert\,
\partial_{w_{k\ell}}
+
\sum_{1\leqslant k<i
\atop
i<\ell\leqslant n}\,
\isqrt\,\left\vert\!
\begin{smallmatrix}
z_k & -w_{ki}
\\
z_\ell & w_{i\ell}
\end{smallmatrix}
\!\right\vert\,
\partial_{w_{k\ell}}
+
\sum_{i<k<\ell\leqslant n}\,
\isqrt\,\left\vert\!
\begin{smallmatrix}
z_k & w_{ik}
\\
z_\ell & w_{i\ell}
\end{smallmatrix}
\!\right\vert\,
\partial_{w_{k\ell}}.
\endaligned
\]

The $\frac{n(n-1)}{2}$ generators of $\mathfrak{g}_2$ are, 
with $i < j$:
\[
\aligned
L_{w_{ij}w_{ij}}^2
&
\,:=\,
\sum_{1\leqslant k<i}\,
\left\vert\!
\begin{smallmatrix}
z_i & w_{ki}
\\
z_j & w_{kj}
\end{smallmatrix}
\!\right\vert\,
\partial_{z_k}
+
z_i\,w_{ij}\,\partial_{z_i}
+
\sum_{i<k<j}\,
\left\vert\!
\begin{smallmatrix}
z_i & -w_{ik}
\\
z_j & w_{kj}
\end{smallmatrix}
\!\right\vert\,
\partial_{z_k}
+
z_j\,w_{ij}\,\partial_{z_j}
+
\sum_{j<k\leqslant n}\,
\left\vert\!
\begin{smallmatrix}
z_i & -w_{ik}
\\
z_j & -w_{jk}
\end{smallmatrix}
\!\right\vert\,
\partial_{z_k}
\\
&
\ \ \ \ \
+
\sum_{1\leqslant k<i}\,
w_{ij}\,w_{ki}\,
\partial_{w_{ki}}
+
\sum_{1\leqslant k<i}\,
w_{ij}\,w_{kj}\,
\partial_{w_{kj}}
+
\sum_{i<\ell<j}\,
w_{ij}\,w_{i\ell}\,
\partial_{w_{i\ell}}
+
w_{ij}\,w_{ij}\,
\partial_{w_{ij}}
\\
&
\ \ \ \ \
+
\sum_{j<\ell\leqslant n}\,
w_{ij}\,w_{i\ell}\,
\partial_{w_{i\ell}}
+
\sum_{i<k<j}\,
w_{ij}\,w_{kj}\,
\partial_{w_{kj}}
+
\sum_{j<\ell\leqslant n}\,
w_{ij}\,w_{j\ell}\,
\partial_{w_{j\ell}}
\\
&
\ \ \ \ \
+
\sum_{1\leqslant k<\ell<i}\,
\left\vert\!
\begin{smallmatrix}
-w_{ki} & -w_{kj}
\\
-w_{\ell i} & -w_{\ell j}
\end{smallmatrix}
\!\right\vert\,
\partial_{w_{k\ell}}
+
\sum_{1\leqslant k<i
\atop
i<\ell<j}\,
\left\vert\!
\begin{smallmatrix}
-w_{ki} & -w_{kj}
\\
w_{i\ell} & -w_{\ell j}
\end{smallmatrix}
\!\right\vert\,
\partial_{w_{k\ell}}
+
\sum_{1\leqslant k<i
\atop
j<\ell\leqslant n}\,
\left\vert\!
\begin{smallmatrix}
-w_{ki} & -w_{kj}
\\
w_{i\ell} & w_{j\ell}
\end{smallmatrix}
\!\right\vert\,
\partial_{w_{k\ell}}
\\
&
\ \ \ \ \ \ \ \ \ \ \ \ \ \ \ \ \ \ \ \ \ \ \ \ \ \ \ \ \ \ \ \ \ \ \
\ \ \ \ \ \ \ \ \ \ \ \ \ \
+
\sum_{i<k<\ell<j}\,
\left\vert\!
\begin{smallmatrix}
w_{ik} & -w_{kj}
\\
w_{i\ell} & -w_{\ell j}
\end{smallmatrix}
\!\right\vert\,
\partial_{w_{k\ell}}
+
\sum_{i<k<j
\atop
j<\ell\leqslant n}\,
\left\vert\!
\begin{smallmatrix}
w_{ik} & -w_{kj}
\\
w_{i\ell} & w_{j\ell}
\end{smallmatrix}
\!\right\vert\,
\partial_{w_{k\ell}}
\\
&
\ \ \ \ \ \ \ \ \ \ \ \ \ \ \ \ \ \ \ \ \ \ \ \ \ \ \ \ \ \ \ \ \ \ \
\ \ \ \ \ \ \ \ \ \ \ \ \ \ \ \ \ \ \ \ \ \ \ \ \ \ \ \ \ \ \ \ \ \ \
\ \ \ \ \ \ \ \ \ \ \ \ \ \ \ \ 
+
\sum_{j<k<\ell\leqslant n}\,
\left\vert\!
\begin{smallmatrix}
w_{ik} & w_{jk}
\\
w_{i\ell} & w_{j\ell}
\end{smallmatrix}
\!\right\vert\,
\partial_{w_{k\ell}}.
\endaligned
\]

\section{$\mathfrak{su}(p,q)$ CR Models: Notation Set Up}
\label{su-p-q-models-notation-set-up}

Start from (2.2) of Zhaohu's file:
for $1 \leqslant a < c \leqslant \ell$:
\[
\Im\,w_{ac}
\,=\,
\sum_{b=\ell+1}^q\,
\Im\,
\big(
z_{ab}\,
\overline{z}_{c,p+q+1-b}
\big)
+
\sum_{b=q+1}^p\,
\Im\,
\big(
z_{ab}\,\overline{z}_{cb}
\big)
+
\sum_{b=p+1}^{p+q-\ell}\,
\Im\,
\big(
z_{ab}\,
\overline{z}_{c,p+q+1-b}
\big),
\]
for $a = c$:
\[
\Im\,w_{aa}
\,=\,
1\sum_{b=\ell+1}^q\,
\Re\,
\big(
z_{ab}\,\overline{z}_{a,p+q+1-b}
\big)
+
\sum_{b=q+1}^p\,
z_{ab}\,\overline{z}_{ab}
+
1\sum_{b=p+1}^{p+q-\ell}\,
\Re\,
\big(
z_{ab}\,
\overline{z}_{a,p+q+1-b}
\big),
\]
and for $1 \leqslant c < a \leqslant \ell$:
\[
\Im\,w_{ac}
\,=\,
\sum_{b=\ell+1}^q\,
\Re\,
\big(
z_{ab}\,
\overline{z}_{c,p+q+1-b}
\big)
+
\sum_{b=q+1}^p\,
\Re\,
\big(
z_{ab}\,\overline{z}_{cb}
\big)
+
\sum_{b=p+1}^{p+q-\ell}\,
\Re\,
\big(
z_{ab}\,
\overline{z}_{c,p+q+1-b}
\big).
\]
When $c = a$, this third line equals the second one.

Change notation:
\[
\begin{array}{ccccccccccc}
z_{\ell+1} & \dots & z_q 
& &
z_{q+1} & \dots & z_p 
& &
z_{p+1} & \dots & z_{p+q-\ell}
\\
\downarrow & \vdots & \downarrow
& &
\downarrow & \vdots & \downarrow
& &
\downarrow & \vdots & \downarrow
\\
u_1 & \dots & u_{q-\ell} 
& &
z_1 & \dots & z_{p-q}
& &
v_{q-\ell} & \dots & v_1
\end{array}
\]
Also, instead of $\Im\, w_{ac}$ when $a \geqslant c$, 
use $\Re\, w_{ac}$, simply by $w_{ac} \longmapsto 
\isqrt\, w_{ac}$.
An advantage will be seen later. 

Thus for $1 \leqslant a < c \leqslant \ell$:
\[
\Im\,w_{ac}
\,=\,
\Im\,
\left\{
\aligned
u_{a,1}\,\overline{v}_{c,1}
&
+\cdots+
u_{a,q-\ell}\,\overline{v}_{c,q-\ell}
\\
+
z_{a,1}\,\overline{z}_{c,1}
&
+\cdots+
z_{a,p-q}\,\overline{z}_{c,p-q}
\\
+
v_{a,1}\,\overline{u}_{c,1}
&
+\cdots+
v_{a,q-\ell}\,\overline{u}_{c,q-\ell}
\endaligned
\right\}
\]
and for $1 \leqslant c \leqslant a \leqslant \ell$:
\[
\Re\,w_{ac}
\,=\,
\Re\,
\left\{
\aligned
u_{a,1}\,\overline{v}_{c,1}
&
+\cdots+
u_{a,q-\ell}\,\overline{v}_{c,q-\ell}
\\
+
z_{a,1}\,\overline{z}_{c,1}
&
+\cdots+
z_{a,p-q}\,\overline{z}_{c,p-q}
\\
+
v_{a,1}\,\overline{u}_{c,1}
&
+\cdots+
v_{a,q-\ell}\,\overline{u}_{c,q-\ell}
\endaligned
\right\}
\,=\,
\Re\,\big\{\text{same}\big\}.
\]

Change integers:
\[
\begin{array}{ccc}
\ell
&
\longleftrightarrow
&
\ell
\\
m\,:=\,q-\ell
&
\longleftrightarrow
&
q
\\
n\,:=\,p-q
&
\longleftrightarrow
&
p.
\end{array}
\]
Zhaohu's constraints:
\[
\aligned
{}
&
p\,>\,q,
&
\ \ \ \ \
&
1\,\leqslant\,\ell\,\leqslant\,q
&
\ \ \ \ \ \ 
&
\longrightarrow
&
\ \ \ \ \ \ 
&
n\,>\,0
&
\ \
\text{arbitrary}
&
\ \ \ \ \
&
m\,\geqslant\,0
&
\ \
\text{arbitrary},
\\
{}
&
p\,=\,q,
&
\ \ \ \ \
&
1\,\leqslant\,\ell
\underset{=\,\,q-1}{\,\leqslant\,p-1}
&
\ \ \ \ \ \ 
&
\longrightarrow
&
\ \ \ \ \ \ 
&
n\,=\,0
&
\ \ \ \ \
&
\ \ \ \ \
&
m\,>\,0
&
\ \
\text{arbitrary}.
\endaligned
\]
Correspondence:
\[
\aligned
\ell
&
\,=\,
\ell,
\\
q
&
\,=\,
\ell+m,
\\
p
&
\,=\,
\ell+m+n,
\endaligned
\]
with $(m,n) \neq (0,0)$.

\medskip

{\bf Summary.}
Coordinates:
\[
\begin{array}{ccc}
z_{a,1} & \dots & z_{a,n}
\\
u_{a,1} & \dots & u_{a,m}
\\
v_{a,1} & \dots & v_{a,m}
\end{array}
\ \ \ \ \ \ \ \ \ \ \ \ \ \ \ \ \ \ \ \
\begin{array}{ccc}
w_{11} & \dots & w_{1\ell}
\\
\vdots & \ddots & \vdots
\\
w_{\ell1} & \dots & w_{\ell\ell}
\end{array}
\ \ \ \ \ \ \ \ \ \ \ \ \ \ \ \ \ \ \ \
\ell\,\geqslant\,1,
\ \ \ \ \
(n,m)
\,\neq\,
(0,0).
\]
Cauchy-Riemann dimension, real codimension, real dimension:
\[
{\rm CRdim}
\,=\,
(n+m+m)\,\ell,
\ \ \ \ \ \ \ \ \ \ \ \ \ \ \ \ \ \ \ \
{\rm codimR}
\,=\,
\ell^2,
\ \ \ \ \ \ \ \ \ \ \ \ \ \ \ \ \ \ \ \
{\rm dimR}
\,=\,
2\,(n+m+m)\,\ell
+
\ell^2.
\]
Equations for $1 \leqslant a < c \leqslant \ell$:
\[
\Im\,w_{ac}
\,=\,
\Im\,
\left\{
\aligned
z_{a,1}\,\overline{z}_{c,1}
&
+\cdots+
z_{a,n}\,\overline{z}_{n,c}
\\
+
u_{a,1}\,\overline{v}_{c,1}
&
+\cdots+
u_{a,m}\,\overline{v}_{c,m}
\\
+
v_{a,1}\,\overline{u}_{c,1}
&
+\cdots+
v_{a,m}\,\overline{u}_{c,m}
\endaligned
\right\},
\]
and for $1 \leqslant c \leqslant a \leqslant \ell$:
\[
\Re\,w_{ac}
\,=\,
\Re\,
\left\{
\aligned
z_{a,1}\,\overline{z}_{c,1}
&
+\cdots+
z_{a,n}\,\overline{z}_{n,c}
\\
+
u_{a,1}\,\overline{v}_{c,1}
&
+\cdots+
u_{a,m}\,\overline{v}_{c,m}
\\
+
v_{a,1}\,\overline{u}_{c,1}
&
+\cdots+
v_{a,m}\,\overline{u}_{c,m}
\endaligned
\right\}.
\]
Formalism with sums:
\[
\aligned
\Im\,w_{ac}
&
\,=\,
\Im\,
\Big(
\sum_{i=1}^n\,
z_{ai}\,\overline{z}_{ci}
+
\sum_{j=1}^m\,
u_{aj}\,\overline{v}_{cj}
+
\sum_{j=1}^m\,
v_{aj}\,\overline{u}_{cj}
\Big)
&
\ \ \ \ \ 
{\scriptstyle{(a\,<\,c)}},
\\
\Re\,w_{ac}
&
\,=\,
\Re\,
\Big(
\sum_{i=1}^n\,
z_{ai}\,\overline{z}_{ci}
+
\sum_{j=1}^m\,
u_{aj}\,\overline{v}_{cj}
+
\sum_{j=1}^m\,
v_{aj}\,\overline{u}_{cj}
\Big)
\,=\,
\Re\,({\rm same})
&
\ \ \ \ \ 
{\scriptstyle{(a\,\geqslant\,c)}}.
\endaligned
\]

The interest of having done above $w_{ac} \longmapsto 
\isqrt\, w_{ac}$ for $a \geqslant c$
is that $\isqrt$ disappears:
\[
\aligned
\Im\,w_{ac}
&
\,=\,
\Im\,
Term_{ac}
&
\ \ \ \ \ 
&
\Longleftrightarrow
&
\ \ \ \ \ 
w_{ac}
-
\overline{w}_{ac}
&
\,=\,
Term_{ac}
-
\overline{Term}_{ac},
&
\ \ \ \ \ 
{\scriptstyle{(a\,<\,c)}},
\\
\Re\,w_{ac}
&
\,=\,
\Re\,
Term_{ac}
&
\ \ \ \ \ 
&
\Longleftrightarrow
&
\ \ \ \ \ 
w_{ac}
+
\overline{w}_{ac}
&
\,=\,
Term_{ac}
+
\overline{Term}_{ac},
&
\ \ \ \ \ 
{\scriptstyle{(a\,\geqslant\,c)}}.
\endaligned
\]

\section{$\mathfrak{su}(p,q)$ CR Models: Expected dimensions}
\label{su-p-q-models-expected-dimensions}

From Zhaohu's construction:
\[
\aligned
\dim\,\mathfrak{g}_{-2}
&
\,=\,
\ell^2
\,=\,
{\rm codim}\,M,
\\
\dim\,\mathfrak{g}_{-1}
&
\,=\,
2\,(n+m+m)\,\ell
\,=\,
{\rm CRdim}\,M,
\\
\dim\,\mathfrak{g}_{0}
&
\,=\,
\dim\,\mathfrak{su}(p,q)
-
2\,{\rm codim}\,M
-
2\,{\rm CRdim}\,M,
\\
\dim\,\mathfrak{g}_{1}
&
\,=\,
\dim\,\mathfrak{g}_{-1}
\\
\dim\,\mathfrak{g}_{2}
&
\,=\,
\dim\,\mathfrak{g}_{-2},
\endaligned
\]
hence:
\[
\aligned
\dim\,\mathfrak{g}_0
&
\,=\,
\dim\,\mathfrak{su}(p,q)
-
2\,codimR
-
2\,CRdim
\\
&
\,=\,
(p+q)^2-1
-
2\,\ell^2
-
4\,(n+m+m)\,\ell
\\
&
\,=\,
(n+m+m)^2
+
2\,\ell^2-1.
\endaligned
\]

\section{$\mathfrak{su}(p,q)$ CR Models: 
Generators of $\mathfrak{g}_{-2}$ and of $\mathfrak{g}_{-1}$}
\label{generators-g-2-g-1}

Simply:
\[
\mathfrak{g}_{-2}
\,=\,
{\rm Span}\,
\big\{
\partial_{w_{\underset{<}{ac}}},\,\,
\isqrt\,\partial_{w_{\underset{\geqslant}{ac}}}
\big\}
\]

Next, $\mathfrak{g}_{-1}$ is generated by:
\[
\aligned
L_{z_{ij}}
\,:=\,
\partial_{z_{ij}}
&
+
\sum_{1\leqslant r<i}\,
z_{rj}\,
\partial_{w_{ri}}
+
\sum_{1\leqslant s<i}\,
z_{sj}\,
\partial_{w_{is}}
+
2\,z_{ij}\,
\partial_{w_{ii}}
\\
{}
&
+
\sum_{i<r\leqslant\ell}\,
z_{rj}\,
\partial_{w_{ri}}
-
\sum_{i<s\leqslant\ell}\,
z_{sj}\,
\partial_{w_{is}}
\ \ \ \ \ \ \ \ \ \ \ \ \ \ \ \ \ \ \ \
{\scriptstyle{(1\,\leqslant\,i\,\leqslant\,\ell,\,\,
1\leqslant\,j\,\leqslant\,n)}},
\endaligned
\]
\[
\aligned
IL_{z_{ij}}
\,=\,
\isqrt\,
\bigg(
\partial_{z_{ij}}
&
-
\sum_{1\leqslant r<i}\,
z_{rj}\,
\partial_{w_{ri}}
-
\sum_{1\leqslant s<i}\,
z_{sj}\,
\partial_{w_{is}}
-
2\,z_{ij}\,
\partial_{w_{ii}}
\\
{}
&
-
\sum_{i<r\leqslant\ell}\,
z_{rj}\,
\partial_{w_{ri}}
+
\sum_{i<s\leqslant\ell}\,
z_{sj}\,
\partial_{w_{is}}
\bigg)
\ \ \ \ \ \ \ \ \ \ \ \ \ \ \ \ \ \ \ \
{\scriptstyle{(1\,\leqslant\,i\,\leqslant\,\ell,\,\,
1\leqslant\,j\,\leqslant\,n)}},
\endaligned
\]
\[
\aligned
L_{u_{ik}}
\,=\,
\partial_{u_{ik}}
&
+
\sum_{1\leqslant r<i}\,
v_{rk}\,
\partial_{w_{ri}}
+
\sum_{1\leqslant s<i}\,
v_{sk}\,
\partial_{w_{is}}
+
2\,v_{ik}\,\partial_{w_{ii}}
\\
&
+
\sum_{i<r\leqslant\ell}\,
v_{rk}\,
\partial_{w_{ri}}
-
\sum_{i<s\leqslant\ell}\,
v_{sk}\,
\partial_{w_{is}}
\ \ \ \ \ \ \ \ \ \ \ \ \ \ \ \ \ \ \ \
{\scriptstyle{(1\,\leqslant\,i\,\leqslant\,\ell,\,\,
1\leqslant\,k\,\leqslant\,m)}},
\endaligned
\]
\[
\aligned
IL_{u_{ik}}
\,=\,
\isqrt\,
\bigg(
\partial_{u_{ik}}
&
-
\sum_{1\leqslant r<i}\,
v_{rk}\,
\partial_{w_{ri}}
-
\sum_{1\leqslant s<i}\,
v_{sk}\,
\partial_{w_{is}}
-
2\,v_{ik}\,\partial_{w_{ii}}
\\
&
-
\sum_{i<r\leqslant\ell}\,
v_{rk}\,
\partial_{w_{ri}}
+
\sum_{i<s\leqslant\ell}\,
v_{sk}\,
\partial_{w_{is}}
\bigg)
\ \ \ \ \ \ \ \ \ \ \ \ \ \ \ \ \ \ \ \
{\scriptstyle{(1\,\leqslant\,i\,\leqslant\,\ell,\,\,
1\leqslant\,k\,\leqslant\,m)}},
\endaligned
\]
\[
\aligned
L_{v_{ik}}
\,=\,
\partial_{v_{ik}}
&
+
\sum_{1\leqslant r<i}\,
u_{rk}\,
\partial_{w_{ri}}
+
\sum_{1\leqslant s<i}\,
u_{sk}\,
\partial_{w_{is}}
+
2\,u_{ik}\,\partial_{w_{ii}}
\\
&
+
\sum_{i<r\leqslant\ell}\,
u_{rk}\,
\partial_{w_{ri}}
-
\sum_{i<s\leqslant\ell}\,
u_{sk}\,
\partial_{w_{is}}
\ \ \ \ \ \ \ \ \ \ \ \ \ \ \ \ \ \ \ \
{\scriptstyle{(1\,\leqslant\,i\,\leqslant\,\ell,\,\,
1\leqslant\,k\,\leqslant\,m)}},
\endaligned
\]
\[
\aligned
IL_{v_{ik}}
\,=\,
\isqrt\,
\bigg(
\partial_{v_{ik}}
&
-
\sum_{1\leqslant r<i}\,
u_{rk}\,
\partial_{w_{ri}}
-
\sum_{1\leqslant s<i}\,
u_{sk}\,
\partial_{w_{is}}
-
2\,u_{ik}\,\partial_{w_{ii}}
\\
&
-
\sum_{i<r\leqslant\ell}\,
u_{rk}\,
\partial_{w_{ri}}
+
\sum_{i<s\leqslant\ell}\,
u_{sk}\,
\partial_{w_{is}}
\bigg)
\ \ \ \ \ \ \ \ \ \ \ \ \ \ \ \ \ \ \ \
{\scriptstyle{(1\,\leqslant\,i\,\leqslant\,\ell,\,\,
1\leqslant\,k\,\leqslant\,m)}}.
\endaligned
\]

\section{$\mathfrak{su}(p,q)$ CR Models: 
Generators of $\mathfrak{g}_0$}
\label{generators-g0}

With the decomposition:
\[
\aligned
\big(n+m+m\big)^2
&
\,=\,
m^2
+
\tfrac{m(m-1)}{2}
+
\tfrac{m(m-1)}{2}
+
n\,m
+
n\,m
+
\tfrac{n(n-1)}{2}
\\
&
\ \ \ \ \
+
m^2
+
\tfrac{m(m-1)}{2}
+
\tfrac{m(m-1)}{2}
+
n\,m
+
n\,m
+
\tfrac{n(n-1)}{2}
\\
&
\ \ \ \ \
+
n
+
m
+
m,
\endaligned
\]
there are $6 + 6 + 3$ families of generators for $\mathfrak{g}_0$
that are independent of any $\partial_{w_{\ell'\ell''}}$.

Firstly:
\[
\aligned
L_{k_1k_2}^1
&
\,:=\,
\sum_{\ell'=1}^\ell\,
u_{\ell'k_1}\,
\partial_{u_{\ell'k_2}}
-
\sum_{\ell'=1}^\ell\,
v_{\ell'k_2}\,
\partial_{v_{\ell'k_1}}
&
\ \ \ \ \ \ \ \ \ \ \ \ \ \ \ \ \ \ \ \
&
{\scriptstyle{(1\,\leqslant\,k_1,\,k_2\,\leqslant\,m)}},
\\
L_{k_1k_2}^2
&
\,:=\,
\sum_{\ell'=1}^\ell\,
v_{\ell'k_1}\,
\partial_{u_{\ell'k_2}}
-
\sum_{\ell'=1}^\ell\,
v_{\ell'k_2}\,
\partial_{u_{\ell'k_1}}
&
\ \ \ \ \ \ \ \ \ \ \ \ \ \ \ \ \ \ \ \
&
{\scriptstyle{(1\,\leqslant\,k_1\,<\,k_2\,\leqslant\,m)}},
\\
L_{k_1k_2}^3
&
\,:=\,
\sum_{\ell'=1}^\ell\,
u_{\ell'k_1}\,
\partial_{v_{\ell'k_2}}
-
\sum_{\ell'=1}^\ell\,
u_{\ell'k_2}\,
\partial_{v_{\ell'k_1}}
&
\ \ \ \ \ \ \ \ \ \ \ \ \ \ \ \ \ \ \ \
&
{\scriptstyle{(1\,\leqslant\,k_1<k_2\,\leqslant\,m)}},
\\
L_{i,k}^4
&
\,:=\,
\sum_{\ell'=1}^\ell\,
u_{\ell'k}\,
\partial_{z_{\ell'i}}
-
\sum_{\ell'=1}^\ell\,
z_{\ell'i}\,
\partial_{v_{\ell'k}}
&
\ \ \ \ \ \ \ \ \ \ \ \ \ \ \ \ \ \ \ \
&
{\scriptstyle{(
1\,\leqslant\,i\,\leqslant\,n,\,\,
1\,\leqslant\,k\,\leqslant\,m)}},
\\
L_{i,k}^5
&
\,:=\,
\sum_{\ell'=1}^\ell\,
v_{\ell'k}\,
\partial_{z_{\ell'i}}
-
\sum_{\ell'=1}^\ell\,
z_{\ell'i}\,
\partial_{u_{\ell'k}}
&
\ \ \ \ \ \ \ \ \ \ \ \ \ \ \ \ \ \ \ \
&
{\scriptstyle{(
1\,\leqslant\,i\,\leqslant\,n,\,\,
1\,\leqslant\,k\,\leqslant\,m)}},
\\
L_{i_1,i_2}^6
&
\,:=\,
\sum_{\ell'=1}^\ell\,
z_{\ell'i_1}\,
\partial_{z_{\ell'i_2}}
-
\sum_{\ell'=1}^\ell\,
z_{\ell'i_2}\,
\partial_{z_{\ell'i_1}}
&
\ \ \ \ \ \ \ \ \ \ \ \ \ \ \ \ \ \ \ \
&
{\scriptstyle{(
1\,\leqslant\,i_1\,<\,i_2\,\leqslant\,n)}},
\endaligned
\]

Secondly:
\[
\aligned
IL_{k_1k_2}^1
&
\,:=\,
\isqrt\,
\sum_{\ell'=1}^\ell\,
u_{\ell'k_1}\,
\partial_{u_{\ell'k_2}}
+
\isqrt\,
\sum_{\ell'=1}^\ell\,
v_{\ell'k_2}\,
\partial_{v_{\ell'k_1}}
&
\ \ \ \ \ \ \ \ \ \ \ \ \ \ \ \ \ \ \ \
&
{\scriptstyle{(1\,\leqslant\,k_1,\,k_2\,\leqslant\,m)}},
\\
IL_{k_1k_2}^2
&
\,:=\,
\isqrt\,
\sum_{\ell'=1}^\ell\,
v_{\ell'k_1}\,
\partial_{u_{\ell'k_2}}
+
\isqrt\,
\sum_{\ell'=1}^\ell\,
v_{\ell'k_2}\,
\partial_{u_{\ell'k_1}}
&
\ \ \ \ \ \ \ \ \ \ \ \ \ \ \ \ \ \ \ \
&
{\scriptstyle{(1\,\leqslant\,k_1\,<\,k_2\,\leqslant\,m)}},
\\
IL_{k_1k_2}^3
&
\,:=\,
\isqrt\,
\sum_{\ell'=1}^\ell\,
u_{\ell'k_1}\,
\partial_{v_{\ell'k_2}}
+
\isqrt\,
\sum_{\ell'=1}^\ell\,
u_{\ell'k_2}\,
\partial_{v_{\ell'k_1}}
&
\ \ \ \ \ \ \ \ \ \ \ \ \ \ \ \ \ \ \ \
&
{\scriptstyle{(1\,\leqslant\,k_1<k_2\,\leqslant\,m)}},
\\
IL_{i,k}^4
&
\,:=\,
\isqrt\,
\sum_{\ell'=1}^\ell\,
u_{\ell'k}\,
\partial_{z_{\ell'i}}
+
\isqrt\,
\sum_{\ell'=1}^\ell\,
z_{\ell'i}\,
\partial_{v_{\ell'k}}
&
\ \ \ \ \ \ \ \ \ \ \ \ \ \ \ \ \ \ \ \
&
{\scriptstyle{(
1\,\leqslant\,i\,\leqslant\,n,\,\,
1\,\leqslant\,k\,\leqslant\,m)}},
\\
IL_{i,k}^5
&
\,:=\,
\isqrt\,
\sum_{\ell'=1}^\ell\,
v_{\ell'k}\,
\partial_{z_{\ell'i}}
+
\isqrt\,
\sum_{\ell'=1}^\ell\,
z_{\ell'i}\,
\partial_{u_{\ell'k}}
&
\ \ \ \ \ \ \ \ \ \ \ \ \ \ \ \ \ \ \ \
&
{\scriptstyle{(
1\,\leqslant\,i\,\leqslant\,n,\,\,
1\,\leqslant\,k\,\leqslant\,m)}},
\\
IL_{i_1,i_2}^6
&
\,:=\,
\isqrt\,
\sum_{\ell'=1}^\ell\,
z_{\ell'i_1}\,
\partial_{z_{\ell'i_2}}
+
\isqrt\,
\sum_{\ell'=1}^\ell\,
z_{\ell'i_2}\,
\partial_{z_{\ell'i_1}}
&
\ \ \ \ \ \ \ \ \ \ \ \ \ \ \ \ \ \ \ \
&
{\scriptstyle{(
1\,\leqslant\,i_1\,<\,i_2\,\leqslant\,n)}},
\endaligned
\]

Thirdly:
\[
\aligned
IL_i^7
&
\,:=\,
\isqrt\,
\sum_{\ell'=1}^\ell\,
z_{\ell'i}\,
\partial_{z_{\ell'i}}
\ \ \ \ \ \ \ \ \ \ \ \ \ \ \ \ \ \ \ \
&
{\scriptstyle{(1\,\leqslant\,i\,\leqslant\,n)}},
\\
IL_k^8
&
\,:=\,
\isqrt\,
\sum_{\ell'=1}^\ell\,
v_{\ell'k}\,
\partial_{u_{\ell'k}}
\ \ \ \ \ \ \ \ \ \ \ \ \ \ \ \ \ \ \ \
&
{\scriptstyle{(1\,\leqslant\,k\,\leqslant\,m)}},
\\
IL_k^9
&
\,:=\,
\isqrt\,
\sum_{\ell'=1}^\ell\,
u_{\ell'k}\,
\partial_{v_{\ell'k}}
\ \ \ \ \ \ \ \ \ \ \ \ \ \ \ \ \ \ \ \
&
{\scriptstyle{(1\,\leqslant\,k\,\leqslant\,m)}}.
\endaligned
\]

In addition, there are 5 families of generators which do depend
on some $\partial_{w_{\ell'\ell'}}$:
\[
\aligned
{}
&
L_{i,j}
\ \ \ \ \ \ \ \ \ \ \ \ \ \ \ \ \ \ \ \
&
{\scriptstyle{(1\,\leqslant\,i\,<\,j\,\leqslant\,\ell)}},
\\
{}
&
L_{i,j}
\ \ \ \ \ \ \ \ \ \ \ \ \ \ \ \ \ \ \ \
&
{\scriptstyle{(1\,\leqslant\,j\,\leqslant\,i\,\leqslant\,\ell)}},
\\
{}
&
IL_{i,j}
\ \ \ \ \ \ \ \ \ \ \ \ \ \ \ \ \ \ \ \
&
{\scriptstyle{(1\,\leqslant\,i\,<\,j\,\leqslant\,\ell)}},
\\
{}
&
IL_{i,i}
\ \ \ \ \ \ \ \ \ \ \ \ \ \ \ \ \ \ \ \
&
{\scriptstyle{(1\,\leqslant\,i\,\leqslant\,\ell)}},
\\
{}
&
IL_{i,j}
\ \ \ \ \ \ \ \ \ \ \ \ \ \ \ \ \ \ \ \
&
{\scriptstyle{(1\,\leqslant\,j\,<\,i\,\leqslant\,\ell)}}.
\endaligned
\]

For $1 \leqslant i < j \leqslant \ell$:
\[
\aligned
L_{i,j}
&
\,:=\,
\sum_{n'=1}^n\,
z_{in'}\,
\partial_{z_{jn'}}
+
\sum_{m'=1}^m\,
u_{im'}\,
\partial_{u_{jm'}}
+
\sum_{m'=1}^m\,
v_{im'}\,
\partial_{v_{jm'}}
\\
&
\ \ \ \ \
+
\sum_{1\leqslant s<i}\,
w_{is}\,
\partial_{w_{js}}
+
\sum_{i\leqslant s\leqslant j}\,
w_{si}\,
\partial_{w_{js}}
+
\sum_{j<s\leqslant\ell}\,
w_{is}\,\partial_{w_{js}}
\\
&
\ \ \ \ \
+
\sum_{1\leqslant r<i}\,
w_{ri}\,
\partial_{w_{rj}}
-
\sum_{i<r<j}\,
w_{ir}\,
\partial_{w_{rj}}
+
\sum_{j\leqslant r\leqslant\ell}\,
w_{ri}\,
\partial_{w_{rj}}.
\endaligned
\]

For $1 \leqslant j \leqslant i \leqslant \ell$:
\[
\aligned
L_{i,j}
&
\,:=\,
\sum_{n'=1}^n\,
z_{in'}\,
\partial_{z_{jn'}}
+
\sum_{m'=1}^m\,
u_{im'}\,
\partial_{u_{jm'}}
+
\sum_{m'=1}^m\,
v_{im'}\,
\partial_{v_{jm'}}
\\
&
\ \ \ \ \
+
\sum_{1\leqslant s\leqslant j}\,
w_{is}\,
\partial_{w_{js}}
-
\sum_{j<s<i}\,
w_{si}\,
\partial_{w_{js}}
+
\sum_{i<s\leqslant\ell}\,
w_{is}\,\partial_{w_{js}}
\\
&
\ \ \ \ \
+
\sum_{1\leqslant r<j}\,
w_{ri}\,
\partial_{w_{rj}}
+
\sum_{j\leqslant r\leqslant i}\,
w_{ir}\,
\partial_{w_{rj}}
+
\sum_{i<r\leqslant\ell}\,
w_{ri}\,
\partial_{w_{rj}}.
\endaligned
\]

For $1 \leqslant i < j \leqslant \ell$:
\[
\aligned
IL_{i,j}
&
\,:=\,
\isqrt\,
\sum_{n'=1}^n\,
z_{in'}\,
\partial_{z_{jn'}}
+
\isqrt\,
\sum_{m'=1}^m\,
u_{im'}\,
\partial_{u_{jm'}}
+
\isqrt\,
\sum_{m'=1}^m\,
v_{im'}\,
\partial_{v_{jm'}}
\\
&
\ \ \ \ \
-
\isqrt\,
\sum_{1\leqslant s<i}\,
w_{si}\,
\partial_{w_{js}}
+
\isqrt\,
\sum_{i< s\leqslant j}\,
w_{is}\,
\partial_{w_{js}}
+
\isqrt\,
\sum_{j<s\leqslant\ell}\,
w_{si}\,\partial_{w_{js}}
\\
&
\ \ \ \ \
-
\isqrt\,
\sum_{1\leqslant r\leqslant i}\,
w_{ir}\,
\partial_{w_{rj}}
-
\isqrt\,
\sum_{i<r<j}\,
w_{ri}\,
\partial_{w_{rj}}
+
\isqrt\,
\sum_{j\leqslant r\leqslant\ell}\,
w_{ir}\,
\partial_{w_{rj}}.
\endaligned
\]

For $1 \leqslant i \leqslant \ell$:
\[
\aligned
IL_{i,i}
&
\,:=\,
\isqrt\,
\sum_{n'=1}^n\,
z_{in'}\,
\partial_{z_{in'}}
+
\isqrt\,
\sum_{m'=1}^m\,
u_{im'}\,
\partial_{u_{im'}}
+
\isqrt\,
\sum_{m'=1}^m\,
v_{im'}\,
\partial_{v_{im'}}
\\
&
\ \ \ \ \
-
\isqrt\,
\sum_{1\leqslant s<i}\,
w_{si}\,
\partial_{w_{is}}
+
\isqrt\,
\sum_{i<s\leqslant\ell}\,
w_{si}\,
\partial_{w_{is}}
\\
&
\ \ \ \ \
-
\isqrt\,
\sum_{1\leqslant s<i}\,
w_{is}\,
\partial_{w_{si}}
+
\isqrt\,
\sum_{i<s\leqslant\ell}\,
w_{si}\,
\partial_{w_{is}}.
\endaligned
\]

For $1 \leqslant j < i \leqslant \ell$:
\[
\aligned
IL_{i,j}
&
\,:=\,
\isqrt\,
\sum_{n'=1}^n\,
z_{in'}\,
\partial_{z_{jn'}}
+
\isqrt\,
\sum_{m'=1}^m\,
u_{im'}\,
\partial_{u_{jm'}}
+
\isqrt\,
\sum_{m'=1}^m\,
v_{im'}\,
\partial_{v_{jm'}}
\\
&
\ \ \ \ \
-
\isqrt\,
\sum_{1\leqslant s\leqslant j}\,
w_{si}\,
\partial_{w_{js}}
+
\isqrt\,
\sum_{j<s\leqslant i}\,
w_{is}\,
\partial_{w_{js}}
+
\isqrt\,
\sum_{i<s\leqslant\ell}\,
w_{si}\,
\partial_{w_{js}}
\\
&
\ \ \ \ \
-
\isqrt\,
\sum_{1\leqslant r<j}\,
w_{ir}\,
\partial_{w_{rj}}
-
\isqrt\,
\sum_{j\leqslant r<i}\,
w_{ri}\,
\partial_{w_{rj}}
+
\isqrt\,
\sum_{i<r\leqslant\ell}\,
w_{ir}\,
\partial_{w_{rj}}.
\endaligned
\]

\section{$\mathfrak{su}(p,q)$ CR Models: 
Generators of $\mathfrak{g}_1$}
\label{generators-g1}

There are 6 families of generators for $\mathfrak{g}_1$:
\[
\aligned
{}
&
L_{z_{ij}z_{ij}}
&
\ \ \ \ \ \ \ \ \ \ \ \ \ \ \ \ \ \ \ \
&
{\scriptstyle{(1\,\leqslant\,i\,\leqslant\,\ell,\,\,
1\,\leqslant\,j\,\leqslant\,n)}},
\\
{}
&
IL_{z_{ij}z_{ij}}
&
\ \ \ \ \ \ \ \ \ \ \ \ \ \ \ \ \ \ \ \
&
{\scriptstyle{(1\,\leqslant\,i\,\leqslant\,\ell,\,\,
1\,\leqslant\,j\,\leqslant\,n)}},
\\
{}
&
L_{u_{ik}u_{ik}}
&
\ \ \ \ \ \ \ \ \ \ \ \ \ \ \ \ \ \ \ \
&
{\scriptstyle{(1\,\leqslant\,i\,\leqslant\,\ell,\,\,
1\,\leqslant\,k\,\leqslant\,m)}},
\\
{}
&
IL_{u_{ik}u_{ik}}
&
\ \ \ \ \ \ \ \ \ \ \ \ \ \ \ \ \ \ \ \
&
{\scriptstyle{(1\,\leqslant\,i\,\leqslant\,\ell,\,\,
1\,\leqslant\,k\,\leqslant\,m)}},
\\
{}
&
L_{v_{ik}v_{ik}}
&
\ \ \ \ \ \ \ \ \ \ \ \ \ \ \ \ \ \ \ \
&
{\scriptstyle{(1\,\leqslant\,i\,\leqslant\,\ell,\,\,
1\,\leqslant\,k\,\leqslant\,m)}},
\\
{}
&
IL_{v_{ik}f_{ik}}
&
\ \ \ \ \ \ \ \ \ \ \ \ \ \ \ \ \ \ \ \
&
{\scriptstyle{(1\,\leqslant\,i\,\leqslant\,\ell,\,\,
1\,\leqslant\,k\,\leqslant\,m)}}.
\endaligned
\]

For $1 \leqslant i \leqslant \ell$ and $1 \leqslant j \leqslant n$:
\[
\aligned
L_{z_{ij}z_{ij}}
&
\,:=\,
\sum_{\ell'=1}^\ell\,
\sum_{n'=1\atop n'\neq j}^n\,
\big(
2\,z_{\ell'j}\,z_{in'}
\big)\,
\partial_{z_{\ell'n'}}
\\
&
\ \ \ \ \
+
\sum_{1\leqslant\ell'<i}\,
\big(
2\,z_{\ell'j}\,z_{ij}
-
w_{\ell'i}
-
w_{i\ell'}
\big)\,
\partial_{z_{\ell'j}}
+
\big(
2\,z_{ij}^2
-
w_{ii}
\big)\,
\partial_{z_{ij}}
+
\sum_{i<\ell'\leqslant\ell}\,
\big(
2\,z_{\ell'j}\,z_{ij}
-
w_{\ell'i}
+
w_{i\ell'}
\big)\,
\partial_{z_{\ell'j}}
\\
&
\ \ \ \ \
+
\sum_{\ell'=1}^\ell\,
\sum_{m'=1}^m\,
\big(
2\,z_{\ell'j}\,
u_{im'}
\big)\,
\partial_{u_{\ell'm'}}
+
\sum_{\ell'=1}^\ell\,
\sum_{m'=1}^m\,
\big(
2\,z_{\ell'j}\,u_{im'}
\big)\,
\partial_{v_{\ell'm'}}
\endaligned
\]
\[
\aligned
{}
&
+
\sum_{1\leqslant r<s<i}\,
\Big(
z_{rj}\,
\big[
-w_{si}+w_{is}
\big]
+
z_{sj}\,
\big[
w_{ri}-w_{ir}
\big]
\Big)\,
\partial_{w_{rs}}
+
\sum_{1\leqslant r<i}\,
\Big(
z_{rj}\,
\big[
w_{ii}
\big]
+
z_{ij}\,
\big[
w_{ri}
-
w_{ir}
\big]
\Big)\,
\partial_{w_{ri}}
\\
{}
&
+
\sum_{1\leqslant r\leqslant i}\,
\Big(
z_{rj}\,
\big[
w_{is}
+
w_{si}
\big]
+
z_{sj}\,
\big[
w_{ri}
-
w_{ir}
\big]
\Big)\,
\partial_{w_{rs}}
+
\sum_{i<r\leqslant\ell}\,
\Big(
z_{ij}\,
\big[
w_{is}
+
w_{si}
\big]
+
z_{sj}\,
\big[
-w_{ii}
\big]
\Big)\,
\partial_{w_{is}}
\\
{}
&
+
\sum_{i<r<s\leqslant\ell}\,
\Big(
z_{rj}\,
\big[
w_{is}
+
w_{si}
\big]
+
z_{sj}\,
\big[
-w_{ri}
-
w_{ir}
\big]
\Big)\,
\partial_{w_{rs}}
\endaligned
\]
\[
\aligned
{}
&
+
\sum_{1\leqslant r<i}\,
\Big(
z_{rj}\,
\big[
-2\,w_{ri}
+
2\,w_{ir}
\big]
\Big)\,
\partial_{w_{rr}}
+
\Big(
z_{ij}\,
\big[
2\,w_{ii}
\big]
\Big)\,
\partial_{w_{ii}}
+
\sum_{i<r\leqslant\ell}\,
\Big(
z_{rj}\,
\big[
2\,w_{ir}
+
2\,w_{ri}
\big]
\Big)\,
\partial_{w_{rr}}
\endaligned
\]
\[
\aligned
{}
&
+
\sum_{1\leqslant s<r<i}\,
\Big(
z_{rj}\,
\big[
-w_{si}
+
w_{is}
\big]
+
z_{sj}\,
\big[
-w_{ri}
+
w_{ir}
\big]
\Big)\,
\partial_{w_{rs}}
+
\sum_{1\leqslant s<i}\,
\Big(
z_{ij}\,
\big[
-w_{si}
+
w_{is}
\big]
+
w_{sj}\,
\big[
w_{ii}
\big]
\Big)\,
\partial_{w_{is}}
\\
{}
&
+
\sum_{i<r\leqslant\ell}\,
\Big(
z_{rj}\,
\big[
-w_{si}
+
w_{is}
\big]
+
z_{sj}\,
\big[
w_{ir}
+
w_{ri}
\big]
\Big)\,
\partial_{w_{rs}}
+
\sum_{i<r\leqslant\ell}\,
\Big(
z_{rj}\,
\big
[w_{ii}
\big]
+
z_{ij}\,
\big[
w_{ir}
+
w_{ri}
\big]
\Big)\,
\partial_{w_{ri}}
\\
{}
&
+
\sum_{i<s<r\leqslant\ell}\,
\Big(
z_{rj}\,
\big[
w_{is}
+
w_{si}
\big]
+
z_{sj}\,
\big[
w_{ir}
+
w_{ri}
\big]
\Big)\,
\partial_{w_{rs}}.
\endaligned
\]

For $1 \leqslant i \leqslant \ell$ and $1 \leqslant j \leqslant n$:
\[
\aligned
IL_{z_{ij}z_{ij}}
&
\,:=\,
\isqrt\,
\sum_{\ell'=1}^\ell\,
\sum_{n'=1\atop n'\neq j}^n\,
\big(
2\,z_{\ell'j}\,z_{in'}
\big)\,
\partial_{z_{\ell'n'}}
\\
&
\ \ \ \ \
+
\isqrt\,
\sum_{1\leqslant\ell'<i}\,
\big(
2\,z_{\ell'j}\,z_{ij}
+
w_{\ell'i}
+
w_{i\ell'}
\big)\,
\partial_{z_{\ell'j}}
+
\isqrt\,
\big(
2\,z_{ij}^2
+
w_{ii}
\big)\,
\partial_{z_{ij}}
\\
{}
&
\ \ \ \ \ 
+
\isqrt\,
\sum_{i<\ell'\leqslant\ell}\,
\big(
2\,z_{\ell'j}\,z_{ij}
+
w_{\ell'i}
-
w_{i\ell'}
\big)\,
\partial_{z_{\ell'j}}
\\
&
\ \ \ \ \
+
\isqrt\,
\sum_{\ell'=1}^\ell\,
\sum_{m'=1}^m\,
\big(
2\,z_{\ell'j}\,
u_{im'}
\big)\,
\partial_{u_{\ell'm'}}
+
\isqrt\,
\sum_{\ell'=1}^\ell\,
\sum_{m'=1}^m\,
\big(
2\,z_{\ell'j}\,u_{im'}
\big)\,
\partial_{v_{\ell'm'}}
\endaligned
\]
\[
\aligned
{}
&
+
\isqrt\,
\sum_{1\leqslant r<s<i}\,
\Big(
z_{rj}\,
\big[
-w_{si}+w_{is}
\big]
+
z_{sj}\,
\big[
w_{ri}-w_{ir}
\big]
\Big)\,
\partial_{w_{rs}}
+
\isqrt\,
\sum_{1\leqslant r<i}\,
\Big(
z_{rj}\,
\big[
w_{ii}
\big]
+
z_{ij}\,
\big[
w_{ri}
-
w_{ir}
\big]
\Big)\,
\partial_{w_{ri}}
\\
{}
&
+
\isqrt\,
\sum_{1\leqslant r\leqslant i}\,
\Big(
z_{rj}\,
\big[
w_{is}
+
w_{si}
\big]
+
z_{sj}\,
\big[
w_{ri}
-
w_{ir}
\big]
\Big)\,
\partial_{w_{rs}}
+
\isqrt\,
\sum_{i<r\leqslant\ell}\,
\Big(
z_{ij}\,
\big[
w_{is}
+
w_{si}
\big]
+
z_{sj}\,
\big[
-w_{ii}
\big]
\Big)\,
\partial_{w_{is}}
\\
{}
&
+
\isqrt\,
\sum_{i<r<s\leqslant\ell}\,
\Big(
z_{rj}\,
\big[
w_{is}
+
w_{si}
\big]
+
z_{sj}\,
\big[
-w_{ri}
-
w_{ir}
\big]
\Big)\,
\partial_{w_{rs}}
\endaligned
\]
\[
\aligned
{}
&
+
\isqrt\,
\sum_{1\leqslant r<i}\,
\Big(
z_{rj}\,
\big[
-2\,w_{ri}
+
2\,w_{ir}
\big]
\Big)\,
\partial_{w_{rr}}
+
\Big(
z_{ij}\,
\big[
2\,w_{ii}
\big]
\Big)\,
\partial_{w_{ii}}
+
\isqrt\,
\sum_{i<r\leqslant\ell}\,
\Big(
z_{rj}\,
\big[
2\,w_{ir}
+
2\,w_{ri}
\big]
\Big)\,
\partial_{w_{rr}}
\endaligned
\]
\[
\aligned
{}
&
+
\isqrt\,
\sum_{1\leqslant s<r<i}\,
\Big(
z_{rj}\,
\big[
-w_{si}
+
w_{is}
\big]
+
z_{sj}\,
\big[
-w_{ri}
+
w_{ir}
\big]
\Big)\,
\partial_{w_{rs}}
\\
{}
&
+
\isqrt\,
\sum_{1\leqslant s<i}\,
\Big(
z_{ij}\,
\big[
-w_{si}
+
w_{is}
\big]
+
w_{sj}\,
\big[
w_{ii}
\big]
\Big)\,
\partial_{w_{is}}
\\
{}
&
+
\isqrt\,
\sum_{i<r\leqslant\ell}\,
\Big(
z_{rj}\,
\big[
-w_{si}
+
w_{is}
\big]
+
z_{sj}\,
\big[
w_{ir}
+
w_{ri}
\big]
\Big)\,
\partial_{w_{rs}}
+
\isqrt\,
\sum_{i<r\leqslant\ell}\,
\Big(
z_{rj}\,
\big
[w_{ii}
\big]
+
z_{ij}\,
\big[
w_{ir}
+
w_{ri}
\big]
\Big)\,
\partial_{w_{ri}}
\\
{}
&
+
\isqrt\,
\sum_{i<s<r\leqslant\ell}\,
\Big(
z_{rj}\,
\big[
w_{is}
+
w_{si}
\big]
+
z_{sj}\,
\big[
w_{ir}
+
w_{ri}
\big]
\Big)\,
\partial_{w_{rs}}.
\endaligned
\]

For $1 \leqslant i \leqslant \ell$ and $1 \leqslant k \leqslant m$:
\[
\aligned
L_{u_{ik}u_{ik}}
&
\,:=\,
\sum_{\ell'=1}^\ell\,
\sum_{n'=1}^n\,
\big(
2\,z_{in'}\,u_{\ell'k}
\big)\,
\partial_{z_{\ell'n'}}
\\
&
\ \ \ \ \
+
\sum_{\ell'=1}^\ell\,
\sum_{m'=1}^m\,
\big(
2\,u_{im'}\,u_{\ell'm'}
\big)\,
\partial_{u_{\ell'm'}}
\\
&
\ \ \ \ \
+
\sum_{\ell'=1}^\ell\,
\sum_{m'=1\atop m'\neq k}^m\,
\big(
2\,u_{\ell'k}\,v_{im'}
\big)\,
\partial_{v_{\ell'm'}}
\\
&
\!\!\!\!\!\!\!\!\!\!\!\!\!\!\!\!\!\!\!\!\!\!\!\!\!
+
\sum_{1\leqslant\ell'<i}\,
\big(
2\,u_{\ell'k}\,v_{ik}
-
w_{\ell'i}
-
w_{i\ell'}
\big)\,
\partial_{v_{\ell'k}}
+
\big(
2\,u_{ik}\,v_{ik}
-
w_{ii}
\big)\,
\partial_{v_{ik}}
+
\sum_{i<\ell'\leqslant\ell}\,
\big(
2\,u_{\ell'k}\,v_{ik}
-
w_{\ell'i}
+
w_{i\ell'}
\big)\,
\partial_{v_{\ell'k}}
\endaligned
\]
\[
\aligned
{}
&
+
\sum_{1\leqslant r<s<i}\,
\Big(
u_{rk}\,
\big[
-w_{si}+w_{is}
\big]
+
u_{sk}\,
\big[
w_{ri}-w_{ir}
\big]
\Big)\,
\partial_{w_{rs}}
+
\sum_{1\leqslant r<i}\,
\Big(
u_{rk}\,
\big[
w_{ii}
\big]
+
u_{ik}\,
\big[
w_{ri}
-
w_{ir}
\big]
\Big)\,
\partial_{w_{ri}}
\\
{}
&
+
\sum_{1\leqslant r\leqslant i}\,
\Big(
u_{rk}\,
\big[
w_{is}
+
w_{si}
\big]
+
u_{sk}\,
\big[
w_{ri}
-
w_{ir}
\big]
\Big)\,
\partial_{w_{rs}}
+
\sum_{i<r\leqslant\ell}\,
\Big(
u_{ik}\,
\big[
w_{is}
+
w_{si}
\big]
+
u_{sk}\,
\big[
-w_{ii}
\big]
\Big)\,
\partial_{w_{is}}
\\
{}
&
+
\sum_{i<r<s\leqslant\ell}\,
\Big(
u_{rk}\,
\big[
w_{is}
+
w_{si}
\big]
+
u_{sk}\,
\big[
-w_{ri}
-
w_{ir}
\big]
\Big)\,
\partial_{w_{rs}}
\endaligned
\]
\[
\aligned
{}
&
+
\sum_{1\leqslant r<i}\,
\Big(
u_{rk}\,
\big[
-2\,w_{ri}
+
2\,w_{ir}
\big]
\Big)\,
\partial_{w_{rr}}
+
\Big(
u_{ik}\,
\big[
2\,w_{ii}
\big]
\Big)\,
\partial_{w_{ii}}
+
\sum_{i<r\leqslant\ell}\,
\Big(
u_{rk}\,
\big[
2\,w_{ir}
+
2\,w_{ri}
\big]
\Big)\,
\partial_{w_{rr}}
\endaligned
\]
\[
\aligned
{}
&
+
\sum_{1\leqslant s<r<i}\,
\Big(
u_{rk}\,
\big[
-w_{si}
+
w_{is}
\big]
+
u_{sk}\,
\big[
-w_{ri}
+
w_{ir}
\big]
\Big)\,
\partial_{w_{rs}}
+
\sum_{1\leqslant s<i}\,
\Big(
u_{ik}\,
\big[
-w_{si}
+
w_{is}
\big]
+
w_{sk}\,
\big[
w_{ii}
\big]
\Big)\,
\partial_{w_{is}}
\\
{}
&
+
\sum_{i<r\leqslant\ell}\,
\Big(
u_{rk}\,
\big[
-w_{si}
+
w_{is}
\big]
+
u_{sk}\,
\big[
w_{ir}
+
w_{ri}
\big]
\Big)\,
\partial_{w_{rs}}
+
\sum_{i<r\leqslant\ell}\,
\Big(
u_{rk}\,
\big
[w_{ii}
\big]
+
u_{ik}\,
\big[
w_{ir}
+
w_{ri}
\big]
\Big)\,
\partial_{w_{ri}}
\\
{}
&
+
\sum_{i<s<r\leqslant\ell}\,
\Big(
u_{rk}\,
\big[
w_{is}
+
w_{si}
\big]
+
u_{sk}\,
\big[
w_{ir}
+
w_{ri}
\big]
\Big)\,
\partial_{w_{rs}}.
\endaligned
\]

For $1 \leqslant i \leqslant \ell$ and $1 \leqslant k \leqslant m$:
\[
\aligned
IL_{u_{ik}u_{ik}}
&
\,:=\,
\isqrt\,
\sum_{\ell'=1}^\ell\,
\sum_{n'=1}^n\,
\big(
2\,z_{in'}\,u_{\ell'k}
\big)\,
\partial_{z_{\ell'n'}}
\\
&
\ \ \ \ \
+
\isqrt\,
\sum_{\ell'=1}^\ell\,
\sum_{m'=1}^m\,
\big(
2\,u_{im'}\,u_{\ell'm'}
\big)\,
\partial_{u_{\ell'm'}}
\\
&
\ \ \ \ \
+
\isqrt\,
\sum_{\ell'=1}^\ell\,
\sum_{m'=1\atop m'\neq k}^m\,
\big(
2\,u_{\ell'k}\,v_{im'}
\big)\,
\partial_{v_{\ell'm'}}
\\
&
\!\!\!\!\!\!\!\!\!\!\!\!\!\!\!\!\!\!\!\!\!\!\!\!\!
\!\!\!\!\!\!\!\!\!\!\!\!\!\!\!\!\!\!\!\!\!\!\!\!\!
+
\isqrt\,
\sum_{1\leqslant\ell'<i}\,
\big(
2\,u_{\ell'k}\,v_{ik}
+
w_{\ell'i}
+
w_{i\ell'}
\big)\,
\partial_{v_{\ell'k}}
+
\isqrt\,
\big(
2\,u_{ik}\,v_{ik}
+
w_{ii}
\big)\,
\partial_{v_{ik}}
+
\isqrt\,
\sum_{i<\ell'\leqslant\ell}\,
\big(
2\,u_{\ell'k}\,v_{ik}
+
w_{\ell'i}
-
w_{i\ell'}
\big)\,
\partial_{v_{\ell'k}}
\endaligned
\]
\[
\aligned
{}
&
+
\isqrt\,
\sum_{1\leqslant r<s<i}\,
\Big(
u_{rk}\,
\big[
-w_{si}+w_{is}
\big]
+
u_{sk}\,
\big[
w_{ri}-w_{ir}
\big]
\Big)\,
\partial_{w_{rs}}
+
\isqrt\,
\sum_{1\leqslant r<i}\,
\Big(
u_{rk}\,
\big[
w_{ii}
\big]
+
u_{ik}\,
\big[
w_{ri}
-
w_{ir}
\big]
\Big)\,
\partial_{w_{ri}}
\\
{}
&
+
\isqrt\,
\sum_{1\leqslant r\leqslant i}\,
\Big(
u_{rk}\,
\big[
w_{is}
+
w_{si}
\big]
+
u_{sk}\,
\big[
w_{ri}
-
w_{ir}
\big]
\Big)\,
\partial_{w_{rs}}
+
\isqrt\,
\sum_{i<r\leqslant\ell}\,
\Big(
u_{ik}\,
\big[
w_{is}
+
w_{si}
\big]
+
u_{sk}\,
\big[
-w_{ii}
\big]
\Big)\,
\partial_{w_{is}}
\\
{}
&
+
\isqrt\,
\sum_{i<r<s\leqslant\ell}\,
\Big(
u_{rk}\,
\big[
w_{is}
+
w_{si}
\big]
+
u_{sk}\,
\big[
-w_{ri}
-
w_{ir}
\big]
\Big)\,
\partial_{w_{rs}}
\endaligned
\]
\[
\aligned
{}
&
+
\isqrt\,
\sum_{1\leqslant r<i}\,
\Big(
u_{rk}\,
\big[
-2\,w_{ri}
+
2\,w_{ir}
\big]
\Big)\,
\partial_{w_{rr}}
+
\Big(
u_{ik}\,
\big[
2\,w_{ii}
\big]
\Big)\,
\partial_{w_{ii}}
+
\isqrt\,
\sum_{i<r\leqslant\ell}\,
\Big(
u_{rk}\,
\big[
2\,w_{ir}
+
2\,w_{ri}
\big]
\Big)\,
\partial_{w_{rr}}
\endaligned
\]
\[
\aligned
{}
&
+
\isqrt\,
\sum_{1\leqslant s<r<i}\,
\Big(
u_{rk}\,
\big[
-w_{si}
+
w_{is}
\big]
+
u_{sk}\,
\big[
-w_{ri}
+
w_{ir}
\big]
\Big)\,
\partial_{w_{rs}}
\\
{}
&
+
\isqrt\,
\sum_{1\leqslant s<i}\,
\Big(
u_{ik}\,
\big[
-w_{si}
+
w_{is}
\big]
+
w_{sk}\,
\big[
w_{ii}
\big]
\Big)\,
\partial_{w_{is}}
\\
{}
&
+
\isqrt\,
\sum_{i<r\leqslant\ell}\,
\Big(
u_{rk}\,
\big[
-w_{si}
+
w_{is}
\big]
+
u_{sk}\,
\big[
w_{ir}
+
w_{ri}
\big]
\Big)\,
\partial_{w_{rs}}
+
\isqrt\,
\sum_{i<r\leqslant\ell}\,
\Big(
u_{rk}\,
\big
[w_{ii}
\big]
+
u_{ik}\,
\big[
w_{ir}
+
w_{ri}
\big]
\Big)\,
\partial_{w_{ri}}
\\
{}
&
+
\isqrt\,
\sum_{i<s<r\leqslant\ell}\,
\Big(
u_{rk}\,
\big[
w_{is}
+
w_{si}
\big]
+
u_{sk}\,
\big[
w_{ir}
+
w_{ri}
\big]
\Big)\,
\partial_{w_{rs}}.
\endaligned
\]

For $1 \leqslant i \leqslant \ell$ and $1 \leqslant k \leqslant m$:
\[
\aligned
L_{v_{ik}v_{ik}}
&
\,:=\,
\sum_{\ell'=1}^\ell\,
\sum_{n'=1}^n\,
\big(
2\,z_{in'}\,v_{\ell'k}
\big)\,
\partial_{z_{\ell'n'}}
\\
&
\ \ \ \ \
+
\sum_{\ell'=1}^\ell\,
\sum_{m'=1}^m\,
\big(
2\,v_{im'}\,v_{\ell'm'}
\big)\,
\partial_{v_{\ell'm'}}
\\
&
\ \ \ \ \
+
\sum_{\ell'=1}^\ell\,
\sum_{m'=1\atop m'\neq k}^m\,
\big(
2\,v_{\ell'k}\,v_{im'}
\big)\,
\partial_{v_{\ell'm'}}
\\
&
\!\!\!\!\!\!\!\!\!\!\!\!\!\!\!\!\!\!\!\!\!\!\!\!\!
+
\sum_{1\leqslant\ell'<i}\,
\big(
2\,v_{\ell'k}\,v_{ik}
-
w_{\ell'i}
-
w_{i\ell'}
\big)\,
\partial_{v_{\ell'k}}
+
\big(
2\,v_{ik}\,v_{ik}
-
w_{ii}
\big)\,
\partial_{v_{ik}}
+
\sum_{i<\ell'\leqslant\ell}\,
\big(
2\,v_{\ell'k}\,v_{ik}
-
w_{\ell'i}
+
w_{i\ell'}
\big)\,
\partial_{v_{\ell'k}}
\endaligned
\]
\[
\aligned
{}
&
+
\sum_{1\leqslant r<s<i}\,
\Big(
v_{rk}\,
\big[
-w_{si}+w_{is}
\big]
+
v_{sk}\,
\big[
w_{ri}-w_{ir}
\big]
\Big)\,
\partial_{w_{rs}}
+
\sum_{1\leqslant r<i}\,
\Big(
v_{rk}\,
\big[
w_{ii}
\big]
+
v_{ik}\,
\big[
w_{ri}
-
w_{ir}
\big]
\Big)\,
\partial_{w_{ri}}
\\
{}
&
+
\sum_{1\leqslant r\leqslant i}\,
\Big(
v_{rk}\,
\big[
w_{is}
+
w_{si}
\big]
+
v_{sk}\,
\big[
w_{ri}
-
w_{ir}
\big]
\Big)\,
\partial_{w_{rs}}
+
\sum_{i<r\leqslant\ell}\,
\Big(
v_{ik}\,
\big[
w_{is}
+
w_{si}
\big]
+
v_{sk}\,
\big[
-w_{ii}
\big]
\Big)\,
\partial_{w_{is}}
\\
{}
&
+
\sum_{i<r<s\leqslant\ell}\,
\Big(
v_{rk}\,
\big[
w_{is}
+
w_{si}
\big]
+
v_{sk}\,
\big[
-w_{ri}
-
w_{ir}
\big]
\Big)\,
\partial_{w_{rs}}
\endaligned
\]
\[
\aligned
{}
&
+
\sum_{1\leqslant r<i}\,
\Big(
v_{rk}\,
\big[
-2\,w_{ri}
+
2\,w_{ir}
\big]
\Big)\,
\partial_{w_{rr}}
+
\Big(
v_{ik}\,
\big[
2\,w_{ii}
\big]
\Big)\,
\partial_{w_{ii}}
+
\sum_{i<r\leqslant\ell}\,
\Big(
v_{rk}\,
\big[
2\,w_{ir}
+
2\,w_{ri}
\big]
\Big)\,
\partial_{w_{rr}}
\endaligned
\]
\[
\aligned
{}
&
+
\sum_{1\leqslant s<r<i}\,
\Big(
v_{rk}\,
\big[
-w_{si}
+
w_{is}
\big]
+
v_{sk}\,
\big[
-w_{ri}
+
w_{ir}
\big]
\Big)\,
\partial_{w_{rs}}
+
\sum_{1\leqslant s<i}\,
\Big(
v_{ik}\,
\big[
-w_{si}
+
w_{is}
\big]
+
w_{sk}\,
\big[
w_{ii}
\big]
\Big)\,
\partial_{w_{is}}
\\
{}
&
+
\sum_{i<r\leqslant\ell}\,
\Big(
v_{rk}\,
\big[
-w_{si}
+
w_{is}
\big]
+
v_{sk}\,
\big[
w_{ir}
+
w_{ri}
\big]
\Big)\,
\partial_{w_{rs}}
+
\sum_{i<r\leqslant\ell}\,
\Big(
v_{rk}\,
\big
[w_{ii}
\big]
+
v_{ik}\,
\big[
w_{ir}
+
w_{ri}
\big]
\Big)\,
\partial_{w_{ri}}
\\
{}
&
+
\sum_{i<s<r\leqslant\ell}\,
\Big(
v_{rk}\,
\big[
w_{is}
+
w_{si}
\big]
+
v_{sk}\,
\big[
w_{ir}
+
w_{ri}
\big]
\Big)\,
\partial_{w_{rs}}.
\endaligned
\]

For $1 \leqslant i \leqslant \ell$ and $1 \leqslant k \leqslant m$:
\[
\aligned
IL_{v_{ik}v_{ik}}
&
\,:=\,
\isqrt\,
\sum_{\ell'=1}^\ell\,
\sum_{n'=1}^n\,
\big(
2\,z_{in'}\,v_{\ell'k}
\big)\,
\partial_{z_{\ell'n'}}
\\
&
\ \ \ \ \
+
\isqrt\,
\sum_{\ell'=1}^\ell\,
\sum_{m'=1}^m\,
\big(
2\,v_{im'}\,v_{\ell'm'}
\big)\,
\partial_{v_{\ell'm'}}
\\
&
\ \ \ \ \
+
\isqrt\,
\sum_{\ell'=1}^\ell\,
\sum_{m'=1\atop m'\neq k}^m\,
\big(
2\,v_{\ell'k}\,v_{im'}
\big)\,
\partial_{v_{\ell'm'}}
\\
&
\!\!\!\!\!\!\!\!\!\!\!\!\!\!\!\!\!\!\!\!\!\!\!\!\!
\!\!\!\!\!\!\!\!\!\!\!\!\!\!\!\!\!\!\!\!\!\!\!\!\!
+
\isqrt\,
\sum_{1\leqslant\ell'<i}\,
\big(
2\,v_{\ell'k}\,v_{ik}
+
w_{\ell'i}
+
w_{i\ell'}
\big)\,
\partial_{v_{\ell'k}}
+
\isqrt\,
\big(
2\,v_{ik}\,v_{ik}
+
w_{ii}
\big)\,
\partial_{v_{ik}}
+
\isqrt\,
\sum_{i<\ell'\leqslant\ell}\,
\big(
2\,v_{\ell'k}\,v_{ik}
+
w_{\ell'i}
-
w_{i\ell'}
\big)\,
\partial_{v_{\ell'k}}
\endaligned
\]
\[
\aligned
{}
&
+
\isqrt\,
\sum_{1\leqslant r<s<i}\,
\Big(
v_{rk}\,
\big[
-w_{si}+w_{is}
\big]
+
v_{sk}\,
\big[
w_{ri}-w_{ir}
\big]
\Big)\,
\partial_{w_{rs}}
+
\isqrt\,
\sum_{1\leqslant r<i}\,
\Big(
v_{rk}\,
\big[
w_{ii}
\big]
+
v_{ik}\,
\big[
w_{ri}
-
w_{ir}
\big]
\Big)\,
\partial_{w_{ri}}
\\
{}
&
+
\isqrt\,
\sum_{1\leqslant r\leqslant i}\,
\Big(
v_{rk}\,
\big[
w_{is}
+
w_{si}
\big]
+
v_{sk}\,
\big[
w_{ri}
-
w_{ir}
\big]
\Big)\,
\partial_{w_{rs}}
+
\isqrt\,
\sum_{i<r\leqslant\ell}\,
\Big(
v_{ik}\,
\big[
w_{is}
+
w_{si}
\big]
+
v_{sk}\,
\big[
-w_{ii}
\big]
\Big)\,
\partial_{w_{is}}
\\
{}
&
+
\isqrt\,
\sum_{i<r<s\leqslant\ell}\,
\Big(
v_{rk}\,
\big[
w_{is}
+
w_{si}
\big]
+
v_{sk}\,
\big[
-w_{ri}
-
w_{ir}
\big]
\Big)\,
\partial_{w_{rs}}
\endaligned
\]
\[
\aligned
{}
&
+
\isqrt\,
\sum_{1\leqslant r<i}\,
\Big(
v_{rk}\,
\big[
-2\,w_{ri}
+
2\,w_{ir}
\big]
\Big)\,
\partial_{w_{rr}}
+
\Big(
v_{ik}\,
\big[
2\,w_{ii}
\big]
\Big)\,
\partial_{w_{ii}}
+
\isqrt\,
\sum_{i<r\leqslant\ell}\,
\Big(
v_{rk}\,
\big[
2\,w_{ir}
+
2\,w_{ri}
\big]
\Big)\,
\partial_{w_{rr}}
\endaligned
\]
\[
\aligned
{}
&
+
\isqrt\,
\sum_{1\leqslant s<r<i}\,
\Big(
v_{rk}\,
\big[
-w_{si}
+
w_{is}
\big]
+
v_{sk}\,
\big[
-w_{ri}
+
w_{ir}
\big]
\Big)\,
\partial_{w_{rs}}
\\
{}
&
+
\isqrt\,
\sum_{1\leqslant s<i}\,
\Big(
v_{ik}\,
\big[
-w_{si}
+
w_{is}
\big]
+
w_{sk}\,
\big[
w_{ii}
\big]
\Big)\,
\partial_{w_{is}}
\\
{}
&
+
\isqrt\,
\sum_{i<r\leqslant\ell}\,
\Big(
v_{rk}\,
\big[
-w_{si}
+
w_{is}
\big]
+
v_{sk}\,
\big[
w_{ir}
+
w_{ri}
\big]
\Big)\,
\partial_{w_{rs}}
+
\isqrt\,
\sum_{i<r\leqslant\ell}\,
\Big(
v_{rk}\,
\big
[w_{ii}
\big]
+
v_{ik}\,
\big[
w_{ir}
+
w_{ri}
\big]
\Big)\,
\partial_{w_{ri}}
\\
{}
&
+
\isqrt\,
\sum_{i<s<r\leqslant\ell}\,
\Big(
v_{rk}\,
\big[
w_{is}
+
w_{si}
\big]
+
v_{sk}\,
\big[
w_{ir}
+
w_{ri}
\big]
\Big)\,
\partial_{w_{rs}}.
\endaligned
\]

\section{$\mathfrak{su}(p,q)$ CR Models: 
Generators of $\mathfrak{g}_2$}
\label{generators-g2}

There are 3 families of generators for $\mathfrak{g}_2$:
\[
\aligned
{}
&
L_{w_{ij}w_{ij}}
&
\ \ \ \ \ \ \ \ \ \ \ \ \ \ \ \ \ \ \ \
&
{\scriptstyle{(1\,\leqslant\,i\,<\,j\,\leqslant\,\ell)}},
\\
{}
&
L_{w_{ii}w_{ii}}
&
\ \ \ \ \ \ \ \ \ \ \ \ \ \ \ \ \ \ \ \
&
{\scriptstyle{(1\,\leqslant\,i\,\leqslant\,\ell)}},
\\
{}
&
IL_{w_{ij}w_{ij}}
&
\ \ \ \ \ \ \ \ \ \ \ \ \ \ \ \ \ \ \ \
&
{\scriptstyle{(1\,\leqslant\,i\,>\,j\,\leqslant\,\ell)}}.
\endaligned
\]

For $1 \leqslant i < j \leqslant \ell$:
\[
\aligned
L_{w_{ij}w_{ij}}
&
\,:=\,
\sum_{1\leqslant\ell'<i}\,
\sum_{n'=1}^n\,
\Big(
z_{in'}
\big[
w_{j\ell'}
+
w_{\ell' j}
\big]
+
z_{jn'}
\big[
-
w_{i\ell'}
-
w_{\ell'i}
\big]
\Big)\,
\partial_{z_{\ell'n'}}
\\
&
\ \ \ \ \
+
\sum_{n'=1}^n\,
\Big(z_{in'}
\big[
w_{ij}
+
w_{ji}
\big]
+
z_{jn'}
\big[
-
w_{ii}
\big]
\Big)\,
\partial_{z_{in'}}
\\
&
\ \ \ \ \
+
\sum_{i<\ell'<j}\,
\sum_{n'=1}^n\,
\Big(
z_{in'}
\big[
w_{j\ell'}
+
w_{\ell' j}
\big]
+
z_{jn'}
\big[
w_{i\ell'}
-
w_{\ell'i}
\big]
\Big)\,\partial_{z_{\ell'n'}}
\\
&
\ \ \ \ \
+
\sum_{n'=1}^n\,
\Big(z_{in'}
\big[
w_{jj}
\big]
+
z_{jn'}
\big[
w_{ij}
-
w_{ji}
\big]
\Big)\,
\partial_{z_{jn'}}
\\
&
\ \ \ \ \
+
\sum_{j<\ell'\leqslant\ell}\,
\sum_{n'=1}^n\,
\Big(
z_{in'}
\big[
-w_{j\ell'}
+
w_{\ell' j}
\big]
+
z_{jn'}
\big[
w_{i\ell'}
-
w_{\ell'i}
\big]
\Big)\,\partial_{z_{\ell'n'}}
\endaligned
\]
\[
\aligned
{}
&
+
\sum_{1\leqslant\ell'<i}\,
\sum_{m'=1}^m\,
\Big(
u_{im'}
\big[
w_{j\ell'}
+
w_{\ell'j}
\big]
+
u_{jm'}
\big[
-
w_{i\ell'}
-
w_{\ell'i}
\big]
\Big)\,
\partial_{u_{\ell'm'}}
\\
&
\ \ \ \ \ 
+
\sum_{m'=1}^m\,
\Big(
u_{im'}
\big[
w_{ij}
+
w_{ji}
\big]
+
u_{jm'}
\big[
-
w_{ii}
\big]
\Big)\,
\partial_{u_{im'}}
\\
&
\ \ \ \ \
+
\sum_{i<\ell'<j}\,
\sum_{m'=1}^m\,
\Big(
u_{im'}
\big[
w_{j\ell'}
+
w_{\ell'j}
\big]
+
u_{jm'}
\big[
w_{i\ell'}
-
w_{\ell'i}
\big]
\Big)\,
\partial_{u_{\ell'm'}}
\\
&
\ \ \ \ \ 
+
\sum_{m'=1}^m\,
\Big(
u_{im'}
\big[
w_{jj}
\big]
+
u_{jm'}
\big[
w_{ij}
-
w_{ji}
\big]
\Big)\,
\partial_{u_{jm'}}
\\
&
\ \ \ \ \
+
\sum_{j<\ell'\leqslant\ell}\,
\sum_{m'=1}^m\,
\Big(
u_{im'}
\big[
-
w_{j\ell'}
+
w_{\ell'j}
\big]
+
u_{jm'}
\big[
w_{i\ell'}
-
w_{\ell'i}
\big]
\Big)\,
\partial_{u_{\ell'm'}}
\endaligned
\]
\[
\aligned
{}
&
+
\sum_{1\leqslant\ell'<i}\,
\sum_{m'=1}^m\,
\Big(
v_{im'}
\big[
w_{j\ell'}
+
w_{\ell'j}
\big]
+
v_{jm'}
\big[
-
w_{i\ell'}
-
w_{\ell'i}
\big]
\Big)\,
\partial_{v_{\ell'm'}}
\\
&
\ \ \ \ \ 
+
\sum_{m'=1}^m\,
\Big(
v_{im'}
\big[
w_{ij}
+
w_{ji}
\big]
+
v_{jm'}
\big[
-
w_{ii}
\big]
\Big)\,
\partial_{v_{im'}}
\\
&
\ \ \ \ \
+
\sum_{i<\ell'<j}\,
\sum_{m'=1}^m\,
\Big(
v_{im'}
\big[
w_{j\ell'}
+
w_{\ell'j}
\big]
+
v_{jm'}
\big[
w_{i\ell'}
-
w_{\ell'i}
\big]
\Big)\,
\partial_{v_{\ell'm'}}
\\
&
\ \ \ \ \ 
+
\sum_{m'=1}^m\,
\Big(
v_{im'}
\big[
w_{jj}
\big]
+
v_{jm'}
\big[
w_{ij}
-
w_{ji}
\big]
\Big)\,
\partial_{v_{jm'}}
\\
&
\ \ \ \ \
+
\sum_{j<\ell'\leqslant\ell}\,
\sum_{m'=1}^m\,
\Big(
v_{im'}
\big[
-
w_{j\ell'}
+
w_{\ell'j}
\big]
+
v_{jm'}
\big[
w_{i\ell'}
-
w_{\ell'i}
\big]
\Big)\,
\partial_{v_{\ell'm'}}
\endaligned
\]
\[
\aligned
{}
&
+
\sum_{1\leqslant r<i}\,
\big(
-2\,w_{ri}\,w_{jr}
+
2\,w_{rj}\,w_{ir}
\big)\,
\partial_{w_{rr}}
\\
&
\ \ \ \ \
+
\sum_{1\leqslant r<i}\,
\big(
w_{ri}\,w_{ij}
-
w_{ir}\,w_{ji}
+
w_{ii}\,w_{jr}
\big)\,
\partial_{w_{ri}}
\\
&
\ \ \ \ \ 
+
\sum_{1\leqslant r<i}\,
\big(
w_{rj}\,w_{ij}
-
w_{ir}\,w_{jj}
+
w_{ji}\,w_{jr}
\big)\,
\partial_{w_{rj}}
\\
&
\ \ \ \ \
+
\sum_{1\leqslant s<i}\,
\big(
-w_{si}\,w_{ji}
+
w_{sj}\,w_{ii}
+
w_{is}\,w_{ij}
\big)\,
\partial_{w_{is}}
\\
&
\ \ \ \ \
+
\big(
2\,w_{ii}\,w_{ij}
\big)\,
\partial_{w_{ii}}
\\
&
\ \ \ \ \
+
\sum_{i<s<j}\,
\big(
-w_{ii}\,w_{js}
+
w_{is}\,w_{ij}
+
w_{si}\,w_{ji}
\big)\,
\partial_{w_{is}}
\\
&
\ \ \ \ \
+
\big(
-w_{ii}\,w_{jj}
+
w_{ij}^2
+
w_{ji}^2
\big)\,
\partial_{w_{ij}}
\endaligned
\]
\[
\aligned
{}
&
+
\sum_{j<s\leqslant\ell}\,
\big(
-w_{ii}\,w_{sj}
+
w_{ij}\,w_{is}
+
w_{ji}\,w_{si}
\big)\,
\partial_{w_{is}}
\\
&
\ \ \ \ \
+
\sum_{i<r<j}\,
\big(
w_{ii}\,w_{rj}
+
w_{ji}\,w_{ir}
+
w_{ij}\,w_{ri}
\big)\,
\partial_{w_{ri}}
\\
&
\ \ \ \ \
+
\sum_{i<r<j}\,
\big(
2\,w_{ir}\,w_{jr}
+
2\,w_{ri}\,w_{rj}
\big)\,
\partial_{w_{rr}}
\\
&
\ \ \ \ \
+
\sum_{i<r<j}\,
\big(
w_{ij}\,w_{rj}
-
w_{ri}\,w_{jj}
+
w_{ji}\,w_{jr}
\big)\,
\partial_{w_{rj}}
\\
&
\ \ \ \ \
+
\sum_{1\leqslant s<i}\,
\big(
-w_{si}\,w_{jj}
+
w_{sj}\,w_{ji}
+
w_{ij}\,w_{js}
\big)\,
\partial_{w_{js}}
\\
&
\ \ \ \ \ \
+
\big(
2\,w_{ij}\,w_{ji}
\big)\,
\partial_{w_{ji}}
\endaligned
\]
\[
\aligned
{}
&
+
\sum_{i<s<j}\,
\big(
w_{is}\,w_{jj}
+
w_{ij}\,w_{js}
+
w_{sj}\,w_{ji}
\big)\,
\partial_{w_{js}}
\\
&
\ \ \ \ \
+
\big(
2\,w_{ij}\,w_{jj}
\big)\,
\partial_{w_{jj}}
\\
&
\ \ \ \ \
+
\sum_{j<s\leqslant\ell}\,
\big(
w_{ij}\,w_{js}
-
w_{ji}\,w_{sj}
+
w_{jj}\,w_{si}
\big)\,
\partial_{w_{js}}
\\
&
\ \ \ \ \
+
\sum_{j<r\leqslant\ell}\,
\big(
-w_{ii}\,w_{jr}
+
w_{ij}\,w_{ri}
+
w_{ir}\,w_{ji}
\big)\,
\partial_{w_{ri}}
\\
&
\ \ \ \ \
+
\sum_{j<r\leqslant\ell}\,
\big(
w_{ij}\,w_{rj}
+
w_{ir}\,w_{jj}
-
w_{ji}\,w_{jr}
\big)\,
\partial_{w_{rj}}
\\
&
\ \ \ \ \
+
\sum_{j<r\leqslant\ell}\,
\big(
2\,w_{ir}\,w_{rj}
-
2\,w_{jr}\,w_{ri}
\big)\,
\partial_{w_{rr}}
\endaligned
\]
\[
\aligned
{}
&
+
\sum_{1\leqslant s<r<i}\,
\Big(
-w_{si}\,w_{jr}
+
w_{sj}\,w_{ir}
-
w_{ri}\,w_{js}
+
w_{rj}\,w_{is}
\Big)\,
\partial_{w_{rs}}
\\
&
\ \ \ \ \
+
\sum_{1\leqslant r<s<i}\,
\Big(
w_{ri}\,w_{sj}
-
w_{rj}\,w_{si}
-
w_{ir}\,w_{js}
+
w_{is}\,w_{jr}
\Big)\,
\partial_{w_{rs}}
\\
&
\ \ \ \ \
+
\sum_{1\leqslant r<i}\,
\sum_{i<s<j}\,
\Big(
w_{ri}\,w_{sj}
+
w_{rj}\,w_{is}
-
w_{ir}\,w_{js}
+
w_{si}\,w_{jr}
\Big)\,
\partial_{w_{rs}}
\\
&
\ \ \ \ \
+
\sum_{1\leqslant r<i}\,
\sum_{j<s\leqslant\ell}\,
\Big(
-w_{ri}\,w_{js}
+
w_{rj}\,w_{is}
-
w_{ir}\,w_{sj}
+
w_{jr}\,w_{si}
\Big)\,
\partial_{w_{rs}}
\endaligned
\]
\[
\aligned
{}
&
+
\sum_{i<r<j}\,
\sum_{1\leqslant s<i}\,
\Big(
-w_{si}\,w_{jr}
+
w_{sj}\,w_{ri}
+
w_{is}\,w_{rj}
+
w_{ir}\,w_{js}
\Big)\,
\partial_{w_{rs}}
\\
&
\ \ \ \ \
+
\sum_{i<s<r<j}\,
\Big(
w_{is}\,w_{jr}
+
w_{ir}\,w_{js}
+
w_{si}\,w_{rj}
+
w_{sj}\,w_{ri}
\Big)\,
\partial_{w_{rs}}
\\
&
\ \ \ \ \
+
\sum_{i<r<s<j}\,
\Big(
-w_{ir}\,w_{sj}
+
w_{is}\,w_{rj}
-
w_{ri}\,w_{js}
+
w_{si}\,w_{jr}
\Big)\,
\partial_{w_{rs}}
\\
&
\ \ \ \ \
+
\sum_{i<r<j}\,
\sum_{j<s\leqslant\ell}\,
\Big(
w_{ir}\,w_{js}
+
w_{is}\,w_{rj}
-
w_{ri}\,w_{sj}
+
w_{jr}\,w_{si}
\Big)\,
\partial_{w_{rs}}
\endaligned
\]
\[
\aligned
{}
&
+
\sum_{j<r\leqslant\ell}\,
\sum_{1\leqslant s<i}\,
\Big(
-w_{si}\,w_{rj}
+
w_{sj}\,w_{ri}
-
w_{is}\,w_{jr}
+
w_{ir}\,w_{js}
\Big)\,
\partial_{w_{rs}}
\\
&
\ \ \ \ \
+
\sum_{j<r\leqslant\ell}\,
\sum_{i<s<j}\,
\Big(
w_{is}\,w_{rj}
+
w_{ir}\,w_{js}
-
w_{si}\,w_{jr}
+
w_{sj}\,w_{ri}
\Big)\,
\partial_{w_{rs}}
\\
&
\ \ \ \ \
+
\sum_{j<s<r\leqslant\ell}\,
\Big(
w_{is}\,w_{rj}
+
w_{ir}\,w_{sj}
-
w_{js}\,w_{ri}
-
w_{jr}\,w_{si}
\Big)\,
\partial_{w_{rs}}
\\
&
\ \ \ \ \
+
\sum_{j<r<s\leqslant\ell}\,
\Big(
w_{ir}\,w_{js}
-
w_{is}\,w_{jr}
-
w_{ri}\,w_{sj}
+
w_{rj}\,w_{si}
\Big)\,
\partial_{w_{rs}}.
\endaligned
\]

For $1 \leqslant i \leqslant \ell$:
\[
\aligned
IL_{w_{ii}w_{ii}}
&
\,:=\,
\isqrt\,
\sum_{1\leqslant\ell'<i}\,
\sum_{n'=1}^n\,
z_{in'}\,
\big(
w_{\ell'i}
+
w_{i\ell'}
\big)\,
\partial_{z_{\ell'n'}}
\\
&
\ \ \ \ \ 
+
\isqrt\,
\sum_{n'=1}^n\,
z_{in'}\,w_{ii}\,
\partial_{z_{in'}}
\\
&
\ \ \ \ \ \ 
\isqrt\,
\sum_{i<\ell'\leqslant\ell}\,
\sum_{n'=1}^n\,
z_{in'}\,
\big(
w_{\ell'i}
-
w_{i\ell'}
\big)\,
\partial_{z_{\ell'n'}}
\endaligned
\]
\[
\aligned
\\
{}
&
+
\isqrt\,
\sum_{1\leqslant\ell'<i}\,
\sum_{m'=1}^m\,
u_{im'}\,
\big(
w_{\ell'i}
+
w_{i\ell'}
\big)\,
\partial_{u_{\ell'm'}}
\\
{}
&
+
\isqrt\,
\sum_{m'=1}^m\,
u_{im'}\,w_{ii}\,
\partial_{u_{im'}}
\\
&
{}
+
\isqrt\,
\sum_{i<\ell'\leqslant\ell}\,
\sum_{m'=1}^m\,
u_{im'}\,
\big(
w_{\ell'i}
-
w_{i\ell'}
\big)\,
\partial_{u_{\ell'm'}}
\endaligned
\]
\[
\aligned
\\
{}
&
+
\isqrt\,
\sum_{1\leqslant\ell'<i}\,
\sum_{m'=1}^m\,
v_{im'}\,
\big(
w_{\ell'i}
+
w_{i\ell'}
\big)\,
\partial_{v_{\ell'm'}}
\\
{}
&
+
\isqrt\,
\sum_{m'=1}^m\,
v_{im'}\,w_{ii}\,
\partial_{v_{im'}}
\\
&
{}
+
\isqrt\,
\sum_{i<\ell'\leqslant\ell}\,
\sum_{m'=1}^m\,
v_{im'}\,
\big(
w_{\ell'i}
-
w_{i\ell'}
\big)\,
\partial_{v_{\ell'm'}}
\endaligned
\]
\[
\aligned
{}
&
+
\isqrt\,
\sum_{1\leqslant r<i}\,
\big(
-w_{ri}^2
+
w_{ir}^2
\big)\,
\partial_{w_{rr}}
\\
{}
&
+
\isqrt\,w_{ii}\,w_{ii}\,
\partial_{w_{ii}}
\\
{}
&
+
\isqrt\,
\sum_{i<r\leqslant\ell}\,
\big(
w_{ri}^2
-
w_{ir}^2
\big)\,
\partial_{w_{rr}}
\\
{}
&
+
\isqrt\,
\sum_{1\leqslant s<i}\,
w_{is}\,w_{ii}\,
\partial_{w_{is}}
\\
{}
&
+
\isqrt\,
\sum_{i<s\leqslant\ell}\,
w_{is}\,w_{ii}\,
\partial_{w_{is}}
\\
{}
&
+
\isqrt\,
\sum_{1\leqslant r<i}\,
w_{ri}\,w_{ii}\,
\partial_{w_{ri}}
\\
{}
&
+
\isqrt\,
\sum_{i<r\leqslant\ell}\,
w_{ri}\,w_{ii}\,
\partial_{w_{ri}}
\endaligned
\]
\[
\aligned
{}
&
+
\isqrt\,
\sum_{1\leqslant r<s<i}\,
\big(
w_{ri}\,w_{is}
-
w_{si}\,w_{ir}
\big)\,
\partial_{w_{rs}}
\\
{}
&
+
\isqrt\,
\sum_{1\leqslant s<r<i}\,
\big(
-w_{si}\,w_{ri}
+
w_{is}\,w_{ir}
\big)\,
\partial_{w_{rs}}
\\
{}
&
+
\isqrt\,
\sum_{i<r<s\leqslant\ell}\,
\big(
-w_{ir}\,w_{si}
+
w_{is}\,w_{ri}
\big)\,
\partial_{w_{rs}}
\\
{}
&
+
\isqrt\,
\sum_{i<s<r\leqslant\ell}\,
\big(
-w_{is}\,w_{ir}
+
w_{si}\,w_{ri}
\big)\,
\partial_{w_{rs}}
\\
{}
&
+
\isqrt\,
\sum_{1\leqslant r<i}\,
\sum_{i<s\leqslant\ell}\,
\big(
w_{ri}\,w_{si}
+
w_{ir}\,w_{is}
\big)\,
\partial_{w_{rs}}
\\
{}
&
+
\sum_{i<r\leqslant\ell}\,
\sum_{1\leqslant s<i}\,
\big(
w_{si}\,w_{ir}
+
w_{is}\,w_{ri}
\big)\,
\partial_{w_{rs}}.
\endaligned
\]

For $1 \leqslant j < i \leqslant \ell$:
\[
\aligned
IL_{w_{ij}w_{ij}}
&
\,:=\,
\isqrt\,
\sum_{1\leqslant\ell'<j}\,
\sum_{n'=1}^n\,
\Big(
z_{in'}\,
\big[
w_{\ell'j}
+
w_{j\ell'}
\big]
+
z_{jn'}\,
\big[
w_{\ell'i}
+
w_{i\ell'}
\big]
\Big)\,
\partial_{z_{\ell'n'}}
\\
&
\ \ \ \ \
+
\isqrt\,
\sum_{n'=1}^n\,
\Big(
z_{in'}\,
\big[
w_{jj}
\big]
+
z_{jn'}\,
\big[
w_{ji}
+
w_{ij}
\big]
\Big)\,
\partial_{z_{jn'}}
\\
&
\ \ \ \ \
+
\isqrt\,
\sum_{j<\ell'<i}\,
\sum_{n'=1}^n\,
\Big(
z_{in'}\,
\big[
-w_{j\ell'}
+
w_{\ell'j}
\big]
+
z_{jn'}\,
\big[
w_{\ell'i}
+
w_{i\ell'}
\big]
\Big)\,
\partial_{z_{\ell'n'}}
\\
&
\ \ \ \ \
+
\isqrt\,
\sum_{n'=1}^n\,
\Big(
z_{in'}\,
\big[
w_{ij}
-
w_{ji}
\big]
+
z_{jn'}
\big[
w_{ii}
\big]
\Big)\,
\partial_{z_{in'}}
\\
&
\ \ \ \ \ 
+
\isqrt\,
\sum_{i<\ell'\leqslant\ell}\,
\sum_{n'=1}^n\,
\Big(
z_{in'}\,
\big[
-w_{j\ell'}
+
w_{\ell'j}
\big]
+
z_{jn'}\,
\big[
-w_{i\ell'}
+
w_{\ell'i}
\big]
\Big)\,
\partial_{z_{\ell'n'}}
\endaligned
\]
\[
\aligned
{}
&
+
\isqrt\,
\sum_{1\leqslant\ell'<j}\,
\sum_{m'=1}^m\,
\Big(
u_{in'}\,
\big[
w_{\ell'j}
+
w_{j\ell'}
\big]
+
u_{jm'}\,
\big[
w_{\ell'i}
+
w_{i\ell'}
\big]
\Big)\,
\partial_{u_{\ell'm'}}
\\
{}
&
+
\isqrt\,
\sum_{m'=1}^m\,
\Big(
u_{im'}\,
\big[
w_{jj}
\big]
+
u_{jm'}\,
\big[
w_{ji}
+
w_{ij}
\big]
\Big)\,
\partial_{u_{jm'}}
\\
{}
&
+
\isqrt\,
\sum_{j<\ell'<i}\,
\sum_{m'=1}^m\,
\Big(
u_{im'}\,
\big[
-w_{j\ell'}
+
w_{\ell'j}
\big]
+
u_{jm'}\,
\big[
w_{\ell'i}
+
w_{i\ell'}
\big]
\Big)\,
\partial_{u_{\ell'm'}}
\\
{}
&
+
\isqrt\,
\sum_{m'=1}^m\,
\Big(
u_{im'}\,
\big[
w_{ij}
-
w_{ji}
\big]
+
u_{jm'}
\big[
w_{ii}
\big]
\Big)\,
\partial_{u_{im'}}
\\
{}
& 
+
\isqrt\,
\sum_{i<\ell'\leqslant\ell}\,
\sum_{m'=1}^m\,
\Big(
u_{im'}\,
\big[
-w_{j\ell'}
+
w_{\ell'j}
\big]
+
u_{jm'}\,
\big[
-w_{i\ell'}
+
w_{\ell'i}
\big]
\Big)\,
\partial_{u_{\ell'm'}}
\endaligned
\]
\[
\aligned
{}
&
+
\isqrt\,
\sum_{1\leqslant\ell'<j}\,
\sum_{m'=1}^m\,
\Big(
v_{in'}\,
\big[
w_{\ell'j}
+
w_{j\ell'}
\big]
+
v_{jm'}\,
\big[
w_{\ell'i}
+
w_{i\ell'}
\big]
\Big)\,
\partial_{v_{\ell'm'}}
\\
{}
&
+
\isqrt\,
\sum_{m'=1}^m\,
\Big(
v_{im'}\,
\big[
w_{jj}
\big]
+
v_{jm'}\,
\big[
w_{ji}
+
w_{ij}
\big]
\Big)\,
\partial_{v_{jm'}}
\\
{}
&
+
\isqrt\,
\sum_{j<\ell'<i}\,
\sum_{m'=1}^m\,
\Big(
v_{im'}\,
\big[
-w_{j\ell'}
+
w_{\ell'j}
\big]
+
v_{jm'}\,
\big[
w_{\ell'i}
+
w_{i\ell'}
\big]
\Big)\,
\partial_{v_{\ell'm'}}
\\
{}
&
+
\isqrt\,
\sum_{m'=1}^m\,
\Big(
v_{im'}\,
\big[
w_{ij}
-
w_{ji}
\big]
+
v_{jm'}
\big[
w_{ii}
\big]
\Big)\,
\partial_{v_{im'}}
\\
{}
& 
+
\isqrt\,
\sum_{i<\ell'\leqslant\ell}\,
\sum_{m'=1}^m\,
\Big(
v_{im'}\,
\big[
-w_{j\ell'}
+
w_{\ell'j}
\big]
+
v_{jm'}\,
\big[
-w_{i\ell'}
+
w_{\ell'i}
\big]
\Big)\,
\partial_{v_{\ell'm'}}
\endaligned
\]
\[
\aligned
{}
&
+
\isqrt\,
\sum_{1\leqslant r<j}\,
\big(
-2\,w_{rj}\,w_{ri}
+
2\,w_{jr}\,w_{ir}
\big)\,
\partial_{w_{rr}}
\\
&
\ \ \ \ \
+
\isqrt\,
\sum_{1\leqslant r<j}\,
\big(
w_{rj}\,w_{ij}
+
w_{ri}\,w_{jj}
-
w_{jr}\,w_{ji}
\big)\,
\partial_{w_{rj}}
\\
&
\ \ \ \ \ 
+
\isqrt\,
\sum_{1\leqslant r<j}\,
\big(
w_{rj}\,w_{ii}
+
w_{ri}\,w_{ij}
+
w_{ji}\,w_{ir}
\big)\,
\partial_{w_{ri}}
\\
&
\ \ \ \ \
+
\isqrt\,
\sum_{1\leqslant s<j}\,
\big(
-w_{sj}\,w_{ji}
+
w_{js}\,w_{ij}
+
w_{jj}\,w_{is}
\big)\,
\partial_{w_{js}}
\\
&
\ \ \ \ \
+
\isqrt\,
\big(
2\,w_{jj}\,w_{ji}
\big)\,
\partial_{w_{jj}}
\\
&
\ \ \ \ \
+
\isqrt\,
\sum_{j<s<i}\,
\big(
-w_{jj}\,w_{si}
+
w_{js}\,w_{ij}
+
w_{ji}\,w_{sj}
\big)\,
\partial_{w_{js}}
\\
&
\ \ \ \ \
+
\isqrt\,
\big(
2\,w_{ij}\,w_{ji}
\big)\,
\partial_{w_{ji}}
\endaligned
\]
\[
\aligned
{}
&
+
\isqrt\,
\sum_{i<s\leqslant\ell}\,
\big(
w_{jj}\,w_{is}
+
w_{ji}\,w_{sj}
+
w_{js}\,w_{ij}
\big)\,
\partial_{w_{js}}
\\
&
\ \ \ \ \
+
\isqrt\,
\sum_{j<r<i}\,
\big(
w_{jj}\,w_{ir}
+
w_{jr}\,w_{ji}
+
w_{rj}\,w_{ij}
\big)\,
\partial_{w_{rj}}
\\
&
\ \ \ \ \
+
\isqrt\,
\sum_{j<r<i}\,
\big(
2\,w_{jr}\,w_{ri}
+
2\,w_{rj}\,w_{ir}
\big)\,
\partial_{w_{rr}}
\\
&
\ \ \ \ \
+
\isqrt\,
\sum_{j<r<i}\,
\big(
-w_{jr}\,w_{ii}
+
w_{ji}\,w_{ir}
+
w_{ri}\,w_{ij}
\big)\,
\partial_{w_{ri}}
\\
&
\ \ \ \ \
+
\isqrt\,
\sum_{1\leqslant s<j}\,
\big(
w_{si}\,w_{ji}
+
w_{js}\,w_{ii}
+
w_{is}\,w_{ij}
\big)\,
\partial_{w_{is}}
\\
&
\ \ \ \ \ \
+
\isqrt\,
\big(
w_{jj}\,w_{ii}
+
w_{ij}^2
+
w_{ji}^2
\big)\,
\partial_{w_{ij}}
\endaligned
\]
\[
\aligned
{}
&
+
\isqrt\,
\sum_{j<s<i}\,
\big(
w_{ji}\,w_{si}
+
w_{sj}\,w_{ii}
+
w_{ij}\,w_{is}
\big)\,
\partial_{w_{is}}
\\
&
\ \ \ \ \
+
\isqrt\,
\big(
2\,w_{ij}\,w_{ii}
\big)\,
\partial_{w_{ii}}
\\
&
\ \ \ \ \
+
\isqrt\,
\sum_{i<s\leqslant\ell}\,
\big(
-w_{ji}\,w_{si}
+
w_{js}\,w_{ii}
+
w_{ij}\,w_{is}
\big)\,
\partial_{w_{is}}
\\
&
\ \ \ \ \
+
\isqrt\,
\sum_{i<r\leqslant\ell}\,
\big(
w_{jj}\,w_{ri}
+
w_{ji}\,w_{jr}
+
w_{ij}\,w_{rj}
\big)\,
\partial_{w_{rj}}
\\
&
\ \ \ \ \
+
\isqrt\,
\sum_{i<r\leqslant\ell}\,
\big(
-w_{ji}\,w_{ir}
+
w_{ij}\,w_{ri}
+
w_{ii}\,w_{rj}
\big)\,
\partial_{w_{ri}}
\\
&
\ \ \ \ \
+
\isqrt\,
\sum_{i<r\leqslant\ell}\,
\big(
-2\,w_{jr}\,w_{ir}
+
2\,w_{rj}\,w_{ri}
\big)\,
\partial_{w_{rr}}
\endaligned
\]
\[
\aligned
{}
&
+
\isqrt\,
\sum_{1\leqslant s<r<j}\,
\Big(
-w_{sj}\,w_{ri}
-
w_{si}\,w_{rj}
+
w_{js}\,w_{ir}
+
w_{jr}\,w_{is}
\Big)\,
\partial_{w_{rs}}
\\
&
\ \ \ \ \
+
\isqrt\,
\sum_{1\leqslant r<s<j}\,
\Big(
w_{rj}\,w_{is}
+
w_{ri}\,w_{js}
-
w_{sj}\,w_{ir}
-
w_{si}\,w_{jr}
\Big)\,
\partial_{w_{rs}}
\\
&
\ \ \ \ \
+
\isqrt\,
\sum_{1\leqslant r<j}\,
\sum_{j<s<i}\,
\Big(
w_{rj}\,w_{is}
+
w_{ri}\,w_{sj}
-
w_{jr}\,w_{si}
+
w_{js}\,w_{ir}
\Big)\,
\partial_{w_{rs}}
\\
&
\ \ \ \ \
+
\isqrt\,
\sum_{1\leqslant r<j}\,
\sum_{i<s\leqslant\ell}\,
\Big(
w_{rj}\,w_{si}
+
w_{ri}\,w_{sj}
+
w_{is}\,w_{jr}
+
w_{js}\,w_{ir}
\Big)\,
\partial_{w_{rs}}
\endaligned
\]
\[
\aligned
{}
&
+
\isqrt\,
\sum_{j<r<i}\,
\sum_{1\leqslant s<j}\,
\Big(
-w_{sj}\,w_{ri}
+
w_{si}\,w_{jr}
+
w_{js}\,w_{ir}
+
w_{rj}\,w_{is}
\Big)\,
\partial_{w_{rs}}
\\
&
\ \ \ \ \
+
\isqrt\,
\sum_{j<s<r<i}\,
\Big(
w_{js}\,w_{ri}
+
w_{jr}\,w_{si}
+
w_{sj}\,w_{ir}
+
w_{rj}\,w_{is}
\Big)\,
\partial_{w_{rs}}
\\
&
\ \ \ \ \
+
\isqrt\,
\sum_{j<r<s<i}\,
\Big(
-w_{jr}\,w_{is}
+
w_{js}\,w_{ir}
-
w_{rj}\,w_{si}
+
w_{ri}\,w_{sj}
\Big)\,
\partial_{w_{rs}}
\\
&
\ \ \ \ \
+
\isqrt\,
\sum_{j<r<i}\,
\sum_{i<s\leqslant\ell}\,
\Big(
-w_{jr}\,w_{si}
+
w_{js}\,w_{ir}
+
w_{is}\,w_{rj}
+
w_{ri}\,w_{sj}
\Big)\,
\partial_{w_{rs}}
\endaligned
\]
\[
\aligned
{}
&
+
\isqrt\,
\sum_{i<r\leqslant\ell}\,
\sum_{1\leqslant s<j}\,
\Big(
w_{sj}\,w_{ir}
+
w_{si}\,w_{jr}
+
w_{ri}\,w_{js}
+
w_{is}\,w_{rj}
\Big)\,
\partial_{w_{rs}}
\\
&
\ \ \ \ \
+
\isqrt\,
\sum_{i<r\leqslant\ell}\,
\sum_{j<s<i}\,
\Big(
-w_{js}\,w_{ir}
+
w_{jr}\,w_{si}
+
w_{sj}\,w_{ri}
+
w_{is}\,w_{rj}
\Big)\,
\partial_{w_{rs}}
\\
&
\ \ \ \ \
+
\isqrt\,
\sum_{i<s<r\leqslant\ell}\,
\Big(
-w_{js}\,w_{ir}
-
w_{jr}\,w_{is}
+
w_{sj}\,w_{ri}
+
w_{si}\,w_{rj}
\Big)\,
\partial_{w_{rs}}
\\
&
\ \ \ \ \
+
\isqrt\,
\sum_{i<r<s\leqslant\ell}\,
\Big(
-w_{jr}\,w_{si}
+
w_{js}\,w_{ri}
-
w_{ir}\,w_{sj}
+
w_{is}\,w_{rj}
\Big)\,
\partial_{w_{rs}}.
\endaligned
\]



\end{document}